\documentclass[11pt]{article}

\usepackage[T1]{fontenc}
\usepackage{lmodern}
\usepackage{microtype}
\usepackage{amsmath,amssymb,amsthm}
\usepackage{mathrsfs}
\usepackage{enumitem}
\usepackage[hidelinks]{hyperref}
\usepackage{xcolor}

\theoremstyle{plain}
\newtheorem{theorem}{Theorem}[section]
\newtheorem{lemma}[theorem]{Lemma}
\newtheorem{proposition}[theorem]{Proposition}
\newtheorem{corollary}[theorem]{Corollary}
\newtheorem*{theorem*}{Theorem}

\theoremstyle{definition}
\newtheorem{definition}[theorem]{Definition}
\newtheorem{remark}[theorem]{Remark}
\newtheorem{notation}[theorem]{Notation}
\newtheorem{fact}[theorem]{Fact}
\newcounter{claim}
\renewcommand{\theclaim}{\arabic{claim}}

\newenvironment{claim}[1][]%
{\refstepcounter{claim}\par\medskip
	\noindent\textit{Claim~\theclaim.}%
	\if\relax\detokenize{#1}\relax\else\ \textit{(#1)}\fi\ \itshape}%
{\par\medskip}

\newcommand{\Pow}{\mathcal{P}}
\newcommand{\bbP}{\mathbb{P}}
\newcommand{\bbQ}{\mathbb{Q}}

\newcommand{\cF}{\mathscr{F}}
\newcommand{\cG}{\mathscr{G}}
\newcommand{\cN}{\mathcal{N}}

\DeclareMathOperator{\Fn}{Fn}
\DeclareMathOperator{\Add}{Add}
\DeclareMathOperator{\dom}{dom}

\DeclareMathOperator{\Fix}{Fix}
\DeclareMathOperator{\Stab}{Stab}
\DeclareMathOperator{\Sym}{Sym}
\DeclareMathOperator{\HS}{HS}
\DeclareMathOperator{\id}{id}
\DeclareMathOperator{\ran}{ran}
\DeclareMathOperator{\Aut}{Aut}

\DeclareMathOperator{\cf}{cf}

\newcommand{\ZF}{\mathsf{ZF}}
\newcommand{\ZFC}{\mathsf{ZFC}}
\newcommand{\GBC}{\mathsf{GBC}}
\newcommand{\ETR}{\mathsf{ETR}}
\newcommand{\DC}{\mathsf{DC}}
\newcommand{\AC}{\mathsf{AC}}
\newcommand{\WO}{\mathsf{WO}}
\newcommand{\PP}{\mathsf{PP}}
\newcommand{\SVC}{\mathsf{SVC}}

\newcommand{\Ord}{\mathrm{Ord}}
\newcommand{\bbR}{\mathbb{R}}
\newcommand{\bbS}{\mathbb{S}}
\newcommand{\cM}{\mathcal{M}}
\newcommand{\supp}{\operatorname{supp}}
\newcommand{\one}{\mathbf{1}}

\newcommand{\tc}{\operatorname{tc}}

\title{A countable-support symmetric approach to separating \texorpdfstring{$\PP$}{PP} from \texorpdfstring{$\AC$}{AC}}
\author{Frank Trevor Gilson}
\date{\today}

\begin{document}
	\maketitle
	
	\begin{abstract}
		We study a countable-support Cohen symmetric seed model and a proposed
		symmetric-iteration approach to separating the Partition Principle
		$\PP$ from the Axiom of Choice.  The previously claimed final construction
		is withdrawn.  Its package forcing targets right inverses, and therefore
		a localized splitting principle stronger than the ordinary localized
		Partition Principle.  For the fixed seed parameter $S=A^\omega$, that
		splitting principle together with $\SVC(S)$ implies $\AC$, contrary to
		the intended preservation of a non-well-orderable Cohen set.  Independently,
		the final non-well-orderability argument confuses stabilization of forcing
		names with an action inside one fixed generic extension, and the limit-stage
		$\SVC(S)$ argument relies on an invalid truncation lemma.
		The retained positive result is the Cohen symmetric seed
		\[
			\cN\models\ZF+\DC+\SVC(S)+\neg\AC,
			\qquad S=(A^\omega)^{\cN}.
		\]
		The package and iteration sections are preserved only as a record of the
		superseded approach.  No model of
		$\ZF+\DC+\PP+\AC_{\WO}+\neg\AC$ is claimed here.
	\end{abstract}

	\section*{Correction and current status}
	The former main theorem and all conclusions depending on the proposed
	class-length iteration are withdrawn.  There are three separate
	load-bearing problems.

	First, let $S=A^\omega$.  The canonical interleaving of sequences gives
	$S\sqcup S\cong S$, and hence a coding
	$\Pow(S)\times S\hookrightarrow\Pow(S)$.  If
	$\mathcal F\subseteq\Pow(S)\setminus\{\varnothing\}$, the membership
	projection
	\[
		\{\langle B,s\rangle:B\in\mathcal F\text{ and }s\in B\}
		\twoheadrightarrow\mathcal F
	\]
	can therefore be coded as a surjection between subsets of $\Pow(S)$.
	A section supplies a choice function for $\mathcal F$.  Thus
	$\PP^{\mathrm{split}}\!\restriction\Pow(S)$ makes $S$ well-orderable.
	Together with $\SVC(S)$, which presents every set as a surjective image of
	$S\times\eta$ for some ordinal $\eta$, this implies $\AC$.

	Second, stabilization of a forcing name by an automorphism does not by
	itself yield an automorphism of its evaluation in the same generic extension.
	The former proof that $A$ remains non-well-orderable in the proposed final
	model uses exactly this unsupported same-generic inference.

	Third, the former limit-stage truncation lemma is false for arbitrary
	forcing names: an outer witnessing condition lying below an earlier stage
	does not ensure that the corresponding subname has the same value after
	truncation.  Consequently the claimed limit-stage persistence of $\SVC(S)$
	is not proved.

	The Cohen seed is independent of these failures.  The proofs of $\DC$ and
	$\SVC(S)$ below have been replaced by, respectively, a standard preservation
	argument and an explicit orbit-presentation argument.  The Ryan--Smith
	localization theorem is also retained with a corrected double-slice proof.
	The consistency problem motivating the paper remains open.
	
	\tableofcontents
	
	\section{Introduction}
	
	The Partition Principle $\PP$ asserts that whenever there is a surjection $A\twoheadrightarrow B$, there is an injection $B\hookrightarrow A$.
	Over $\ZFC$ this is immediate, since every surjection admits a right inverse.
	The aim of this manuscript was to separate $\PP$ from the full Axiom of
	Choice by a symmetric-iteration construction.  The proposed construction
	does not establish that separation.
	
	\paragraph{Status of the former main theorem.}
	The assertion that the displayed class-length iteration yields a model of
	\[
		\ZF+\DC+\PP+\AC_{\WO}+\neg\AC
	\]
	is withdrawn.  Sections~\ref{subsec:Qf} and~\ref{sec:iterate-packages}
	are retained as historical descriptions of the failed approach, not as
	certified constructions.
	
	\paragraph{Metatheory.}
	All forcing and symmetry constructions are carried out over a fixed ground universe $V$.
	Set-sized stages are treated in the usual set-theoretic background.
	To define and analyze the $\Ord$-length iteration as a definable class forcing (in particular, to carry out the relevant class recursions),
	we work in a background class theory supporting the relevant recursions, namely $\GBC+\mathsf{ETR}$,
	using a class well-order of $V$ only where it is explicitly invoked for bookkeeping;
	see Section~\ref{sec:prelims} and
	Remark~\ref{rem:metatheory}.
	
	\begin{remark}[The axiom $\AC_{\WO}$]
		\label{rem:ACwo}
		Throughout this paper, $\AC_{\WO}$ denotes \emph{choice for well-orderable sets}:
		every well-orderable family of nonempty sets admits a choice function
		(equivalently: for every family $\langle X_i:i\in I\rangle$ of nonempty sets indexed by a well-orderable set $I$,
		there is a choice function $c$ with $c(i)\in X_i$ for all $i\in I$).
		Equivalently, for every ordinal $\lambda$, the principle $AC_\lambda$ holds.
		This is the form of $\AC_{\WO}$ used when applying Ryan--Smith localization (Theorem~\ref{thm:rs-localization}).
	\end{remark}
	
	\paragraph{Why the proposed route fails.}
	Let $\cN$ be the Cohen symmetric seed model over $\Add(\omega,\omega_1)$
	constructed in Section~\ref{sec:cohen-seed}, let
	$A=\{c_\alpha:\alpha<\omega_1\}$, and fix $S=(A^\omega)^{\cN}$.
	The proposed successor packages schedule sections for surjections between
	subsets of $\Pow(S)$ and right inverses for bounded presentations of
	$\AC_{\WO}$.
	
	The local forcing facts about an individual package remain elementary, but
	the intended use of those packages is not viable.  Scheduling sections for
	all relevant surjections forces
	$\PP^{\mathrm{split}}\!\restriction\Pow(S)$, not merely
	$\PP\!\restriction\Pow(S)$.  As shown in the correction above, the former
	principle together with $\SVC(S)$ implies $\AC$.  The proposed iteration
	therefore cannot simultaneously establish its splitting target and retain
	$A$ as a witness to $\neg\AC$.
	
	\paragraph{What is claimed.}
	The positive theorem retained here is the countable-support Cohen symmetric
	seed $\cN\models\ZF+\DC+\SVC(S)+\neg\AC$.  We also retain the
	reduction facts in Section~\ref{sec:reduction-blueprint}, with the corrected
	proof of Ryan--Smith localization.  No final model, Ordering Principle,
	Kinna--Wagner consequence, or other corollary formerly derived from the
	withdrawn main theorem is claimed.
	
	\paragraph{Organization.}
	Section~\ref{sec:prelims} fixes the metatheoretic conventions and explains how the present manuscript
	carries out the specific \(\Ord\)-length recursion.
	Section~\ref{sec:cohen-seed} constructs the seed model $\cN$ and establishes
	$\DC$, $\neg\AC$, and $\SVC(S)$.
	Section~\ref{sec:reduction-blueprint} recalls Ryan--Smith localization and isolates the parameter $\Pow(S)$ together with the localized principles.
	Subsections~\ref{subsec:Qf} and~\ref{subsec:ACWO-package} record the
	package forcings used in the unsuccessful approach.
	Section~\ref{sec:iterate-packages} is retained only as a superseded
	historical construction.
	
	\section{Preliminaries and metatheoretic conventions}
	\label{sec:prelims}
	
	We follow the forcing, automorphism, symmetric-extension, and bounded-stage
	countable-support symmetric-iteration
	based on Karagila's iteration scheme \cite{KaragilaIteratingSymmExt} and Karagila and Schilhan \cite{KaragilaSchilhan}.
	This section fixes the metatheory and notation used throughout, and explains how
	the present paper carries out the specific \(\Ord\)-length recursion needed for the
	package iteration.
	
	\subsection{Background theory and conventions}
	
	\begin{remark}[Metatheory vs.\ object theory]\label{rem:meta-pp}
		All forcing notions, names, automorphism groups, filters, and supports are constructed
		over a fixed background set universe $V$.
		\begin{itemize}
			\item For the set-forcing and symmetric-extension statements used at \emph{bounded}
			stages, it suffices to work over $V\models\ZF$.
			\item Whenever we invoke preservation of dependent choice for a \emph{bounded}
			stage symmetric model, we assume $V\models\ZFC$.
			Throughout this paper, $\DC$ means $\DC_\omega$ (dependent choice for
			$\omega$-sequences).
			\item For the class-length recursion through $\Ord$ (bookkeeping and the definition of the
			iteration itself), we work in a background \emph{two-sorted} structure $(V,\mathcal C)$
			satisfying a class theory adequate for transfinite recursion along $\Ord$
			(e.g.\ $(V,\mathcal C)\models\GBC+\ETR_{\Ord}$, or simply $\GBC+\ETR$ in the standard sense).
			When convenient we fix a class well-order of $V$ (equivalently, assume Global Choice in the
			background), but we use such Global Choice only for bookkeeping/canonical enumerations and
			explicitly indicate when it is invoked.
		\end{itemize}
		
		When we speak of a “$V$-generic” $G\subseteq\bbP_{\Ord}$, this is understood in the metatheory
		(equivalently, in an outer universe) as a subclass meeting every dense class $D\in\mathcal C$.
		We emphasize that we do \emph{not} rely on a global forcing theorem for $\bbP_{\Ord}$:
		valuations of \emph{set-sized} $\bbP_{\Ord}$-names localize to bounded forcing stages
		(Lemma~\ref{lem:Ord-stage-localization}), and set-sized final \emph{hereditarily symmetric}
		names descend to bounded-stage hereditarily symmetric names
		(Lemma~\ref{lem:Ord-HS-stage-descent}).
		Accordingly, set-theoretic arguments are carried out
		inside initial-segment extensions $V[G_\alpha]$ and their bounded-stage symmetric submodels.
		Finally, when we say “countable”, “$\omega_1$-complete”, etc., this is always with respect to the
		background universe $V$.
	\end{remark}
	
	\subsection{Restrictions, supports, and stage generics}
	
	\begin{remark}[Forcing vs.\ group restrictions]\label{rem:restriction-conventions}
		We follow standard notation discipline for forcing versus group restrictions and do
		\emph{not} overload restriction/projection notation.
		\begin{itemize}
			\item For forcing, $p\upharpoonright\beta$ denotes the usual restriction of a
			condition $p\in\bbP_\lambda$ to an initial segment $\beta\le\lambda$ (and when
			needed, $\pi_{\beta,\lambda}:\bbP_\lambda\to\bbP_\beta$ denotes the canonical
			projection map).
			\item For groups, $\rho_{\beta,\lambda}:\cG_\lambda\to\cG_\beta$ denotes the
			restriction homomorphism.
		\end{itemize}
	\end{remark}
	
	\begin{notation}[Stage generics]\label{not:stage-generics}
		Fix (in the metatheory) a $V$-generic subclass $G\subseteq\bbP_{\Ord}$ for the class-length iteration forcing.
		For each set ordinal $\alpha$, write $G_\alpha:=G\cap\bbP_\alpha$ for the induced $V$-generic filter on the
		initial-segment forcing.
	\end{notation}
	
	\subsection{Countable-support iteration facts proved internally or cited}
	\label{subsec:imported-facts}
	
	\begin{remark}[Countable-support iteration facts]\label{rem:iteration-api}
		We use the following countable-support iteration facts, proved internally or cited from \cite{KaragilaIteratingSymmExt} and \cite{KaragilaSchilhan}.
		Each time we invoke an item below, we will explicitly point to the corresponding hypothesis in the current construction.
		
		\begin{enumerate}[label=(\alph*)]
			\item \textbf{Bounded-stage recursion and coherence.}
			For each set ordinal \(\lambda\), the countable-support symmetric iteration provides the bounded-stage forcing \(\bbP_\lambda\), the group \(\cG_\lambda\), and the projection system $\pi^\gamma_\beta : \cG_\beta \to \cG_\gamma$ satisfying coherence. For limit \(\lambda\), $\cG_\lambda$ is the inverse limit (Definition~\ref{def:limit-stage-group}). Cf.\ \cite[Def.~4.4, Lemma~4.5]{KaragilaSchilhan} for the general $I$-support framework with $I = [\lambda]^{\le\omega}$.
			
			\item \textbf{Framework and actions.}
			At each stage \(\lambda\) the iteration provides a single ambient group
			\(\cG_\lambda\le\Aut(\bbP_\lambda)\) acting by automorphisms on the full initial-segment forcing
			\(\bbP_\lambda\); see \cite[\S 5]{KaragilaIteratingSymmExt}.
			
			\item \textbf{Limit filters.}
			The limit-stage filter used in this manuscript is the modified filter \(\tilde{\cF}^*_\lambda\)
			(Definition~\ref{def:modified-limit-filter}).
			Its normality and \(\omega_1\)-completeness are proved here (Lemma~\ref{lem:filter-properties}).
			
			\item \textbf{Bounded-stage $\DC$.}
			For any set ordinal \(\lambda\), assuming \(V\models\ZFC\), the stage filter $\tilde{\cF}^*_\lambda$ is $\omega_1$-complete (Lemma~\ref{lem:filter-properties}(iv)). The iterand forcings are countably closed (Lemma~\ref{lem:Qf-closed}), hence the stage symmetric model $M_\lambda \models \DC$ by Lemma~\ref{lem:countably-closed-symm-DC} (using Lemma~\ref{lem:filter-properties}(iv) for $\omega_{1}$-completeness of $\tilde{\mathcal{F}}^{*}_{\lambda}$).
		\end{enumerate}
		
		\smallskip
		We stress that external sources are used here only at bounded stages.
		In particular, this paper does \emph{not} import any general theorem asserting \(\DC\) or
		\(\ZF\) for the final class-length symmetric model.
		The only facts cited from \cite{KaragilaIteratingSymmExt,KaragilaSchilhan} are the projection/coherence infrastructure for the forcing/group recursion and the general framework.
		Bounded-stage $\DC$ follows from Lemma~\ref{lem:countably-closed-no-reals-DC}(2) and Lemma~\ref{lem:filter-properties}(iv) as noted above.
		For the present strengthened filter \(\tilde{\cF}^*_\lambda\), these hypotheses are verified
		in the present paper (notably Lemma~\ref{lem:filter-properties}).
		The final-model arguments are then proved directly here: \(\DC\) is obtained by combining item~(d) with the stage-descent lemma for set-sized final hereditarily symmetric names (Lemma~\ref{lem:Ord-HS-stage-descent}), while \(\ZF\) is verified by reducing each set-sized final-name argument to a common bounded stage via Lemmas~\ref{lem:Ord-stage-localization} and~\ref{lem:Ord-HS-stage-descent}.
	\end{remark}
	
	\medskip
	Choice principles and the Ryan--Smith localization theorem are stated (and used)
	in Section~\ref{sec:reduction-blueprint}.
	
	\section{The Cohen symmetric seed model}
	\label{sec:cohen-seed}
	
	We fix a ground model \(V\models\ZFC\). Let
	\[
	\bbP \;=\; \Add(\omega,\omega_1)\;:=\;\Fn(\omega_1\times\omega,2,{<}\omega),
	\]
	ordered by reverse inclusion.
	For each \(\alpha<\omega_1\), let \(\dot c_\alpha\) be the canonical \(\bbP\)-name for the \(\alpha\)-th Cohen real,
	\[
	\dot c_\alpha \;:=\;
	\{\langle \check n,p\rangle : p\in\bbP\ \wedge\ (\alpha,n)\in\dom(p)\ \wedge\ p(\alpha,n)=1\}.
	\]
	Fix \(G_0\subseteq\bbP\) \(V\)-generic and write \(c_\alpha:=\dot c_\alpha^{G_0}\in 2^\omega\).
	
	\begin{lemma}[Cohen reals are pairwise distinct]
		\label{lem:cohen-reals-distinct}
		For $\alpha \neq \beta < \omega_1$, the Cohen reals $c_\alpha$ and $c_\beta$ are distinct.
	\end{lemma}
	
	\begin{proof}
		Fix $\alpha\neq\beta<\omega_1$.
		Let
		\[
		D_{\alpha,\beta}=\Bigl\{p\in\bbP:\exists n\in\omega\ \bigl[(\alpha,n),(\beta,n)\in\dom(p)\ \wedge\ p(\alpha,n)\neq p(\beta,n)\bigr]\Bigr\}.
		\]
		Given any $p$, choose $n$ such that neither $(\alpha,n)$ nor $(\beta,n)$ is in $\dom(p)$,
		and extend $p$ to $q\le p$ with $q(\alpha,n)=0$ and $q(\beta,n)=1$.
		Then $q\in D_{\alpha,\beta}$.
		Hence \(D_{\alpha,\beta}\) is dense, so \(G_0\) meets it and therefore \(c_\alpha\neq c_\beta\).
	\end{proof}
	
	\subsection{The forcing, group action, and the countable-support filter}
	
	\begin{definition}[Automorphisms of \(\Add(\omega,\omega_1)\)]
		Let \(\cG:=\Sym(\omega_1)\) be the full permutation group of \(\omega_1\).
		Each \(\pi\in\cG\) induces an automorphism of \(\bbP\) by permuting the \(\omega_1\)-coordinate:
		\[
		\dom(\pi p)=\{(\pi(\alpha),n):(\alpha,n)\in\dom(p)\}
		\quad\text{and}\quad
		(\pi p)(\pi(\alpha),n)=p(\alpha,n).
		\]
		This extends to \(\bbP\)-names by the standard recursion:
		\[
		\pi\dot x\;:=\;\{\langle \pi\dot y,\pi p\rangle:\langle \dot y,p\rangle\in\dot x\}.
		\]
	\end{definition}
	
	\begin{definition}[Countable-support filter of subgroups]
		\label{def:cohen-filter}
		For \(E\subseteq \omega_1\), write
		\[
		\Fix(E)\;:=\;\{\pi\in\cG : \pi\!\restriction E=\mathrm{id}_E\}.
		\]
		Let \(\cF\) be the normal filter of subgroups of \(\cG\) generated by the family
		\(\{\Fix(E):E\in[\omega_1]^{\le\omega}\}\).
		A (normal) filter of subgroups is upward closed and closed under finite intersections (and conjugation).
	\end{definition}
	
	\begin{lemma}[Basis for the Cohen filter]\label{lem:cohen-filter-basis}
		If $H\in\cF$, then there is $E\in[\omega_1]^{\le\omega}$ such that $\Fix(E)\le H$.
	\end{lemma}
	
	\begin{proof}
		Unwinding “$\cF$ is the normal filter generated by $\{\Fix(E)\}$”, there are
		$E_0,\dots,E_{k-1}\in[\omega_1]^{\le\omega}$ and $\pi_0,\dots,\pi_{k-1}\in\cG$ such that
		\[
		\bigcap_{i<k}\pi_i\,\Fix(E_i)\,\pi_i^{-1}\ \le\ H.
		\]
		But $\pi_i\Fix(E_i)\pi_i^{-1}=\Fix(\pi_i[E_i])$, and finite intersections of $\Fix(\cdot)$’s are
		$\Fix(\bigcup \cdot)$.
		Let $E:=\bigcup_{i<k}\pi_i[E_i]$ (finite union of countable sets is countable). Then
		$\Fix(E)\le H$.
	\end{proof}
	
	\begin{corollary}[External $\omega_1$-completeness of $\cF$]\label{cor:cohen-filter-omega1complete}
		In the ground model $V$, the filter $\cF$ is $\omega_1$-complete.
	\end{corollary}
	
	\begin{proof}
		Let $\delta<\omega_1$ and let $\langle H_\xi:\xi<\delta\rangle\in V$
		with each $H_\xi\in\cF$.  By Lemma~\ref{lem:cohen-filter-basis},
		choose in $V$ countable $E_\xi\subseteq\omega_1$ with
		$\Fix(E_\xi)\le H_\xi$.  Since $V\models\ZFC$ and $\delta$ is
		countable in $V$, the set $E=\bigcup_{\xi<\delta}E_\xi$ is countable.
		Then
		$\Fix(E)\le\bigcap_{\xi<\delta}H_\xi$, so the intersection belongs
		to $\cF$.
	\end{proof}
	
	\begin{definition}[Symmetric and hereditarily symmetric names]
		For a \(\bbP\)-name \(\dot x\), define its stabilizer
		\[
		\Stab(\dot x)\;:=\;\{\pi\in\cG:\pi\dot x=\dot x\}.
		\]
		For \(E\subseteq\omega_1\), we say that \(E\) \emph{supports} \(\dot x\) if \(\Fix(E)\le\Stab(\dot x)\).
		We call \(\dot x\) \emph{symmetric} if \(\Stab(\dot x)\in\cF\), equivalently if \(\dot x\) has some
		\(V\)-countable support \(E\in[\omega_1]^{\le\omega}\) (by Lemma~\ref{lem:cohen-filter-basis}).
		Let \(\HS\) be the class of \emph{hereditarily symmetric} names: \(\dot x\in\HS\iff \dot x\) is symmetric and
		every name appearing in \(\dot x\) lies in \(\HS\).
	\end{definition}
	
	Fix \(G_0\subseteq \bbP\) \(V\)-generic. The associated symmetric extension is
	\[
	\cN \;:=\;\HS^{G_0} \;=\;\{\dot x^{G_0}:\dot x\in\HS\}.
	\]
	
	\subsection{CCC and dependent choice}
	
	\begin{lemma}[CCC]
		\label{lem:add-ccc}
		\(\Add(\omega,\omega_1)\) is ccc.
		In particular, it preserves all cardinals and cofinalities.
	\end{lemma}
	
	\begin{proof}
		Let $\{p_\xi:\xi<\omega_1\}\subseteq \Add(\omega,\omega_1)$ be uncountable.
		Write $D_\xi=\dom(p_\xi)\subseteq \omega_1\times\omega$, so each $D_\xi$ is finite.
		
		By thinning out, assume $|D_\xi|=n$ for all $\xi$.
		By the $\Delta$--system lemma (applied to the finite sets $D_\xi$), thin out again
		to an uncountable set $I\subseteq\omega_1$ such that $\{D_\xi:\xi\in I\}$ forms a
		$\Delta$--system with root $R$.
		There are only finitely many possible values for $p_\xi\!\restriction R$.
		Thin out once more so that $p_\xi\!\restriction R=p_\eta\!\restriction R$ for all
		$\xi,\eta\in I$.
		Now for $\xi\neq\eta$ in $I$, the conditions $p_\xi$ and $p_\eta$ agree on
		$D_\xi\cap D_\eta=R$, hence $p_\xi\cup p_\eta$ is a condition extending both.
		Thus $p_\xi$ and $p_\eta$ are compatible. Therefore there is no uncountable antichain,
		i.e.\ $\Add(\omega,\omega_1)$ is ccc.
	\end{proof}
	
	\begin{lemma}[ZF for symmetric extensions]\label{lem:HS-models-ZF}
		Let \(\langle\bbP,\cG,\cF\rangle\) be a symmetric system, i.e.\ \(\cG\le\Aut(\bbP)\) and
		\(\cF\) is a normal filter of subgroups of \(\cG\).
		If \(H\subseteq\bbP\) is \(V\)-generic, then
		\(\HS^H=\{\dot x^H:\dot x\in\HS\}\) is a transitive model of \(\ZF\) and
		\(V\subseteq \HS^H \subseteq V[H]\).
	\end{lemma}
	
	\begin{proof}
		This is the standard symmetric extension theorem; see \cite[Lemma~15.51]{JechSetTheory}.
	\end{proof}
	
	\begin{lemma}[HS witnesses from HS-sets]\label{lem:HS-witness-from-HSset}
		Let $\dot A\in\HS$ and let $\varphi(v,\vec{\dot b})$ be any formula with parameters $\vec{\dot b}\in\HS$.
		If $p\Vdash \exists v\in\dot A\,\varphi(v,\vec{\dot b})$, then there exist $q\le p$ and $\dot a\in\HS$
		such that $q\Vdash \dot a\in\dot A\wedge \varphi(\dot a,\vec{\dot b})$.
	\end{lemma}
	
	\begin{proof}
		Fix $s\le p$. Since $p\Vdash\exists v\in\dot A\,\varphi(v,\vec{\dot b})$, there are $t\le s$ and a name $\dot v$
		with $t\Vdash \dot v\in\dot A\wedge \varphi(\dot v,\vec{\dot b})$.
		From $t\Vdash \dot v\in\dot A$, by the definition of forcing for membership there is $u\le t$ and
		$\langle \dot a,r\rangle\in\dot A$ such that $u\le r$ and $u\Vdash \dot v=\dot a$.
		Then $u\Vdash \varphi(\dot a,\vec{\dot b})$ as well, and since $\dot A\in\HS$ every name appearing in $\dot A$
		(in particular $\dot a$) is hereditarily symmetric.
		Thus $u$ witnesses the conclusion.
	\end{proof}
	
	\begin{remark}[Metatheoretic countability vs.\ internal countability]
		\label{rem:meta-vs-internal-countable}
		Throughout the construction of HS names, ``countable'' means \emph{countable in the ground model $V$}
		(e.g.\ $E\in[\omega_1]^{\le\omega}$ is evaluated in $V$).
		This notion is used only to ensure that
		$\Fix(E)\in\cF$ (Definition~\ref{def:cohen-filter}), hence that the relevant names are symmetric.
		In particular, when we form supports such as $E=\bigcup_{n<\omega}E_n$ in the \emph{metatheory},
		we use that $V\models\ZFC$, so a countable union of countable sets (in $V$) is countable (in $V$).
		We do \emph{not} claim that these supports are countable in the symmetric model $\cN$.
	\end{remark}
	
	\begin{theorem}[\(\DC\) in the symmetric model]
		\label{thm:dc-in-cohen-symm}
		\(\cN\models \ZF+\DC\).
	\end{theorem}
	
	\begin{proof}
		By Lemma~\ref{lem:HS-models-ZF}, the symmetric extension
		$\cN=\HS^{G_0}$ is a transitive model of $\ZF$.
		For $p\in\bbP$, let
		\[
			\operatorname{coord}(p)
			:=\{\alpha<\omega_1:\exists n\ ((\alpha,n)\in\dom(p))\}.
		\]
		This set is finite, and $\Fix(\operatorname{coord}(p))$ fixes $p$;
		hence every condition is tenacious.  By
		Lemma~\ref{lem:add-ccc}, $\bbP$ is $\omega_1$-cc, and by
		Corollary~\ref{cor:cohen-filter-omega1complete}, the symmetry filter is
		externally $\omega_1$-complete.  The standard preservation theorem for
		dependent choice in symmetric extensions
		\cite[Lemma~3.3]{KaragilaPreservingDC} therefore applies and gives
		$\cN\models\DC$.
	\end{proof}
	
	\subsection{Countable supports for \(\HS\)-names}
	
	\begin{notation}[Coordinate restriction to $E$]\label{not:PE-restriction}
		Let \(E\subseteq\omega_1\).
		Set
		\[
		\bbP_E \;:=\;\Add(\omega,E)\;=\;\Fn(E\times\omega,2,{<}\omega),
		\]
		viewed as a (complete) subposet of \(\bbP=\Add(\omega,\omega_1)\) via the inclusion
		\(\bbP_E\hookrightarrow\bbP\).
		
		For \(p\in\bbP\), define the restriction of \(p\) to \(E\) by
		\[
		p\!\restriction E \;:=\; p\cap\bigl((E\times\omega)\times 2\bigr).
		\]
		Equivalently, \(\dom(p\!\restriction E)=\dom(p)\cap(E\times\omega)\) and
		\((p\!\restriction E)(\alpha,n)=p(\alpha,n)\) for \((\alpha,n)\in\dom(p)\cap(E\times\omega)\).
		Then \(p\!\restriction E\in\bbP_E\) and \(p\le p\!\restriction E\) (reverse inclusion order).
	\end{notation}
	
	\begin{remark}[Symmetry supports]\label{rem:symmetry-support}
		Recall that \(E\subseteq\omega_1\) \emph{supports} a \(\bbP\)-name \(\tau\) if \(\Fix(E)\le\Stab(\tau)\).
		We refer to such \(E\) as a \emph{symmetry support} to distinguish it from forcing supports.
	\end{remark}
	
	\begin{lemma}[Every \(\HS\)-name has a countable support]
		\label{lem:hs-countable-support}
		For every \(\dot{x}\in\HS\) there exists \(E\in[\omega_1]^{\le\omega}\) such that
		\(\Fix(E)\le\Stab(\dot{x})\).
	\end{lemma}
	
	\begin{proof}
		If \(\dot x\in\HS\), then in particular \(\dot x\) is \(\cF\)--symmetric, so
		\(\Stab(\dot x)\in\cF\).
		By Lemma~\ref{lem:cohen-filter-basis}, there exists a countable \(E\in[\omega_1]^{\le\omega}\) with
		\(\Fix(E)\le\Stab(\dot x)\).
		This is exactly the assertion that $\dot x$ has countable support.
	\end{proof}
	
	\begin{lemma}[$\bbP_E$-names are hereditarily symmetric]\label{lem:PE-names-HS}
		Let $E\in[\omega_1]^{\le\omega}$ and $\bbP_E=\Add(\omega,E)\le\bbP=\Add(\omega,\omega_1)$.
		If $\tau$ is a $\bbP_E$-name, view $\tau$ as a $\bbP$-name via the inclusion
		$\bbP_E\hookrightarrow\bbP$ (i.e.\ the same set-theoretic name, whose conditions all lie in $\bbP_E\subseteq\bbP$).
		Then $\tau\in\HS$, and $E$ is a support for $\tau$ (i.e.\ $\Fix(E)\le\Stab(\tau)$).
	\end{lemma}
	
	\begin{proof}
		Fix $\pi\in\Fix(E)$.
		Since $\pi$ fixes $E$ pointwise, it acts trivially on $\bbP_E$:
		for every $p\in\bbP_E$ we have $\pi p=p$.
		We show by induction on the rank of $\tau$
		that $\pi\tau=\tau$.
		Indeed, if $\langle\sigma,p\rangle\in\tau$, then $p\in\bbP_E$
		so $\pi p=p$, and by induction $\pi\sigma=\sigma$, hence
		$\langle\pi\sigma,\pi p\rangle=\langle\sigma,p\rangle\in\tau$; similarly every member of $\pi\tau$
		lies in $\tau$.
		Thus $\pi\tau=\tau$.
		
		Therefore $\Fix(E)\le\Stab(\tau)$. Since $E$ is countable, $\Fix(E)\in\cF$, so $\tau$ is symmetric.
		The same argument applies to every subname of $\tau$ (all of whose conditions also lie in $\bbP_E$),
		so $\tau$ is hereditarily symmetric, i.e.\ $\tau\in\HS$.
	\end{proof}
	
	Lemma~\ref{lem:support-localization}, Remark~\ref{rem:support-vs-syntactic}, and Lemma~\ref{lem:PE-evaluation-in-subextension} record standard support and subextension facts for the Cohen seed,
	and will not be invoked explicitly later;
	they are included to fix conventions and prevent common
	confusions about symmetry supports versus $\bbP_E$-names.
	
	\begin{lemma}[Support localization for $\Add(\omega,\omega_1)$]
		\label{lem:support-localization}
		Let \(E\subseteq\omega_1\) and let \(\tau\) be a \(\bbP\)-name with \(\Fix(E)\le\Stab(\tau)\).
		Let \(\varphi\) be any formula and let \(\vec\sigma\) be a tuple of \(\bbP\)-names such that
		\(\Fix(E)\le\Stab(\sigma_i)\) for each \(i\).
		Then for every \(p\in\bbP\),
		\[
		p\Vdash \varphi(\tau,\vec\sigma)
		\quad\Longleftrightarrow\quad
		p\!\restriction E\Vdash \varphi(\tau,\vec\sigma).
		\]
	\end{lemma}
	
	\begin{proof}
		The implication “\(\Leftarrow\)” holds since \(p\le p\!\restriction E\).
		For “\(\Rightarrow\)”, suppose \(p\Vdash \varphi(\tau,\vec\sigma)\) but
		\(p\!\restriction E\nVdash \varphi(\tau,\vec\sigma)\). Then there is \(q\le p\!\restriction E\)
		with \(q\Vdash \neg\varphi(\tau,\vec\sigma)\).
		Let
		\[
		\begin{split}
			&F_p:=\{\alpha\in\omega_1\setminus E:\exists n\in\omega\ ((\alpha,n)\in\dom(p))\},\\
			&F_q:=\{\alpha\in\omega_1\setminus E:\exists n\in\omega\ ((\alpha,n)\in\dom(q))\}.
		\end{split}
		\]
		Choose \(\pi\in\Fix(E)\) such that \(\pi[F_p]\cap F_q=\varnothing\) and \(\pi\) fixes every element of \(F_q\).
		Then \(\pi p\) is compatible with \(q\), so pick \(r\le \pi p,q\).
		Since \(\pi\in\Fix(E)\le\Stab(\tau)\) and \(\Fix(E)\le\Stab(\sigma_i)\) for each \(i\), we have
		\(\pi\tau=\tau\) and \(\pi\sigma_i=\sigma_i\) for all \(i\).
		By automorphism invariance of the forcing relation,
		\(\pi p\Vdash \varphi(\pi\tau,\pi\vec\sigma)\), i.e.\ \(\pi p\Vdash \varphi(\tau,\vec\sigma)\), hence \(r\Vdash \varphi(\tau,\vec\sigma)\).
		But \(r\le q\) and \(q\Vdash \neg\varphi(\tau,\vec\sigma)\), contradiction.
	\end{proof}
	
	\begin{remark}[Symmetry support vs.\ forcing support]\label{rem:support-vs-syntactic}
		The condition $\Fix(E)\le\Stab(\tau)$ (``$\tau$ has symmetry support $E$'')
		means that $\tau$ is invariant under permutations fixing $E$ pointwise.
		This \emph{does not} imply that $\tau$ is a $\bbP_E$--name, nor that $\tau^{G_0}$ lies in
		$V[G_E]$.
		Indeed, letting \(\dot A:=\{\langle \dot c_\alpha,\mathbf 1\rangle:\alpha<\omega_1\}\) be the
		canonical name for the Cohen set $A=\{c_\alpha:\alpha<\omega_1\}$, we have
		$\Stab(\dot A)=\Sym(\omega_1)$, so $\dot A$ has symmetry support $\emptyset$.
		Nevertheless, $\dot A$ is not a $\bbP_E$--name for any countable $E$, and $A\notin V[G_E]$
		for every countable $E\subseteq\omega_1$.
	\end{remark}
	
	\begin{lemma}[Subextension evaluation for $\bbP_E$-names]\label{lem:PE-evaluation-in-subextension}
		Let $E\subseteq\omega_1$, let $\bbP_E=\Add(\omega,E)\le\bbP=\Add(\omega,\omega_1)$,
		and let $\tau$ be a $\bbP_E$--name (viewed as a $\bbP$--name via the inclusion).
		If $G_0\subseteq\bbP$ is $V$--generic and $G_E:=G_0\cap\bbP_E$, then
		\[
		\tau^{G_0}=\tau^{G_E}\in V[G_E].
		\]
	\end{lemma}
	
	\begin{proof}
		By induction on the rank of $\tau$.
		Using the inclusion $\bbP_E\hookrightarrow\bbP$, the valuation recursion for $\tau^{G_0}$
		only consults conditions from $G_0$ that appear in $\tau$, and these conditions all lie in
		$\bbP_E$, hence belong to $G_E$.
		Applying the inductive hypothesis to subnames of $\tau$ yields $\tau^{G_0}=\tau^{G_E}$.
	\end{proof}
	
	\subsection{The canonical Cohen set \(A\) and failure of choice}
	
	Recall that for each \(\alpha<\omega_1\), \(\dot c_\alpha\) denotes the canonical \(\bbP\)-name for the
	\(\alpha\)-th Cohen real and \(c_\alpha:=\dot c_\alpha^{G_0}\in 2^\omega\).
	
	\begin{definition}[The Cohen set]\label{def:cohen-set-A}
		Define the canonical \(\bbP\)-name
		\[
		\dot A \;:=\; \{\langle \dot c_\alpha,\mathbf 1_{\bbP}\rangle:\alpha<\omega_1\},
		\]
		and set
		\[
		A \;:=\; \dot A^{G_0} \;=\; \{c_\alpha:\alpha<\omega_1\}\subseteq 2^\omega.
		\]
	\end{definition}
	
	\begin{lemma}\label{lem:A-in-N}
		\(A\in\cN\).
	\end{lemma}
	
	\begin{proof}
		By Lemma~\ref{lem:PE-names-HS}, each \(\dot c_\alpha\) lies in \(\HS\) (indeed it has symmetry support \(\{\alpha\}\)).
		For any \(\pi\in\cG\), \(\pi\dot A=\dot A\), hence \(\Stab(\dot A)=\cG\in\cF\) and \(\dot A\) is symmetric.
		Since all members of \(\dot A\) are hereditarily symmetric, \(\dot A\in\HS\).
		Therefore \(A=\dot A^{G_0}\in \HS^{G_0}=\cN\).
	\end{proof}
	
	\begin{proposition}[Non-well-orderability]
		\label{prop:A-not-wo}
		In \(\cN\), the set \(A\) is not well-orderable.  In particular, \(\cN\models \neg\WO\) and hence \(\neg\AC\).
	\end{proposition}
	
	\begin{proof}
		We show that in \(\cN\) there is no injection from \(A\) into any ordinal.
		Since
		\(\AC\) implies that every set is well-orderable (equivalently, injects into an
		ordinal), this will imply that \(A\) is not well-orderable in \(\cN\), hence
		\(\cN\models\neg\WO\) and therefore \(\cN\models\neg\AC\).
		
		Suppose toward a contradiction that there exist an ordinal \(\theta\), a name \(\dot f\in\HS\), and a condition \(p_0\in\bbP\) such that
		\[
		p_0\Vdash\text{``\(\dot f:\dot A\to\check\theta\) is an injection.''}
		\]
		By Lemma~\ref{lem:hs-countable-support} fix a countable support \(E\in[\omega_1]^{\le\omega}\)
		for \(\dot f\), i.e.\ \(\Fix(E)\le\Stab(\dot f)\).
		
		Let \(\supp(p)\subseteq\omega_1\) denote the finite set of ordinals \(\alpha\) such that
		\((\alpha,n)\in\dom(p)\) for some \(n\in\omega\).
		Choose \(\alpha\in\omega_1\setminus(E\cup\supp(p_0))\).
		Since \(p_0\) forces that \(\dot f(\dot c_\alpha)\) is an ordinal \(<\theta\), there are
		\(q\le p_0\) and \(\xi<\theta\) such that
		\[
		q\Vdash \dot f(\dot c_\alpha)=\check\xi.
		\]
		Now choose \(\beta\in\omega_1\setminus\bigl(E\cup\supp(q)\bigr)\) with \(\beta\neq\alpha\),
		and let \(\pi\in\Sym(\omega_1)\) be the transposition swapping \(\alpha\) and \(\beta\)
		and fixing every other ordinal.
		Then \(\pi\in\Fix(E)\), so \(\pi\dot f=\dot f\),
		and also \(\pi\dot c_\alpha=\dot c_\beta\).
		
		By the standard automorphism invariance of the forcing relation,
		\[
		q\Vdash \dot f(\dot c_\alpha)=\check\xi
		\quad\Longrightarrow\quad
		\pi q\Vdash (\pi\dot f)(\pi\dot c_\alpha)=\check\xi,
		\]
		so using \(\pi\dot f=\dot f\) and \(\pi\dot c_\alpha=\dot c_\beta\) we obtain
		\[
		\pi q\Vdash \dot f(\dot c_\beta)=\check\xi.
		\]
		Because \(q\) mentions \(\alpha\) but not \(\beta\), and \(\pi\) fixes every ordinal in
		\(\supp(q)\setminus\{\alpha\}\), the conditions \(q\) and \(\pi q\) are compatible;
		indeed,
		they have disjoint requirements on the \(\alpha\)-row versus the \(\beta\)-row, and agree
		on their common domain.
		Let \(r:=q\cup \pi q\), so \(r\in\bbP\) and \(r\le q,\pi q\).
		Then
		\[
		r\Vdash \dot f(\dot c_\alpha)=\check\xi \ \wedge\  \dot f(\dot c_\beta)=\check\xi.
		\]
		
		Finally, strengthen \(r\) to a condition \(r'\le r\) forcing \(\dot c_\alpha\neq\dot c_\beta\):
		choose \(n\in\omega\) with \((\alpha,n),(\beta,n)\notin\dom(r)\) and put
		\[
		r':=r\cup\{((\alpha,n),0),((\beta,n),1)\}.
		\]
		Then \(r'\Vdash \dot c_\alpha(n)\neq \dot c_\beta(n)\), hence \(r'\Vdash \dot c_\alpha\neq\dot c_\beta\).
		But \(r'\le r\) still forces \(\dot f(\dot c_\alpha)=\dot f(\dot c_\beta)\), contradicting the
		injectivity of \(\dot f\) on \(\dot A\).
		
		This contradiction shows that no injection \(A\to\theta\) exists in \(\cN\) for any ordinal \(\theta\).
		Therefore \(A\) is not well-orderable in \(\cN\), and consequently \(\cN\models\neg\WO\) and
		\(\cN\models\neg\AC\).
	\end{proof}
	
	\begin{remark}[Countable subsets of \(A\)]\label{rem:A-dedekind-infinite}
		In contrast with the finite-support generalized Cohen model (where the analogous \(A\) can be Dedekind-finite; cf.\ \cite[\S2.4]{ransomBPI}),
		the present countable-support filter admits many canonical countable subsets of \(A\).
		Fix \(E\in[\omega_1]^{\le\omega}\cap V\) and set
		\[
		\dot A_E \;:=\; \{\langle \dot c_\alpha,\mathbf 1_{\bbP}\rangle:\alpha\in E\}.
		\]
		Then \(\Fix(E)\le\Stab(\dot A_E)\), hence \(\dot A_E\in\HS\) and
		\[
		A_E:=\dot A_E^{G_0}=\{c_\alpha:\alpha\in E\}\in\cN.
		\]
		If \(E\) is infinite and \(b:\omega\to E\) is a bijection in \(V\), then the sequence
		\(\langle c_{b(n)}:n\in\omega\rangle\in\cN\) (Lemma~\ref{lem:code-countable-support-by-Aomega})
		witnesses that \(A\) is Dedekind-infinite in \(\cN\).
	\end{remark}
	
	\begin{lemma}[Coding a countable set of Cohen reals by an element of $A^\omega$]
		\label{lem:code-countable-support-by-Aomega}
		Let $E\in[\omega_1]^{\le\omega}$ be countable in $V$.
		Fix in $V$ a surjection
		$e:\omega\twoheadrightarrow E$. Then the sequence
		\[
		s_E:=\langle c_{e(n)}:n\in\omega\rangle
		\]
		is an element of $A^\omega$ in $\cN$.
	\end{lemma}
	
	\begin{proof}
		Each $\dot c_{e(n)}$ is hereditarily symmetric with support $\{e(n)\}$, hence
		the sequence name $\langle \dot c_{e(n)}:n\in\omega\rangle$ is hereditarily
		symmetric with support $E$.
		Therefore $s_E\in\cN$, and clearly $s_E\in A^\omega$.
	\end{proof}
	
	\begin{lemma}[Well-orders are rigid]
		\label{lem:wo-rigid}
		If $(X,\prec)$ is a well-order and $h:X\to X$ is a bijection such that
		\[
		x\prec y \iff h(x)\prec h(y)\qquad(x,y\in X),
		\]
		then $h=\mathrm{id}_X$.
	\end{lemma}
	
	\begin{proof}
		Suppose \(h\) is an order automorphism of \((X,\prec)\) and \(h\neq \mathrm{id}_X\).
		Let \(x_0\) be the \(\prec\)-least element of \(X\) with \(h(x_0)\neq x_0\).
		If \(h(x_0)\prec x_0\), then by minimality of \(x_0\) we have \(h(h(x_0))=h(x_0)\).
		But order preservation gives \(h(h(x_0))\prec h(x_0)\), contradiction.
		If \(x_0\prec h(x_0)\), apply the previous case to \(h^{-1}\), which is also an order automorphism of \((X,\prec)\),
		and satisfies \(h^{-1}(x_0)\prec x_0\), again a contradiction.
		Therefore \(h=\mathrm{id}_X\).
	\end{proof}
	
	\subsection{A safe $\SVC$ seed in $\cN$}
	
	\begin{definition}[\texorpdfstring{$\SVC(S)$}{SVC(S)}]\label{def:SVC}
		Let \(S\) be a set.
		We write \(\SVC(S)\) for the statement:
		\[
		\forall X\,\exists\eta\in\Ord\,\exists e\ (e:S\times\eta\twoheadrightarrow X).
		\]
	\end{definition}
	
	\begin{lemma}[Orbit presentation by $A^\omega$]
		\label{lem:orbit-presentation-Aomega}
		For every set $X\in\cN$ there are an ordinal $\eta$ and a partial
		function $g\in\cN$ such that
		\[
			\dom(g)\subseteq (A^\omega)^{\cN}\times\eta
			\qquad\text{and}\qquad
			\ran(g)=X.
		\]
	\end{lemma}
	
	\begin{proof}
		Fix $\dot X\in\HS$ with $X=\dot X^{G_0}$, and let
		$E_X\subseteq\omega_1$ be a countable support for $\dot X$.
		Put $H=\Fix(E_X)$.  Working in the ground model, enumerate the literal
		pairs of the set-name $\dot X$ as
		\[
			\dot X=\{\langle\dot\sigma_i,p_i\rangle:i<\eta\}
		\]
		for an ordinal $\eta$.  For each $i<\eta$, choose a countable support
		$E_i$ for $\dot\sigma_i$, put
		$D_i=E_i\cup\operatorname{coord}(p_i)\cup\{0\}$, and fix a
		surjection $d_i:\omega\twoheadrightarrow D_i$.
	
		For $i<\eta$ and $\pi\in H$, let
		\[
			\dot s_{i,\pi}
			=\langle\dot c_{\pi(d_i(n))}:n<\omega\rangle
		\]
		be the canonical sequence name.  It is hereditarily symmetric, with
		support $\pi[D_i]$.  Using canonical names for ordered pairs, form the
		name $\dot g$ whose graph contains, with condition $\pi p_i$, the entry
		mapping $(\dot s_{i,\pi},\check i)$ to $\pi\dot\sigma_i$, for every
		$i<\eta$ and $\pi\in H$.  This is a set-name in the ground model.
		If $\theta\in H$, its action sends the entry indexed by $(i,\pi)$ to
		the entry indexed by $(i,\theta\pi)$.  Hence $E_X$ supports $\dot g$,
		and all constituent names are hereditarily symmetric, so
		$\dot g\in\HS$.
	
		Let $g=\dot g^{G_0}$.  To see that $g$ is single-valued, suppose two
		active entries indexed by $(i,\pi)$ and $(j,\rho)$ have the same input.
		Then $i=j$, and equality of the two Cohen sequences, together with
		Lemma~\ref{lem:cohen-reals-distinct}, gives
		$\pi(d_i(n))=\rho(d_i(n))$ for every $n$.  Since $d_i$ maps onto
		$D_i$, the permutation $\rho^{-1}\pi$ fixes $E_i$ pointwise.
		Consequently
		$\pi\dot\sigma_i=\rho\dot\sigma_i$, so the two outputs agree.
	
		Every output lies in $X$: from
		$\langle\dot\sigma_i,p_i\rangle\in\dot X$ we have
		$p_i\Vdash\dot\sigma_i\in\dot X$, and automorphism invariance gives
		$\pi p_i\Vdash\pi\dot\sigma_i\in\dot X$ for $\pi\in H$.
		Conversely, if $x\in X$, some $i<\eta$ satisfies
		$p_i\in G_0$ and $x=\dot\sigma_i^{G_0}$.  The graph entry indexed by
		$(i,\mathrm{id})$ is then active and has output $x$.  Thus
		$\ran(g)=X$ and
		$\dom(g)\subseteq(A^\omega)^{\cN}\times\eta$.
	\end{proof}
	
	\begin{proposition}[$\SVC(A^\omega)$ in the Cohen symmetric model]
		\label{prop:svc-Aomega}
		Let $S=(A^\omega)^{\cN}$.  Then $\cN\models\SVC(S)$.
	\end{proposition}
	
	\begin{proof}
		Fix $X\in\cN$.  The empty set is immediate.  If $X\ne\varnothing$,
		apply Lemma~\ref{lem:orbit-presentation-Aomega} and choose one
		$x_0\in X$.  Extend the resulting partial surjection $g$ to
		$e:S\times\eta\twoheadrightarrow X$ by assigning the value $x_0$
		outside $\dom(g)$.  This uses one witness from one fixed nonempty set,
		not an indexed choice principle.
	\end{proof}
	
	\begin{definition}[The fixed seed parameter \(S\)]
		\label{def:fixed-ST}
		Work in \(\cN\) and let
		\[
			S:=(A^\omega)^{\cN}.
		\]
		Thereafter $S$ denotes this fixed seed set.  In a larger model it is not
		to be reinterpreted as that model's possibly larger set of
		$\omega$-sequences from $A$.
	\end{definition}
	
	\section{Reduction blueprint: global \texorpdfstring{$\PP$}{PP} from localized \texorpdfstring{$\PP$}{PP}}
	\label{sec:reduction-blueprint}
	
	We now pin down the exact hypotheses and the parameter set to which we will localize~$\PP$.
	
	\begin{definition}[Localized \texorpdfstring{$\PP$}{PP}]
		\label{def:PP-local}
		Let $T$ be a set.
		\begin{enumerate}
			\item $\PP\!\restriction T$ is the statement that for all $X,Y\subseteq T$, if there is a surjection
			$f:Y\twoheadrightarrow X$, then there is an injection $i:X\hookrightarrow Y$.
			\item $\PP(T)$ is the statement that for every set $X$, if there is a surjection $f:T\twoheadrightarrow X$,
			then there is an injection $i:X\hookrightarrow T$.
		\end{enumerate}
	\end{definition}
	
	\noindent
	The notation $\PP\!\restriction T$ follows \cite[\S 3.8]{SmithLocalReflections} and is unrelated to forcing-condition restrictions $p\restriction E$.
	
	\begin{definition}[\texorpdfstring{$\SVC^{+}$}{SVC+}]
		\label{def:SVCplus}
		For a set $T$, $\SVC^{+}(T)$ denotes the statement that for every set $X$ there exists an ordinal $\eta$
		and an injection $j:X\hookrightarrow T\times\eta$.
	\end{definition}
	
	\begin{fact}[Ryan--Smith]
		\label{fact:SVC-to-SVCplus}
		For every set $S$, $\SVC(S)$ implies $\SVC^{+}(\Pow(S))$.
		\cite[\S 2.1, Fact]{SmithLocalReflections}
	\end{fact}
	
	\begin{remark}[The fixed parameter \(S\) and model-relative powersets]
		\label{rem:fixed-T}
		By Fact~\ref{fact:SVC-to-SVCplus}, the seed model \(\cN\) satisfies
		\[
		\cN\models \SVC^+(\Pow(S)).
		\]
		However, whenever \(\Pow(S)\) appears below in a statement about a stage model \(M\) or about the final model \(\cM\), it means the powerset of the fixed set \(S\) as computed in that ambient model.
	\end{remark}
	
	\begin{theorem}[Ryan--Smith localization of \texorpdfstring{$\PP$}{PP}]
		\label{thm:rs-localization}
		Assume $\SVC^{+}(T)$ (Definition~\ref{def:SVCplus}).
		Then
		\[
		\PP \quad\Longleftrightarrow\quad \bigl(\PP\!\restriction T \ \wedge\  \AC_{\WO}\bigr).
		\]
		\cite[Proposition~3.17]{SmithLocalReflections} (stated there with parameter ``$S$''; 
		we apply it with $S := T$).
		Here $\PP\!\restriction T$ is as in Definition~\ref{def:PP-local}, and $\AC_{\WO}$ is as in Remark~\ref{rem:ACwo}.
		As this argument is critical to our goals, we demonstrate a proof here.
	\end{theorem}
	
	\begin{proof}
		($\Rightarrow$) If $\PP$ holds, then $\PP\!\restriction T$ is immediate.
		Moreover, $\PP$ implies $\AC_{\WO}$ (e.g.\ via the consequence $(\forall X)\,\aleph(X)=\aleph^\ast(X)$;
		see \cite[Proposition~3.17 and footnote~4]{SmithLocalReflections} for references).
		
		($\Leftarrow$) Assume $\PP\!\restriction T$ and $\AC_{\WO}$, and let $\varphi:Y\twoheadrightarrow X$ be a surjection.
		We must produce an injection $X\hookrightarrow Y$.
		Apply $\SVC^{+}(T)$ to the set
		\[
		Z \;:=\; (X\times\{0\})\ \cup\ (Y\times\{1\}),
		\]
		obtaining an ordinal $\eta$ and an injection $j:Z\hookrightarrow T\times\eta$.
		Define injections
		\[
		\iota_X:X\to T\times\eta,\quad \iota_X(x)=j(\langle x,0\rangle),
		\qquad
		\iota_Y:Y\to T\times\eta,\quad \iota_Y(y)=j(\langle y,1\rangle),
		\]
		and let
		\[
		A:=\ran(\iota_X)\subseteq T\times\eta,\qquad B:=\ran(\iota_Y)\subseteq T\times\eta.
		\]
		Then $\iota_X:X\to A$ and $\iota_Y:Y\to B$ are bijections, hence have inverses $\iota_X^{-1}:A\to X$ and
		$\iota_Y^{-1}:B\to Y$. Define
		\[
		f \;:=\;
		\iota_X\circ \varphi \circ \iota_Y^{-1} \;:\; B \twoheadrightarrow A.
		\]
		It therefore suffices to show that whenever $A,B\subseteq T\times\eta$ and $f:B\twoheadrightarrow A$,
		there is an injection $g:A\hookrightarrow B$.
		
		\begin{claim}[Slicing lemma for $T\times\eta$]
			\label{claim:slicing-lemma}
			Let $\eta$ be an ordinal. Assume $\PP\!\restriction T$ and $\AC_{\WO}$.
			If $A,B\subseteq T\times\eta$ and $f:B\twoheadrightarrow A$, then there is an injection $g:A\hookrightarrow B$.
		\end{claim}
		\begin{proof}[Proof of Claim~\ref{claim:slicing-lemma}]
			Treat $f$ as a partial map on $T\times\eta$ with domain $B$.
			For each $\langle t,\alpha\rangle\in A$, define
			\[
			\varepsilon_{t,\alpha}
			\;:=\;
			\min\bigl\{\varepsilon<\eta : (\exists s\in T)\ \bigl(\langle s,\varepsilon\rangle\in B\ \wedge\ f(\langle s,\varepsilon\rangle)=\langle t,\alpha\rangle\bigr)\bigr\}.
			\]
			(Nonempty since $f$ is surjective and $\min$ exists since $\eta$ is an ordinal.)
			
			For $\varepsilon,\alpha<\eta$, define the \emph{slices}
			\[
			A^\varepsilon_\alpha
			\;:=\;
			\{\,t\in T : \langle t,\alpha\rangle\in A\ \wedge\ \varepsilon_{t,\alpha}=\varepsilon\,\}\ \subseteq\ T,
			\]
			and
			\[
			B^\varepsilon_\alpha
			\;:=\;
			\{\,s\in T : \langle s,\varepsilon\rangle\in B\ \wedge\
			f(\langle s,\varepsilon\rangle)=\langle t,\alpha\rangle \text{ for some } t\in A^\varepsilon_\alpha\,\}\ \subseteq\ T.
			\]
			Define $h^\varepsilon_\alpha:B^\varepsilon_\alpha\to A^\varepsilon_\alpha$ by
			\[
			h^\varepsilon_\alpha(s)=t \quad\text{iff}\quad f(\langle s,\varepsilon\rangle)=\langle t,\alpha\rangle.
			\]
			This is well-defined (since $f$ is a function), and it is a surjection: if $t\in A^\varepsilon_\alpha$,
			then by definition of $\varepsilon_{t,\alpha}=\varepsilon$ there exists $s$ with
			$\langle s,\varepsilon\rangle\in B$ and $f(\langle s,\varepsilon\rangle)=\langle t,\alpha\rangle$, hence
			$s\in B^\varepsilon_\alpha$ and $h^\varepsilon_\alpha(s)=t$.
			
			By $\PP\!\restriction T$ applied to the surjection $h^\varepsilon_\alpha:B^\varepsilon_\alpha\twoheadrightarrow A^\varepsilon_\alpha$,
			there exists an injection $i^\varepsilon_\alpha:A^\varepsilon_\alpha\hookrightarrow B^\varepsilon_\alpha$ whenever
			$A^\varepsilon_\alpha\neq\varnothing$.
			
			Let
			\[
			E \;:=\;
			\{\,\langle \varepsilon,\alpha\rangle\in\eta\times\eta : A^\varepsilon_\alpha\neq\varnothing\,\}.
			\]
			Then $E$ is well-orderable (e.g.\ by the lexicographic order on $\eta\times\eta$), so by $\AC_{\WO}$ we may choose
			a function $c$ with
			\[
			c(\varepsilon,\alpha)\in
			\{\,i : i \text{ is an injection } A^\varepsilon_\alpha\hookrightarrow B^\varepsilon_\alpha\,\}
			\quad\text{for each }\langle \varepsilon,\alpha\rangle\in E.
			\]
			
			Define $g:A\to T\times\eta$ by
			\[
			g(\langle t,\alpha\rangle)
			\;:=\;
			\bigl\langle\, c(\varepsilon_{t,\alpha},\alpha)(t),\ \varepsilon_{t,\alpha}\,\bigr\rangle.
			\]
			Then $g(\langle t,\alpha\rangle)\in B$: indeed, $c(\varepsilon_{t,\alpha},\alpha)(t)\in B^{\varepsilon_{t,\alpha}}_\alpha$,
			so by definition of $B^{\varepsilon_{t,\alpha}}_\alpha$ we have
			$\bigl\langle c(\varepsilon_{t,\alpha},\alpha)(t),\varepsilon_{t,\alpha}\bigr\rangle\in B$.
			
			Finally, $g$ is injective. If
			$g(\langle t,\alpha\rangle)=g(\langle t',\alpha'\rangle)$, then
			$\varepsilon_{t,\alpha}=\varepsilon_{t',\alpha'}=:\varepsilon$ from equality of second coordinates. Let
			$s:=c(\varepsilon,\alpha)(t)=c(\varepsilon,\alpha')(t')$ be the common first coordinate.
			Then $s\in B^\varepsilon_\alpha\cap B^\varepsilon_{\alpha'}$. But if $\alpha\neq\alpha'$ this is impossible:
			$s\in B^\varepsilon_\alpha$ implies $f(\langle s,\varepsilon\rangle)$ has second coordinate $\alpha$, while
			$s\in B^\varepsilon_{\alpha'}$ implies the second coordinate is $\alpha'$, contradicting that $f$ is a function.
			Hence $\alpha=\alpha'$, and then $c(\varepsilon,\alpha)$ is injective, so $t=t'$.
			Thus $\langle t,\alpha\rangle=\langle t',\alpha'\rangle$.
		\end{proof}
		
		By Claim~\ref{claim:slicing-lemma} there is an injection $g:A\hookrightarrow B$.
		Then
		\[
		i \;:=\; \iota_Y^{-1}\circ g \circ \iota_X \;:\;
		X \to Y
		\]
		is an injection, as required. Therefore $\PP$ holds.
	\end{proof}
	
	\begin{remark}
		Ryan--Smith formulate $\PP$ as the equivalence $|X|\le|Y|\iff |X|\le^*|Y|$.
		Since $|X|\le|Y|\Rightarrow |X|\le^*|Y|$ holds in $\ZF$, this is equivalent to our
		surjection-to-injection formulation of $\PP$.
	\end{remark}
	
	\begin{remark}[On the hypotheses in Theorem~\ref{thm:rs-localization}]
		Theorem~\ref{thm:rs-localization} is exactly \cite[Proposition~3.17]{SmithLocalReflections}
		(with parameter $T$).
		We therefore keep both target principles \(\PP\!\restriction \Pow(S)\) and \(\AC_{\WO}\) explicit in the iteration.
		
		Ryan--Smith also considers the single-parameter instance $\PP(T)$ (Definition~\ref{def:PP-local}).
		In general, small-choice hypotheses together with $\PP(T)$ do \emph{not} suffice to recover global $\PP$:
		Cohen's first model satisfies $\SVC^{+}(\mathbb{R})\wedge \PP(\mathbb{R})\wedge\neg\PP$
		\cite[Proposition~3.22]{SmithLocalReflections}.
		Moreover, even adding $\AC_{\WO}$ does not allow one
		to replace a localized requirement of the form $\PP\!\restriction T$ by a single-parameter assertion
		$\PP(T)$: the Feferman-style model $N_{\aleph_1}$ of Truss satisfies
		$\AC_{\WO}\wedge \SVC^{+}(\Pow(\mathbb{R}))\wedge \PP(\Pow(\mathbb{R}))\wedge\neg\PP$
		\cite[Proposition~3.23]{SmithLocalReflections}.
		
		Finally, whether $\AC_{\WO}$ follows from $\SVC^{+}(T)\wedge \PP\!\restriction T$ is open in general;
		see \cite[Question~3.18]{SmithLocalReflections}.
		This is why $\AC_{\WO}$ remains an explicit forcing target below.
	\end{remark}
	
	\subsection{The localized $\PP$-package forcing $\bbQ_f$}
	\label{subsec:Qf}

	\paragraph{Status.}
	The individual forcing calculations in this subsection are retained, but
	the use of $\bbQ_f$ as a separation package is superseded.  Iterating these
	forcings schedules sections rather than merely reverse injections; together
	with the retained $\SVC(S)$ hypothesis, the resulting localized splitting
	principle implies $\AC$.
	
	Fix a transitive model \(M\models\ZF+\DC\) extending the Cohen seed \(\cN\)
	(in practice, \(M\) will be an intermediate stage of our symmetric iteration).
	Retain the fixed parameter \(S\) from Definition~\ref{def:fixed-ST}, and recall from
	Remark~\ref{rem:fixed-T} that \(\Pow(S)\) is interpreted model-relatively in the ambient model under discussion.
	
	\begin{definition}[Localized splitting principle]
		\label{def:PPsplitT}
		Let \(M\) be a transitive model containing the fixed parameter \(S\).
		Write
		\[
		\PP^{\mathrm{split}}\!\restriction \Pow(S)
		\]
		for the assertion that for all \(X,Y\subseteq \Pow(S)^M\) and all surjections
		\(f:Y\twoheadrightarrow X\) in \(M\), there exists a right inverse
		\(s:X\to Y\) with \(f\circ s=\id_X\) (necessarily injective).
	\end{definition}
	
	\begin{remark}[Why localized splitting does not solve the problem]
		\label{rem:localized-splitting}
		The global statement “every surjection splits” (i.e.\ every surjection has a right inverse)
		is equivalent to \(\AC\), so we are \emph{not} attempting to force global splittings.
	
		The proposed construction worked with the formally localized principle
		\[
		\PP^{\mathrm{split}}\!\restriction \Pow(S),
		\]
		where $S=(A^\omega)^{\cN}$ is the fixed seed parameter and
		$\Pow(S)$ is computed in the ambient model.  Although this concerns only
		surjections with domain and codomain contained in $\Pow(S)$, it is not
		compatible with the intended application: as explained in the correction
		note, $\SVC(S)$ together with this localized splitting principle implies
		$\AC$.
	
		The implication from localized splitting to
		$\PP\!\restriction\Pow(S)$ remains true, but it cannot be used here
		while preserving $\neg\AC$.
	
		\smallskip
		\noindent
		For the choice principles on well-orderable sets, Ryan--Smith isolate two related facts.
		First, under \(\SVC(S)\), the bounded global equivalence
		\[
		\AC_{\WO}\quad\Longleftrightarrow\quad \AC_{<\aleph^*(S)}(S)
		\]
		is given in \cite[Cor.~3.3]{SmithLocalReflections}.
		Second, the corresponding splitting formulation for fixed well-orderable targets is given in
		\cite[Prop.~3.2]{SmithLocalReflections}. This is the Ryan--Smith reduction used in
		Proposition~\ref{prop:final-ACWO}, together with the bounded scheduling horizon
		\(\lambda<\aleph^*(S)\).
	
		(For discussion of the hypotheses in Ryan--Smith's localization (including limitations on weakening
		them and the redundancy question for \(\AC_{\WO}\)), see
		\cite[Prop.~3.17, Question~3.18, and Props.~3.22--3.23]{SmithLocalReflections}.)
	\end{remark}
	
	\begin{definition}[$\PP$-package for a fixed surjection]
		\label{def:Qf}
		Let \(X,Y\subseteq \Pow(S)^M\) and let \(f:Y\twoheadrightarrow X\) be a surjection in \(M\).
		Define \(\bbQ_f\) to be the poset whose conditions are \emph{countable partial
		injective sections of \(f\)}:
		\[
		p\in\bbQ_f \iff
		\begin{cases}
			p\text{ is a function }p:\dom(p)\to Y,\\
			\dom(p)\subseteq X\text{ is countable},\\
			p\text{ is injective},\\
			\forall x\in\dom(p)\ \bigl(f(p(x))=x\bigr).
		\end{cases}
		\]
		The order is extension: \(q\leq p\) iff \(q\supseteq p\).
	\end{definition}
	
	\begin{remark}[Size and chain condition]
		\label{rem:Qf-size}
		The forcing $\mathbb Q_f$ can be large and is not expected to be ccc in general
		(e.g.\ $|\mathbb Q_f|$ may be comparable to $|X|^{\le\omega}\cdot|Y|^{\le\omega}$).
		This is harmless for our purposes: what we use is countable closure (hence no new reals and
		$\DC$ preservation), not any chain condition.
	\end{remark}
	
	\begin{lemma}[Basic density]
		\label{lem:Qf-dense}
		For each $x\in X$, the set
		\[
		D_x:=\{p\in\bbQ_f : x\in\dom(p)\}
		\]
		is dense in $\bbQ_f$.
	\end{lemma}
	
	\begin{proof}
		Fix $p\in\bbQ_f$ and $x\in X$ with $x\notin\dom(p)$.
		Let $y\in Y$ be such that $f(y)=x$ (exists since $f$ is surjective).
		
		We claim $y\notin\ran(p)$: if $y=p(x')$ for some $x'\in\dom(p)$ then
		$x=f(y)=f(p(x'))=x'$, contradicting $x\notin\dom(p)$.
		
		Thus $q:=p\cup\{(x,y)\}$ is still a partial function, still injective, and still satisfies
		$f(q(x))=x$; hence $q\in\bbQ_f$ and $q\le p$ with $x\in\dom(q)$.
	\end{proof}
	
	\begin{lemma}[Countable closure]
		\label{lem:Qf-closed}
		$\bbQ_f$ is countably closed (i.e.\ $\omega$-closed) in $M$.
	\end{lemma}
	
	\begin{proof}
		Let $\langle p_n:n<\omega\rangle$ be a decreasing sequence in $\bbQ_f$.
		Put $p:=\bigcup_{n<\omega}p_n$.
		Since each $p_{n+1}\supseteq p_n$, this union is a function.
		It is injective, and for every $x\in\dom(p)$ we have $f(p(x))=x$ because this holds already in
		some $p_n$ containing $x$.
		
		It remains to see that $\dom(p)=\bigcup_{n<\omega}\dom(p_n)$ is countable in $M$.
		Here we use $M\models\DC$, hence $M\models\AC_\omega$, so a countable union of countable sets
		is countable in $M$.
		(This is the only place in the closure argument where $\DC$ is used.)
		Therefore $p\in\bbQ_f$, and clearly $p\le p_n$ for all $n$.
	\end{proof}
	
	\begin{remark}[Closure depends on $\DC$]
		\label{rem:Qf-closure-DC}
		The countable closure of $\bbQ_f$ is proved inside $M$ using $\DC$ (equivalently $\AC_\omega$),
		via “countable union of countable sets is countable”.
		In a bare $\ZF$ context this need not hold,
		so Lemma~\ref{lem:Qf-closed} should be read as relative to our standing hypothesis $M\models\ZF+\DC$.
	\end{remark}
	
	\begin{proposition}[What $\bbQ_f$ forces]
		\label{prop:Qf-adds-section}
		If $G\subseteq\bbQ_f$ is $M$-generic, then
		\[
		s_G:=\bigcup G
		\]
		is a total injective function $s_G:X\to Y$ with $f\circ s_G=\id_X$.
		In particular, $M[G]\models ``X\hookrightarrow Y"$.
	\end{proposition}
	
	\begin{proof}
		By Lemma~\ref{lem:Qf-dense}, $G$ meets $D_x$ for every $x\in X$, so $x\in\dom(s_G)$ for all $x$
		and $s_G$ is total on $X$. Injectivity and $f(s_G(x))=x$ are preserved under unions of
		compatible conditions, hence hold for $s_G$.
	\end{proof}
	
	\begin{lemma}[Countably closed forcing adds no new $\omega$-sequences and preserves $\DC$]
		\label{lem:countably-closed-no-reals-DC}
		Let $M\models\ZF+\DC$ be transitive and let $\mathbb{U}\in M$ be countably closed in $M$
		(i.e.\ every decreasing $\omega$-sequence in $M$ has a lower bound in $M$).
		Then:
		\begin{enumerate}[label=(\arabic*)]
			\item $\mathbb{U}$ adds no new $\omega$-sequences of ordinals.
			\item Forcing with $\mathbb{U}$ over $M$ preserves $\DC$.
		\end{enumerate}
	\end{lemma}
	
	\begin{proof}
		\noindent\emph{(1) No new $\omega$-sequences of ordinals.}
		Let $u_0\in\mathbb{U}$ and $\dot a\in M$ with $u_0\Vdash \dot a:\check\omega\to\Ord$.
		For each $n<\omega$, let
		\[
		D_n:=\Bigl\{u\le u_0:\exists\alpha\in\Ord\ (u\Vdash \dot a(\check n)=\check\alpha)\Bigr\}.
		\]
		Each $D_n$ is dense below $u_0$.
		Define the set of \emph{states}
		\[
		\Sigma:=\Bigl\{(n,u,\vec\alpha): n<\omega,\ u\in\mathbb{U},\ u\le u_0,\ \vec\alpha\in\Ord^n,\
		u\Vdash \forall k<n\ \dot a(\check k)=\check{\vec\alpha(k)}\Bigr\}.
		\]
		
		Let $R$ be the relation on $\Sigma$ given by
		\[
		(n,u,\vec\alpha)\ R\ (n+1,v,\vec\beta)
		\]
		iff $v\le u$, $\vec\beta\in\Ord^{n+1}$ extends $\vec\alpha$, and
		$v\Vdash \dot a(\check n)=\check{\vec\beta(n)}$.
		Then $R$ is total on $\Sigma$ by density of $D_n$.
		By $\DC$ in $M$, there is an $R$-chain $\langle (n,u_n,\vec\alpha_n):n<\omega\rangle\in M$
		starting from some $(0,u_0,\emptyset)\in\Sigma$ (note $\Ord^0=\{\emptyset\}$).
		By countable closure, there is $u_\omega\in\mathbb{U}$ with $u_\omega\le u_n$ for all $n$.
		
		Let $g:=\bigcup_{n<\omega}\vec\alpha_n$.
		Since the chain is an element of $M$,
		Replacement in $M$ yields $g\in M$ and $g:\omega\to\Ord$.
		By construction, $u_\omega\Vdash \forall n\ \dot a(\check n)=\check{g(n)}$, hence
		$u_\omega\Vdash \dot a=\check g$. Thus $\mathbb{U}$ adds no new $\omega$-sequences of ordinals.
		
		\smallskip
		\noindent\emph{(2) Preservation of $\DC$.}
		Let $u_0\in\mathbb{U}$ and $\dot X,\dot R,\dot x_0\in M$ with
		\[
		u_0\Vdash \ \dot x_0\in\dot X\ \wedge\ (\forall x\in\dot X)(\exists y\in\dot X)\ x\,\dot R\,y.
		\]
		Define the state space
		\[
		\Gamma:=\Bigl\{(n,u,\dot x): n<\omega,\ u\in\mathbb{U},\ u\le u_0,\ \dot x\in M,\
		u\Vdash \dot x\in\dot X\Bigr\}.
		\]
		Define $S$ on $\Gamma$ by
		\[
		(n,u,\dot x)\ S\ (n+1,v,\dot y)
		\]
		iff $v\le u$ and $v\Vdash \dot y\in\dot X\ \wedge\ \dot x\,\dot R\,\dot y$.
		Totality of $S$ follows from the premise forced by $u_0$.
		By $\DC$ in $M$, pick an $S$-chain $\langle (n,u_n,\dot x_n):n<\omega\rangle\in M$
		starting with $(0,u_0,\dot x_0)\in\Gamma$.
		By countable closure, let $u_\omega\le u_n$ for all $n$.
		
		Let
		\[
		\dot f:=\bigl\{\langle \langle \check n,\dot x_n\rangle,\one_{\mathbb{U}}\rangle : n<\omega\bigr\}.
		\]
		Since the chain is in $M$, Replacement gives $\dot f\in M$.
		Then $u_\omega\Vdash \dot f:\check\omega\to\dot X$ and
		$u_\omega\Vdash \forall n\ \dot f(\check n)\,\dot R\,\dot f(\check{n+1})$,
		which is exactly $\DC$ in the extension.
	\end{proof}
	
	\begin{lemma}[Countably closed symmetric extensions preserve $\DC$]
		\label{lem:countably-closed-symm-DC}
		Let $M\models\ZF+\DC$ be transitive, let $\mathbb{P}\in M$ be countably closed in $M$,
		and let $\mathcal{F}$ be an $\omega_{1}$-complete normal filter on $\Aut(\mathbb{P})^{M}$.
		If $G\subseteq\mathbb{P}$ is $M$-generic, then the symmetric extension
		$\HS_{\mathcal{F}}^{G}$ satisfies $\DC$.
	\end{lemma}
	
	\begin{proof}
		Since $\mathbb{P}$ is countably closed in $M$, forcing with $\mathbb{P}$ over $M$
		preserves $\DC$ by Lemma~\ref{lem:countably-closed-no-reals-DC}(2); hence
		$M[G]\models\DC$.
		
		Let $\dot A,\dot R\in\HS_{\mathcal{F}}$ and $p\in G$ with
		$p\Vdash$``$\dot R$ is total on $\dot A$ and $\dot x_{0}\in\dot A$''.
		Work in $M[G]$ and put $A:=\dot A^{G}$, $R:=\dot R^{G}$, $x_{0}:=\dot x_{0}^{G}$.
		By $\DC$ in $M[G]$ choose $\langle a_{n}:n<\omega\rangle\in M[G]$ with
		$a_{0}=x_{0}$ and $a_{n}Ra_{n+1}$ for all $n$.
		
		Since $A\in\HS_{\mathcal{F}}^{G}$ and each $a_{n}\in A$, we have
		$a_{n}\in\HS_{\mathcal{F}}^{G}$; hence for each $n$ we may fix
		$\dot x_{n}\in\HS_{\mathcal{F}}$ with $\dot x_{n}^{G}=a_{n}$.
		Choose $H_{n}\in\mathcal{F}$ with $H_{n}\le\Stab(\dot x_{n})$.
		
		Put $H:=\bigcap_{n<\omega}H_{n}$. By $\omega_{1}$-completeness of $\mathcal{F}$,
		$H\in\mathcal{F}$. Define
		\[
		\dot s:=\Bigl\{\bigl\langle\langle\check n,\dot x_{n}\rangle,1_{\mathbb{P}}\bigr\rangle:n<\omega\Bigr\}.
		\]
		Since every element of $M[G]$ has a $\mathbb{P}$-name in $M$, we may choose each
		$\dot x_{n}\in M\cap\HS_{\mathcal{F}}$. The sequence $\langle\dot x_{n}:n<\omega\rangle$ is then an
		$\omega$-sequence of elements of $M$ in $M[G]$; by
		Lemma~\ref{lem:countably-closed-no-reals-DC}(1), it already belongs to $M$.
		Hence $\dot s\in M$.
		
		Since $H\le\Stab(\dot x_{n})$ for every $n$, we have $H\le\Stab(\dot s)$,
		hence $\dot s\in\HS_{\mathcal{F}}$. Moreover
		$1_{\mathbb{P}}\Vdash$``$\dot s:\check\omega\to\dot A$ and
		$\forall n\ \dot s(\check n)\mathbin{\dot R}\dot s(\check{n+1})$''.
		Thus $\dot s^{G}\in\HS_{\mathcal{F}}^{G}$ witnesses $\DC$ for $(A,R)$.
	\end{proof}
	
	\begin{corollary}[No new $\omega$-sequences; $\DC$ preservation]
		\label{cor:Qf-no-reals-DC}
		$\bbQ_f$ adds no new $\omega$-sequences of ordinals (in particular, no new reals).
		Moreover, over a $\ZF+\DC$ ground, forcing with $\bbQ_f$ preserves $\DC$.
	\end{corollary}
	
	\begin{proof}
		Apply Lemma~\ref{lem:countably-closed-no-reals-DC} with $\mathbb{U}=\bbQ_f$.
		Countable closure holds by Lemma~\ref{lem:Qf-closed}.
	\end{proof}
	
	\subsection{The $\AC_{\WO}$-package forcing $\mathbb{R}_f$}
	\label{subsec:ACWO-package}
	
	Fix a transitive model \(M\models\ZF+\DC\) extending the Cohen seed \(\cN\) (again,
	in practice \(M\) will be an intermediate stage of our iteration). Retain the
	fixed parameter \(S\) from Definition~\ref{def:fixed-ST}; as in Remark~\ref{rem:fixed-T},
	\(\Pow(S)\) will always be interpreted relative to the ambient model under discussion.
	
	\begin{definition}[The Lindenbaum bound $\aleph^*(S)$]
		\label{def:aleph-star}
		For a set $S$ in $M$, let $\aleph^*(S)^M$ denote the least ordinal $\kappa$ such that
		there is no surjection $S\twoheadrightarrow \kappa$ in $M$.
		Equivalently, $\kappa=\sup\{\alpha:\exists\text{ a surjection }S\twoheadrightarrow\alpha\}+1$.
	\end{definition}
	
	\begin{remark}
		In $\ZF$ the ordinal $\aleph^*(S)$ exists for every set $S$ (it is the Lindenbaum number of $S$);
		see e.g.\ Jech \cite[\S\textup{(Hartogs/Lindenbaum)}]{JechSetTheory}.
	\end{remark}
	
	\begin{remark}[Why the bound $\lambda<\aleph^*(S)$ is the correct target]
		\label{rem:ACwo-bound}
		Assuming $\SVC(S)$, Ryan--Smith shows that $\AC_{\WO}$ is equivalent to the bounded scheme
		$\AC_{<\aleph^*(S)}$ (more precisely, $\AC_{<\aleph^*(S)}(S)$ in his notation);
		see \cite[Cor.~3.3]{SmithLocalReflections}.  For fixed $\lambda$, the corresponding
		“splitting” reformulations of $\AC_\lambda$ appear in \cite[Prop.~3.2]{SmithLocalReflections}.
		Thus, in the presence of $\SVC(S)$ it suffices to schedule splittings only for surjections onto
		ordinals $\lambda<\aleph^*(S)^M$.
		In our application $\SVC(S)$ holds in each stage model; see
		Lemma~\ref{lem:svc-persists} below.
	\end{remark}
	
	\begin{definition}[$\AC_{\WO}$-package for a fixed surjection onto an ordinal]
		\label{def:Rf}
		Let $\lambda$ be an ordinal with $\lambda<\aleph^*(S)^M$, let $Y$ be a set in $M$,
		and let $f:Y\twoheadrightarrow \lambda$ be a surjection in $M$.
		Define $\mathbb{R}_f$ to be the poset of \emph{countable partial right inverses} to $f$:
		\[
		r\in\mathbb{R}_f \iff
		\begin{cases}
			r\text{ is a function }r:\dom(r)\to Y,\\
			\dom(r)\subseteq \lambda\text{ is countable},\\
			\forall \xi\in\dom(r)\ \bigl(f(r(\xi))=\xi\bigr).
		\end{cases}
		\]
		The order is extension: $s\le r$ iff $s\supseteq r$.
	\end{definition}
	
	\begin{remark}[Total sections are automatically injective]
		\label{rem:Rf-injective}
		If $s:\lambda\to Y$ satisfies $f\circ s=\id_\lambda$, then $s$ is injective:
		if $s(\xi_1)=s(\xi_2)$ then
		$\xi_1=f(s(\xi_1))=f(s(\xi_2))=\xi_2$.
	\end{remark}
	
	\begin{remark}[No localization to \(\Pow(S)\) for \(\AC_{\WO}\)]
		\label{rem:Rf-not-localized}
		Unlike the localized \(\PP\)-packages \(\bbQ_f\) (whose definitions quantify only over
		\(X,Y\subseteq \Pow(S)^M\)), the \(\AC_{\WO}\)-package is necessarily global:
		\(\AC_{\WO}\) concerns \emph{all} well-ordered index sets.
		The reduction in Remark~\ref{rem:ACwo-bound} shows that, assuming \(\SVC(S)\), it suffices
		to schedule splittings only for surjections onto ordinals \(\lambda<\aleph^*(S)^M\),
		which is what makes the bookkeeping manageable.
	\end{remark}
	
	\begin{lemma}[Basic density]
		\label{lem:Rf-dense}
		For each $\xi<\lambda$, the set
		\[
		D_\xi:=\{r\in\mathbb{R}_f : \xi\in\dom(r)\}
		\]
		is dense in $\mathbb{R}_f$.
	\end{lemma}
	
	\begin{proof}
		Fix $r\in\mathbb{R}_f$ and $\xi<\lambda$ with $\xi\notin\dom(r)$.
		Choose any $y\in Y$ with $f(y)=\xi$ (exists since $f$ is surjective).
		Then $s:=r\cup\{(\xi,y)\}$ is still a function, still satisfies $f(s(\xi))=\xi$,
		and has countable domain; hence $s\in\mathbb{R}_f$ and $s\le r$ with $\xi\in\dom(s)$.
	\end{proof}
	
	\begin{lemma}[Countable closure]
		\label{lem:Rf-closed}
		Assume $M\models\ZF+\DC$. Then $\mathbb{R}_f$ is countably closed.
	\end{lemma}
	
	\begin{proof}
		Let $\langle r_n:n<\omega\rangle$ be a decreasing sequence in $\mathbb{R}_f$.
		Put $r:=\bigcup_{n<\omega}r_n$.
		This is a function, and it satisfies
		$f(r(\xi))=\xi$ for every $\xi\in\dom(r)$ since this already holds in some $r_n$
		containing $\xi$.
		
		Finally, $\dom(r)=\bigcup_{n<\omega}\dom(r_n)$ is a countable union of countable
		sets. Since $M\models\DC$ (hence $M\models AC_\omega$), $\dom(r)$ is countable in $M$.
		Therefore $r\in\mathbb{R}_f$, and clearly $r\le r_n$ for all $n$.
	\end{proof}
	
	\begin{remark}[Closure depends on $\DC$]
		\label{rem:Rf-closure-DC}
		The countable closure of $\mathbb{R}_f$ (Lemma~\ref{lem:Rf-closed}) is proved inside $M$
		using $\DC$ (equivalently $\AC_\omega$), via “countable union of countable sets is countable”.
		In bare $\ZF$ this need not hold.
	\end{remark}
	
	\begin{proposition}[What $\mathbb{R}_f$ forces: a right inverse]
		\label{prop:Rf-adds-section}
		If $G\subseteq\mathbb{R}_f$ is $M$-generic, then
		\[
		s_G:=\bigcup G
		\]
		is a total function $s_G:\lambda\to Y$ with $f\circ s_G=\id_\lambda$
		(hence $s_G$ is injective).
		In particular, \(M[G]\models\) “\(f\) has a right inverse.”
	\end{proposition}
	
	\begin{proof}
		By Lemma~\ref{lem:Rf-dense}, $G$ meets $D_\xi$ for every $\xi<\lambda$, so
		$\xi\in\dom(s_G)$ for all $\xi<\lambda$ and $s_G$ is total on $\lambda$.
		The equation $f(s_G(\xi))=\xi$ is preserved under unions of compatible conditions,
		so $f\circ s_G=\id_\lambda$.
	\end{proof}
	
	\begin{corollary}[No new reals;
		DC preservation (package-level)]
		\label{cor:Rf-no-reals-DC}
		$\mathbb{R}_f$ adds no new reals (indeed, no new $\omega$-sequences of ordinals).
		Moreover, forcing with $\mathbb{R}_f$ over a $\ZF+\DC$ ground preserves $\DC$.
	\end{corollary}
	
	\begin{proof}
		Apply Lemma~\ref{lem:countably-closed-no-reals-DC} with $\mathbb{U}=\mathbb{R}_f$.
		Countable closure holds by Lemma~\ref{lem:Rf-closed}.
	\end{proof}
	
	\begin{remark}[Package products]
		\label{rem:package-product}
		When we say ``finite product of package posets'' we mean an arbitrary finite-support
		product of finitely many factors (not necessarily just two), typically consisting of
		one $\AC_{\WO}$-package $\mathbb{R}_f$ together with finitely many localized
		$\PP$-packages $\bbQ_{g_0},\dots,\bbQ_{g_{n-1}}$
		scheduled at the same stage.
	\end{remark}
	
	\section{The superseded package iteration over $\cN$}
	\label{sec:iterate-packages}

	\noindent
	\fbox{\parbox{0.94\textwidth}{
	\textbf{Status of this section.}
	This section records the former package-iteration proposal for historical
	clarity.  It is not an established construction.  Its claimed final
	conclusions, including preservation of $\DC$ and $\neg\AC$, persistence
	of $\SVC(S)$, and the final $\PP$ theorem, are withdrawn.  No statement
	in this section should be used as a premise without an independent proof.
	}}
	\medskip
	
	\begin{remark}[Roadmap of Section~\ref{sec:iterate-packages}]
		\label{rem:sec5-roadmap}
		This section describes the former class-length countable-support
		\emph{symmetric iteration} that was intended to force local splitting
		principles while preserving $\DC$ and maintaining $\neg\AC$.
		Here is the dependency flow.
		
		\begin{enumerate}[label=\textbf{(\arabic*)}]
			\item \textbf{Coordinates and products at successor stages.}
			At each successor stage $\alpha+1$ the forcing factor is a \emph{finite product} of orbit packages,
			and package coordinates are tracked by the triple system of
			Definition~\ref{def:coords}.
			This is the bookkeeping substrate for all later “fix a coordinate”
			arguments.
			
			\item \textbf{Diagonal infrastructure.}
			Subsection~\ref{subsec:iteration-diag} defines iteration-level diagonal lifts $\widehat{\pi}^{\,D}$
			and the associated subgroups $\Delta_\lambda^\uparrow(E,D)$, where $D\subseteq\mathrm{Coords}_{<\lambda}$ is
			\emph{countable} and $(E,D)$ is required to be \emph{globally admissible}
			(Definition~\ref{def:global-admissible}).
			The successor “diagonal kernel” elements
			(Definition~\ref{def:succ-diag-kernel} and Lemma~\ref{lem:succ-diag-kernel-in-kernel})
			supply the cancellation needed to define coherent lifts and to verify group membership.
			
			\item \textbf{The modified limit filter.}
			Subsection~\ref{subsec:limit-filter} defines the limit-stage filter $\tilde{\cF}^*_\lambda$
			(Definition~\ref{def:modified-limit-filter}) generated by pushforwards together with
			\emph{globally admissible} $\Delta^\uparrow$-subgroups.
			The key structural fact is the
			\emph{core lemma} (Lemma~\ref{lem:core-lemma}): every $H\in\tilde{\cF}^*_\lambda$ contains some
			$\Delta_\lambda^\uparrow(E,D)$ with globally admissible parameters.
			Conjugation is handled by
			Lemma~\ref{lem:delta-up-general-conj} together with the admissible-hull step
			(Lemma~\ref{lem:admissible-hull}), so normal-filter generation never “escapes” admissibility.
			
			\item \textbf{Bookkeeping by names (not interpretations).}
			Subsection~\ref{subsec:bookkeeping} fixes the class bookkeeping $\mathcal B$ of codes for potential
			$\bbP_\alpha$-names (Definition~\ref{def:bookkeeping-class}).
			At stage $\alpha$ we test the
			\emph{current} interpretation in $M_\alpha$ and add the corresponding orbit package factor(s) only
			when a code currently yields a relevant surjection instance not yet split.
			
			\item \textbf{The \(\Ord\)-length recursion.}
			Subsection~\ref{subsec:iteration-api-level} defines the specific \(\Ord\)-length iteration
			\(\langle \bbP_\alpha,\cG_\alpha,\tilde{\cF}^*_\alpha:\alpha\in\Ord\rangle\) used in this paper,
			starting from \(\bbP_0=\Add(\omega,\omega_1)\) and setting
			\(\bbP_{\alpha+1}=\bbP_\alpha*\dot{\bbS}_\alpha\) at successor stages, where
			\(\dot{\bbS}_\alpha\) is the finite product of the scheduled orbit-package factors.
			The bounded-stage successor/limit infrastructure comes from Remark~\ref{rem:iteration-api},
			while the transfinite recursion through \(\Ord\) is carried out explicitly in the
			background metatheory from Remark~\ref{rem:meta-pp}.
			By Lemma~\ref{lem:stage-iterand-HS}, the canonical \(\bbP_\alpha\)-names for these
			orbit-package iterands are hereditarily symmetric, so the two-step presentation is
			well-defined in the stage ground \(M_\alpha\).
			
		\item \textbf{Formerly claimed outcome.}
			Subsection~\ref{subsec:iteration-yields} records the withdrawn claims about the proposed final model \(\cM\):
			\(\ZF\) and \(\DC\) are preserved along the iteration; all surjections \(Y\twoheadrightarrow X\) with
			\(X,Y\subseteq \Pow(S)\) split (hence \(\PP^{\mathrm{split}}\!\restriction \Pow(S)\) in the ambient model under discussion);
			\(\SVC(S)\) persists; the bounded Ryan--Smith \(S\)-presentation of \(\AC_{\WO}\) is forced via scheduled right inverses, and hence \(\AC_{\WO}\) holds; the Ryan--Smith localization theorem upgrades
			\((\PP\!\restriction \Pow(S)\ \wedge\ \AC_{\WO})\) to full \(\PP\); and finally \(\neg\AC\) holds because the Cohen set \(A\) is not well-orderable (Proposition~\ref{prop:final-notAC}) using the core lemma plus a nontrivial diagonal lift acting on \(A\).
		\end{enumerate}
		
		\smallskip
		The material below is retained to identify the superseded mechanism and
		the exact locations of the failed preservation arguments.
	\end{remark}
	
	We now record the proposed symmetric iteration whose successor packages target the two local principles
	\(\PP^{\mathrm{split}}\!\restriction \Pow(S)\) and \(\AC_{\WO}\) over the Cohen symmetric seed model \(\cN\),
	while maintaining \(\neg\AC\) (witnessed by the non-well-orderability of \(A\)).
	All technical bounded-stage iteration infrastructure is imported via the API of
	Remark~\ref{rem:iteration-api}.
	The \(\Ord\)-length iteration used in this section is
	then defined explicitly in the present manuscript by recursion in the background
	metatheory from Remark~\ref{rem:meta-pp}.
	The purpose of retaining this section is to isolate the \emph{package
	systems} and bookkeeping targets that led to the withdrawn claim.
	
	\begin{remark}[Modularity of the construction]
		\label{rem:modularity}
		The argument is deliberately modular: the seed model technology (Section~\ref{sec:cohen-seed}),
		the localized \(\PP\)-packages (Subsection~\ref{subsec:Qf}), the \(\AC_{\WO}\)-packages
		(Subsection~\ref{subsec:ACWO-package}), and the symmetric-iteration ``API''
		(Remark~\ref{rem:iteration-api}) are separated so that if a hypothesis later needs to be
		strengthened or a package definition adjusted, it can typically be done locally without
		rewriting unrelated parts of the paper.
	\end{remark}
	
	\subsection{Package symmetric systems at a stage}
	\label{subsec:stage-systems}
	
	Fix a transitive \(M\models\ZF+\DC\) extending \(\cN\) (in practice \(M=M_\alpha\) at an
	intermediate stage of the iteration). Keep the fixed parameter \(S\) from
	Definition~\ref{def:fixed-ST}. Whenever \(\Pow(S)\) appears below, it is interpreted in the
	ambient model under discussion (Remark~\ref{rem:fixed-T}).
	
	\begin{remark}[An ambient \(\cG\)-action on \(S\) and \(\Pow(S)\)]
		\label{rem:ambient-action-T}
		Recall \(\cG=\Sym(\omega_1)\) (Definition~\ref{def:cohen-filter}). Using the Cohen reals
		\(c_\alpha\in A\), we obtain an induced action of \(\cG\) on \(A\) by
		\(\pi\cdot c_\alpha:=c_{\pi(\alpha)}\), hence on \(S=A^\omega\) by coordinatewise action, and on
		\(\Pow(S)\) by \(\pi\cdot X:=\{\pi\cdot s:s\in X\}\).
	
		We regard this as an action fixing all ordinals pointwise. More generally, whenever \(x\) is a set
		with \(\tc(x)\subseteq V(\Pow(S))\cup\Ord\), we extend the action recursively by
		\(\pi\cdot x:=\{\pi\cdot y:y\in x\}\). In particular, for a function \(f:X\to Y\) with
		\(\tc(f)\subseteq V(\Pow(S))\cup\Ord\) (e.g.\ \(X,Y\subseteq \Pow(S)\times\eta\) for an ordinal \(\eta\)), we
		identify \(f\) with its graph and set
		\[
		\pi\cdot f:=\{(\pi\cdot a,\pi\cdot b):(a,b)\in f\}.
		\]
	\end{remark}
	
	\begin{lemma}[Kernel acts trivially on \(\Pow(S)\)]
		\label{lem:ker-trivial-on-T}
		Let \(\alpha\) be any stage and let \(k\in\ker(\pi^0_\alpha)\), where \(\pi^0_\alpha\) is the stage-\(0\)
		projection of Definition~\ref{def:stage0-projection}.
		Then \(k\) fixes \(A\), \(S=A^\omega\), and every element of \(\Pow(S)\) pointwise for the ambient
		\(\cG\)-action of Remark~\ref{rem:ambient-action-T}.
		Consequently, \(k\) fixes every set/function whose transitive closure is contained in \(V(\Pow(S))\cup\Ord\);
		in particular, if \(f:X\to Y\) with \(X,Y\subseteq \Pow(S)\), or
		\(f:S\times\eta\twoheadrightarrow\lambda\) with \(\eta,\lambda\in\Ord\), then \(k\cdot f=f\).
	\end{lemma}
	
	\begin{proof}
		Put \(\sigma:=\pi^0_\alpha(k)\in\cG=\Sym(\omega_1)\). Since \(k\in\ker(\pi^0_\alpha)\), we have \(\sigma=\id\).
		By Remark~\ref{rem:ambient-action-T}, the ambient action on \(A\) and on \(S=A^\omega\) is induced from the action
		\(\alpha\mapsto\sigma(\alpha)\) on indices of the Cohen reals, so \(\sigma=\id\) implies that \(k\) fixes \(A\) and hence \(S\) pointwise.
		It follows immediately that every subset of \(S\) is fixed pointwise, hence so is every element of \(\Pow(S)\).
	
		For the final assertion, use the recursive extension of the action to \(V(\Pow(S))\) from
		Remark~\ref{rem:ambient-action-T} and argue by induction on rank that every \(x\) with
		\(\tc(x)\subseteq V(\Pow(S))\cup\Ord\) is fixed pointwise (using that \(k\) fixes all ordinals pointwise).
		In particular, any such function \(f\) is fixed.
	\end{proof}
	
	\begin{definition}[Orbit-symmetrized package systems]
		\label{def:package-systems}
		Let \(M\models\ZF+\DC\) be a transitive stage model and let \(\cG=\Sym(\omega_1)\) act on
		\(\Pow(S)\) as in Remark~\ref{rem:ambient-action-T}. Fix a surjection instance \(f\) in \(M\).
	
		Write \([f]:=\{\pi\cdot f:\pi\in\cG\}\) for the \(\cG\)-orbit of \(f\) (a set in \(M\) by Replacement).
	
		\begin{enumerate}[label=(\alph*)]
			\item \textbf{Orbit \(\PP^{\mathrm{split}}\!\restriction \Pow(S)\)-package.}
			If \(f:Y\twoheadrightarrow X\) with \(X,Y\subseteq \Pow(S)^M\), define
			\[
			\bbQ_{[f]}:=\prod^{\mathrm{cs}}_{g\in[f]}\bbQ_g,
			\]
			where \(\bbQ_g\) is the section forcing from Definition~\ref{def:Qf}.
			A condition \(\vec p\in\bbQ_{[f]}\) is a function with countable support
			\(\supp(\vec p):=\dom(\vec p)\subseteq[f]\) such that \(\vec p(g)\in\bbQ_g\) for all \(g\in\dom(\vec p)\),
			ordered by coordinatewise extension.
	
			\item \textbf{Orbit \(\AC_{\WO}\)-package (in \(S\)-presentation).}
			If \(\eta\) is an ordinal, \(\lambda<\aleph^*(S)^M\), and \(f:S\times\eta\twoheadrightarrow\lambda\) is a surjection in \(M\),
			define
			\[
			\bbR_{[f]}:=\prod^{\mathrm{cs}}_{g\in[f]}\bbR_g,
			\]
			where \(\bbR_g\) is the partial right-inverse forcing from Definition~\ref{def:Rf}.
		\end{enumerate}
	
		In both cases \(\cG\) acts on the orbit package by permuting the factors together with the
		\emph{canonical} transport induced by \(\pi\) on the defining parameters.
		More precisely, for \(\pi\in\cG\) and \(g\in[f]\) let
		\[
		A_{\pi,g}:\bbS_g\to \bbS_{\pi\cdot g}
		\]
		be the transport-of-structure isomorphism induced by \(\pi\) (i.e.\ apply \(\pi\) to every
		ordinal and every element of \(\Pow(S)\) appearing in a condition, using the ambient \(\cG\)-action
		on \(\Pow(S)\) from Remark~\ref{rem:ambient-action-T} and the fact that \(\bbS_g\) is defined
		functorially from \(g\)).
		Then the orbit action on conditions is
		\[
		(\pi\cdot \vec p)(g)\ :=\ A_{\pi,\pi^{-1}\cdot g}\bigl(\vec p(\pi^{-1}\cdot g)\bigr).
		\]
	\end{definition}
	
	\begin{lemma}[Orbit packages are countably closed and preserve \(\DC\)]
		\label{lem:orbit-package-closed}
		Work in \(M\models\ZF+\DC\).
		Let \(f\in M\) be a surjection instance and let \(\bbS_{[f]}\) be its
		orbit package. Then:
		\begin{enumerate}
			\item \(\bbS_{[f]}\) is countably closed;
			\item forcing with \(\bbS_{[f]}\) adds no new \(\omega\)-sequences of ordinals;
			\item \(M^{\bbS_{[f]}}\models\DC\).
		\end{enumerate}
	\end{lemma}
	\begin{proof}
		Let \(\langle \vec p_n:n<\omega\rangle\) be a descending sequence in \(\bbS_{[f]}\).
		Let \(C:=\bigcup_{n<\omega}\supp(\vec p_n)\subseteq [f]\).
		Each \(\supp(\vec p_n)\) is countable and
		\(M\models\DC\), so \(C\) is countable.
		
		For each \(g\in C\), consider the descending sequence \(\langle \vec p_n(g):n<\omega\rangle\) in the
		factor \(\bbS_g\) (interpreting \(\vec p_n(g)\) as the top condition whenever
		\(g\notin\supp(\vec p_n)\)).
		By Lemma~\ref{lem:Qf-closed} and Lemma~\ref{lem:Rf-closed}, each factor
		\(\bbS_g\) is countably closed, so choose \(q_g\in\bbS_g\) with \(q_g\le \vec p_n(g)\) for all \(n\).
		Define \(\vec q\) by \(\supp(\vec q)=C\) and \(\vec q(g)=q_g\).  Then \(\vec q\in\bbS_{[f]}\) and
		\(\vec q\le \vec p_n\) for all \(n\), proving countable closure.
		
		Items (2) and (3) follow from Lemma~\ref{lem:countably-closed-no-reals-DC}.
	\end{proof}
	
	\begin{remark}[Limit-stage \(\omega_1\)-completeness]
		The countable completeness of the diagonal stage filters is what makes the bounded-stage
		countable-support limit construction compatible with the hereditary-symmetry and
		\(\DC\)-preservation arguments.
		\(\omega_1\)-completeness is a key
		property of the stage filters at limits (see the limit-filter construction and the theorem
		on limit filter properties), and we reflect that requirement here as item~(c) in
		Remark~\ref{rem:iteration-api}.
	\end{remark}
	
	\begin{lemma}[Well-definedness of the orbit package actions]
		\label{lem:package-actions-wd}
		In each case of Definition~\ref{def:package-systems}, if $\vec p$ is a condition of the orbit package
		$\bbS_{[f]}$ (i.e.\ $\bbQ_{[f]}$ or $\bbR_{[f]}$) and $\pi\in\Sym(\omega_1)$, then $\pi\cdot \vec p$ is again a condition
		of the same forcing, and the action respects the order.
	\end{lemma}
	
	\begin{proof}
		By construction, each $A_{\pi,g}$ is a poset isomorphism
		$\bbS_g\cong \bbS_{\pi\cdot g}$, and the transport maps satisfy the cocycle law
		\[
		A_{\pi\circ\rho,g} \ =\ A_{\pi,\rho\cdot g}\circ A_{\rho,g}
		\]
		for all $\pi,\rho\in\Sym(\omega_1)$ and $g\in[f]$.
		Therefore $(\pi\circ\rho)\cdot\vec p=\pi\cdot(\rho\cdot\vec p)$ and
		$\id\cdot\vec p=\vec p$, so this is a well-defined action by poset
		automorphisms.
		Since the action only permutes coordinates and applies
		coordinatewise isomorphisms, it preserves the ideal of supports.
	\end{proof}
	
	\begin{lemma}[Countable \(\cG\)-support of successor coordinates]
		\label{lem:coord-countable-support}
		Let \(\alpha\) be a stage and let \(d\in\mathrm{Coords}_{\alpha+1}\setminus\mathrm{Coords}_\alpha\) be a new package
		coordinate at the successor stage \(\alpha+1\) (Definition~\ref{def:coords}).
		Then there exists a
		countable \(E_d\subseteq\omega_1\) such that \(\Fix(E_d)\le\Stab_{\cG}(d)\).
	\end{lemma}
	\begin{proof}
		Write \(d=\langle \alpha,i,g\rangle\) with \(i<m_\alpha\) and \(g\in[f_{\alpha,i}]\).
		By \(\mathbf{Support}_\alpha\) from Remark~\ref{rem:simultaneous-induction},
		applied at the predecessor stage \(\alpha\) to the canonical \(\HS_\alpha\)-name \(\dot f\) for \(f\),
		there is a countable globally admissible pair \((E_f,D_f)\) such that
		\[
		\Delta^\uparrow_\alpha(E_f,D_f)\le \Stab_{\cG_\alpha}(\dot f).
		\]
		Fix \(\pi\in\Fix(E_f)\).
		By Lemma~\ref{lem:diag-lift-in-group}, the diagonal lift \(\widehat\pi^{\,D_f}\in\cG_\alpha\) is defined, and by the
		choice of \((E_f,D_f)\) it fixes \(\dot f\).
		
		Let \(\widehat\pi^{\,\varnothing}\in\cG_\alpha\) denote the standard lift.
		Then
		\[
		k:=\bigl(\widehat\pi^{\,\varnothing}\bigr)^{-1}\cdot \widehat\pi^{\,D_f}\in\ker(\pi^0_\alpha),
		\]
		so by Lemma~\ref{lem:ker-trivial-on-T} the element \(k\) fixes \(f\) (since \(\tc(f)\subseteq V(\Pow(S))\cup\Ord\)).
		As \(\widehat\pi^{\,D_f}\) fixes \(f\), it follows that \(\widehat\pi^{\,\varnothing}\) fixes \(f\) as well, i.e.
		\(\pi\cdot f=f\).  Hence \(\Fix(E_f)\le\Stab_{\cG}(f)\).
		
		Now if \(g=\sigma\cdot f\) for some \(\sigma\in\cG\), then \(\Fix(\sigma``E_f)\le\Stab_{\cG}(g)\).
		Taking \(E_d:=\sigma``E_f\) yields the desired countable support for \(d\).
	\end{proof}
	
	\begin{remark}[How Lemma~\ref{lem:coord-countable-support} fits the induction package]
		\label{rem:no-circularity}
		The only inductive input in Lemma~\ref{lem:coord-countable-support} is
		\(\mathbf{Support}_\alpha\) at the predecessor stage \(\alpha\), applied to the
		canonical \(\HS_\alpha\)-name \(\dot f\), together with \(\mathbf{Diag}_\alpha\) for the
		existence of diagonal lifts in \(\cG_\alpha\).
		No property of stage \(\alpha+1\) is used.
	\end{remark}
	
	\begin{definition}[Diagonal-cancellation automorphisms]
		\label{def:diag-cancel}
		Fix an orbit package $\bbS_{[f]}$ (either $\bbQ_{[f]}$ or $\bbR_{[f]}$).
		Let $E\subseteq\omega_1$ be countable and let $D\subseteq [f]$ be countable such that
		\[
		\forall g\in D\ \bigl(\Fix(E)\le \Stab(g)\bigr),
		\]
		i.e.\ every $\pi\in\Fix(E)$ fixes each $g\in D$ as an orbit element.
		For $\pi\in\Fix(E)$ define an automorphism $\widehat\pi^{\,D}\in\Aut(\bbS_{[f]})$ by:
		for each $\vec p\in\bbS_{[f]}$ and each $g\in[f]$,
		\[
		(\widehat\pi^{\,D}\cdot \vec p)(g):=
		\begin{cases}
			\vec p(g), & g\in D,\\
			(\pi\cdot \vec p)(g), & g\notin D,
		\end{cases}
		\]
		where $\pi\cdot\vec p$ is the orbit action from Definition~\ref{def:package-systems} (cf.\ Lemma~\ref{lem:package-actions-wd}).
		
		Let
		\[
		\Delta(E,D):=\{\widehat\pi^{\,D}:\pi\in\Fix(E)\}\ \le\ \Aut(\bbS_{[f]}).
		\]
	\end{definition}
	
	\begin{remark}[Support convention for diagonal lifts]
		\label{rem:diag-support-convention}
		Conditions $\vec p \in \bbS_{[f]}$ are partial functions with countable support
		$\supp(\vec p)\subseteq [f]$.
		Fix an admissible pair $(E,D)$ and $\pi\in\Fix(E)$, and write $\vec p':=\widehat{\pi}^{\,D}\cdot\vec p$.
		Then for each $g\in [f]$ we have the coordinate formula
		\[
		\vec p'(g)=
		\begin{cases}
			\vec p(g), & g\in D,\\[2pt]
			A_{\pi,\pi^{-1}\cdot g}\bigl(\vec p(\pi^{-1}\cdot g)\bigr), & g\notin D,
		\end{cases}
		\]
		i.e.\ $\widehat{\pi}^{\,D}$ agrees with the standard orbit action off $D$ and is the identity on $D$.
		Equivalently, for $h\in [f]\setminus D$ we may write
		\[
		\vec p'(\pi\cdot h)=A_{\pi,h}\bigl(\vec p(h)\bigr).
		\]
		In particular,
		\[
		\supp(\vec p') \;=\;
		(\supp(\vec p)\cap D)\ \cup\ \pi\cdot(\supp(\vec p)\setminus D),
		\]
		which is countable since $\pi$ is a bijection.
	\end{remark}
	
	\begin{definition}[Orbit-package stage filter (diagonal base)]
		\label{def:orbit-stage-filter-diag}
		Let \(\cG_{[f]}\le\Aut(\bbS_{[f]})\) be the subgroup generated by all diagonal-cancellation maps
		\(\widehat\pi^{\,D}\) (equivalently,
		\(\cG_{[f]}=\langle \Delta(E,D):(E,D)\text{ admissible}\rangle\)).
		
		For countable $E\subseteq\omega_1$ and countable $D\subseteq[f]$ such that
		$\Fix(E)\le\Stab(g)$ for all $g\in D$, define
		\[
		\Delta^\uparrow(E,D)\ :=\ \big\langle\ \Delta(E',D')\ \bigm|\ 
		\begin{array}{l}
			E\subseteq E'\in[\omega_1]^{\le\omega},\ D\subseteq D'\subseteq[f]\ \text{countable, and}\\[2pt]
			\Fix(E')\le\Stab(g)\ \ \forall g\in D'
		\end{array}
		\big\rangle,
		\]
		where $\Delta(E',D')$ is as in Definition~\ref{def:diag-cancel}.
		Call a pair \((E,D)\) \emph{admissible} if \(E\subseteq\omega_1\) is countable, \(D\subseteq [f]\) is countable, and
		\(\Fix(E)\leq\Stab(g)\) for every \(g\in D\).
		
		Let $\cF_{[f]}$ be the filter of subgroups of $\cG_{[f]}$ generated by the family
		$\{\Delta^\uparrow(E,D)\mid (E,D)\text{ admissible}\}$.
		Since this family is downward directed (if $(E_i,D_i)$ are admissible then $(E_0\cup E_1,\,D_0\cup D_1)$ is admissible and
		$\Delta^\uparrow(E_0\cup E_1,\,D_0\cup D_1)\le \Delta^\uparrow(E_i,D_i)$ for $i=0,1$),
		the generated filter is just the upward closure of this family;
		equivalently,
		\[
		H\in\cF_{[f]}\quad\Longleftrightarrow\quad
		\exists\,\text{admissible}(E,D)\ \bigl(\Delta^\uparrow(E,D)\le H\bigr).
		\]
	\end{definition}
	
	\begin{remark}
		Lemmas~\ref{lem:diag-conj-standard} and~\ref{lem:diag-conj-diagonal} show that these diagonal generators are closed
		under conjugation by the successor-stage groups (generated by standard and diagonal lifts), hence yield the required
		normal filters in the symmetric-iteration construction.
	\end{remark}
	
	\begin{lemma}[Monotonicity of $\Delta^\uparrow$]\label{lem:diag-up-mono}
		If $(E_0,D_0)$ and $(E_1,D_1)$ are admissible and $E_0\subseteq E_1$ and
		$D_0\subseteq D_1$, then
		\[
		\Delta^\uparrow(E_1,D_1)\ \le\ \Delta^\uparrow(E_0,D_0).
		\]
	\end{lemma}
	
	\begin{proof}
		Any pair $(E',D')$ eligible in the definition of $\Delta^\uparrow(E_1,D_1)$
		satisfies $E_1\subseteq E'$ and $D_1\subseteq D'$, hence also $E_0\subseteq E'$
		and $D_0\subseteq D'$.
		Thus the generating family for $\Delta^\uparrow(E_1,D_1)$
		is a subfamily of the generating family for $\Delta^\uparrow(E_0,D_0)$, so the
		generated subgroup is smaller.
	\end{proof}
	
	\begin{lemma}[$\omega_1$-completeness of the diagonal stage filter]
		\label{lem:omega1complete-diag-filter}
		$\cF_{[f]}$ is $\omega_1$-complete.
	\end{lemma}
	
	\begin{proof}
		Let $\langle H_n:n<\omega\rangle$ be members of $\cF_{[f]}$.
		By Definition~\ref{def:orbit-stage-filter-diag}, for each $n$ choose an admissible
		pair $(E_n,D_n)$ such that
		\[
		\Delta^\uparrow(E_n,D_n)\ \le\ H_n.
		\]
		Let $E:=\bigcup_{n<\omega}E_n$ and $D:=\bigcup_{n<\omega}D_n$; both are countable.
		Moreover $(E,D)$ is admissible: if $g\in D$, pick $n$ with $g\in D_n$, then
		$\Fix(E)\subseteq\Fix(E_n)\le\Stab(g)$.
		By Lemma~\ref{lem:diag-up-mono}, for each $n$ we have
		\[
		\Delta^\uparrow(E,D)\ \le\ \Delta^\uparrow(E_n,D_n)\ \le\ H_n.
		\]
		Hence $\Delta^\uparrow(E,D)\le \bigcap_{n<\omega}H_n$, so $\bigcap_{n<\omega}H_n\in\cF_{[f]}$.
	\end{proof}
	
	\begin{lemma}[Diagonal-cancellation gives a subgroup of automorphisms]
		\label{lem:diag-cancel-subgroup}
		In the setup of Definition~\ref{def:diag-cancel}:
		\begin{enumerate}[label=(\roman*)]
			\item For each $\pi\in\Fix(E)$, $\widehat\pi^{\,D}$ is an automorphism of $\bbS_{[f]}$.
			\item The map $\pi\mapsto \widehat\pi^{\,D}$ is a group monomorphism $\Fix(E)\hookrightarrow \Aut(\bbS_{[f]})$.
			In particular, $\Delta(E,D)$ is a subgroup of $\Aut(\bbS_{[f]})$.
		\end{enumerate}
	\end{lemma}
	
	\begin{proof}
		Fix $\pi\in\Fix(E)$. Since $\pi\cdot(\cdot)$ is an automorphism of $\bbS_{[f]}$ and $\widehat\pi^{\,D}$
		agrees with it on all coordinates outside $D$ while acting as the identity on coordinates in $D$,
		it preserves supports, preserves coordinatewise extension, and is bijective with inverse
		$\widehat{\pi^{-1}}^{\,D}$.
		This proves (i).
		
		For (ii), $\widehat\pi^{\,D}\circ \widehat\sigma^{\,D}=\widehat{\pi\sigma}^{\,D}$ follows
		coordinatewise from the definition (on $D$ both sides are the identity; off $D$ both sides agree with the rigidified orbit action
		(Lemma~\ref{lem:package-actions-wd})).
		Injectivity is immediate.
	\end{proof}
	
	\begin{lemma}[Conjugation of diagonal-cancellation groups by standard lifts]
		\label{lem:diag-conj-standard}
		Let $(E,D)$ be admissible.
		For $\sigma\in\Sym(\omega_1)$ write
		\[
		\sigma\cdot E:=\{\sigma(\xi):\xi\in E\}\subseteq\omega_1,
		\qquad
		\sigma\cdot D:=\{\sigma\cdot g:g\in D\}\subseteq[f].
		\]
		Then $(\sigma\cdot E,\sigma\cdot D)$ is admissible. Moreover, for every $\pi\in\Fix(E)$,
		\[
		\widehat{\sigma}^{\,\emptyset}\circ \widehat{\pi}^{\,D}\circ(\widehat{\sigma}^{\,\emptyset})^{-1}
		\;=\;
		\widehat{\sigma\pi\sigma^{-1}}^{\,\sigma\cdot D}.
		\]
		Consequently,
		\[
		\widehat{\sigma}^{\,\emptyset}\,\Delta(E,D)\,(\widehat{\sigma}^{\,\emptyset})^{-1}
		=\Delta(\sigma\cdot E,\sigma\cdot D),
		\qquad
		\widehat{\sigma}^{\,\emptyset}\,\Delta^\uparrow(E,D)\,(\widehat{\sigma}^{\,\emptyset})^{-1}
		=\Delta^\uparrow(\sigma\cdot E,\sigma\cdot D).
		\]
	\end{lemma}
	
	\begin{proof}
		\emph{Admissibility.}
		Let $\tau\in\Fix(\sigma\cdot E)$ and let $g\in\sigma\cdot D$, say $g=\sigma\cdot h$ with $h\in D$.
		Then $\sigma^{-1}\tau\sigma\in\Fix(E)\le\Stab(h)$, hence
		\[
		\tau\cdot g=\tau\cdot(\sigma\cdot h)=\sigma\cdot\bigl((\sigma^{-1}\tau\sigma)\cdot h\bigr)=\sigma\cdot h=g.
		\]
		So $\Fix(\sigma\cdot E)\le\Stab(g)$ for all $g\in\sigma\cdot D$.
		
		\emph{Conjugation formula.}
		Fix $\pi\in\Fix(E)$ and $\vec p\in\bbS_{[f]}$, and let $g\in[f]$. Put $h:=\sigma^{-1}\cdot g$.
		If $g\in\sigma\cdot D$ (equivalently $h\in D$), then $\widehat{\pi}^{\,D}$ acts as the identity at $h$,
		so the conjugate acts as the identity at $g$, i.e. it fixes the $g$-coordinate of $\vec p$.
		If $g\notin\sigma\cdot D$ (equivalently $h\notin D$), then at $h$ the map $\widehat{\pi}^{\,D}$ agrees with the
		standard orbit action by $\pi$.
		Using that $\widehat{\sigma}^{\,\emptyset}$ is the standard orbit action by $\sigma$
		(i.e.\ $D=\emptyset$), a direct coordinate computation gives that the conjugate agrees with the standard orbit action by
		$\sigma\pi\sigma^{-1}$ at $g$.
		Thus the conjugate acts as the identity on $\sigma\cdot D$ and as $\sigma\pi\sigma^{-1}$
		off $\sigma\cdot D$, exactly $\widehat{\sigma\pi\sigma^{-1}}^{\,\sigma\cdot D}$.
		
		The subgroup identity for $\Delta(E,D)$ follows by taking images of $\Fix(E)$ under $\pi\mapsto\widehat{\pi}^{\,D}$ and
		using $\sigma\Fix(E)\sigma^{-1}=\Fix(\sigma\cdot E)$.
		Finally, for $\Delta^\uparrow(E,D)$: conjugating any generator $\Delta(E',D')$ occurring in the definition of
		$\Delta^\uparrow(E,D)$ yields $\Delta(\sigma\cdot E',\sigma\cdot D')$, which is a generator for
		$\Delta^\uparrow(\sigma\cdot E,\sigma\cdot D)$.
		Applying the same argument with $\sigma^{-1}$ gives equality.
	\end{proof}
	
	\begin{lemma}[Conjugation of diagonal generators by diagonal lifts]
		\label{lem:diag-conj-diagonal}
		Fix an admissible pair $(E_0,D_0)$ and let $\sigma\in\Fix(E_0)$.
		Write $\widehat{\sigma}^{\,D_0}\in\Delta(E_0,D_0)$ for the corresponding diagonal lift.
		
		Then for every admissible pair $(E,D)$, letting
		\[
		E^\ast:=E\cup E_0,\qquad D^\ast:=D\cup D_0,
		\]
		we have
		\[
		\widehat{\sigma}^{\,D_0}\cdot \Delta^\uparrow(E^\ast,D^\ast)\cdot (\widehat{\sigma}^{\,D_0})^{-1}
		\;=\;
		\Delta^\uparrow(\sigma\cdot E^\ast,\sigma\cdot D^\ast).
		\]
		In particular, $\widehat{\sigma}^{\,D_0}\,H\,(\widehat{\sigma}^{\,D_0})^{-1}\in\cF_{[f]}$ for every $H\in\cF_{[f]}$.
	\end{lemma}
	
	\begin{proof}
		Let $(E,D)$ be admissible and set $E^\ast:=E\cup E_0$ and $D^\ast:=D\cup D_0$.
		Then $(E^\ast,D^\ast)$ is admissible since $\Fix(E^\ast)\le\Fix(E)$ and $\Fix(E^\ast)\le\Fix(E_0)$.
		
		We first claim that if $(E',D')$ is admissible with $E^\ast\subseteq E'$ and $D^\ast\subseteq D'$, then
		\[
		\widehat{\sigma}^{\,D_0}\cdot \Delta(E',D')\cdot (\widehat{\sigma}^{\,D_0})^{-1}
		\;=\;
		\Delta(\sigma\cdot E',\sigma\cdot D').
		\]
		Indeed, $D_0\subseteq D'$ by assumption, so every $\widehat{\pi}^{\,D'}\in\Delta(E',D')$ is the identity on all
		coordinates in $D_0$.
		Since $\widehat{\sigma}^{\,D_0}$ agrees with the standard lift $\widehat{\sigma}^{\,\emptyset}$
		off $D_0$, the same coordinate computation as in Lemma~\ref{lem:diag-conj-standard} gives
		\[
		\widehat{\sigma}^{\,D_0}\cdot \widehat{\pi}^{\,D'}\cdot (\widehat{\sigma}^{\,D_0})^{-1}
		\;=\;
		\widehat{\sigma\pi\sigma^{-1}}^{\,\sigma\cdot D'}.
		\]
		Taking images of $\Fix(E')$ under $\pi\mapsto\widehat{\pi}^{\,D'}$ and using
		$\sigma\Fix(E')\sigma^{-1}=\Fix(\sigma\cdot E')$ yields the subgroup identity above.
		
		Now conjugate the generating family for $\Delta^\uparrow(E^\ast,D^\ast)$:
		each generator $\Delta(E',D')$ (with $E^\ast\subseteq E'$ and $D^\ast\subseteq D'$) is sent to
		$\Delta(\sigma\cdot E',\sigma\cdot D')$, which is a generator for
		$\Delta^\uparrow(\sigma\cdot E^\ast,\sigma\cdot D^\ast)$.
		Applying the same argument with $\sigma^{-1}$ gives equality.
		
		Finally, if $H\in\cF_{[f]}$, choose admissible $(E,D)$ with $\Delta^\uparrow(E,D)\le H$.
		Then $\Delta^\uparrow(E^\ast,D^\ast)\le \Delta^\uparrow(E,D)\le H$ by monotonicity, so
		\[
		\Delta^\uparrow(\sigma\cdot E^\ast,\sigma\cdot D^\ast)
		\;=\;
		\widehat{\sigma}^{\,D_0}\cdot \Delta^\uparrow(E^\ast,D^\ast)\cdot (\widehat{\sigma}^{\,D_0})^{-1}
		\;\le\;
		\widehat{\sigma}^{\,D_0}\,H\,(\widehat{\sigma}^{\,D_0})^{-1},
		\]
		which shows $\widehat{\sigma}^{\,D_0}\,H\,(\widehat{\sigma}^{\,D_0})^{-1}\in\cF_{[f]}$.
	\end{proof}
	
	\subsection{Iteration-level diagonal automorphisms}
	\label{subsec:iteration-diag}
	
	The diagonal-cancellation construction of Definition~\ref{def:diag-cancel} is formulated 
	for a single orbit package $\bbS_{[f]}$.
	For the $\neg\AC$ argument at limit stages, 
	we need diagonal lifts acting on the \emph{full iteration forcing} $\bbP_\lambda$, 
	with protection sets $D$ that may span multiple packages across multiple stages.
	
	\begin{remark}[Stagewise induction package and dependency order]
		\label{rem:simultaneous-induction}
		For each stage \(\lambda\), the arguments in this subsection verify the following assertions
		simultaneously:
		
		\begin{description}
			\item[\(\mathbf{Diag}_\lambda\)]
			Every iteration-level diagonal lift \(\widehat{\pi}^{\,D}\) with \((E,D)\) globally admissible
			at stage \(\lambda\) and \(\pi\in\Fix(E)\) belongs to \(\cG_\lambda\).
			
			\item[\(\mathbf{Support}_\lambda\)]
			Every relevant hereditarily symmetric \(\bbP_\lambda\)-name \(\tau\) admits a countable globally admissible
			pair \((E,D)\) at stage \(\lambda\) such that
			\[
			\Delta_\lambda^\uparrow(E,D)\subseteq \Stab_{\cG_\lambda}(\tau).
			\]
			
			\item[\(\mathbf{Hull}_\lambda\)]
			Every countable coordinate set \(D\subseteq \mathrm{Coords}_{<\lambda}\) admits a countable admissible hull:
			for every countable \(E\subseteq\omega_1\), there is countable \(E^\ast\supseteq E\) such that
			\((E^\ast,D)\) is globally admissible at stage \(\lambda\).
			
			\item[\(\mathbf{Core}_\lambda\)]
			Every \(H\in\tilde{\cF}^*_\lambda\) contains some
			\(\Delta_\lambda^\uparrow(E,D)\) with \((E,D)\) countable and globally admissible at stage \(\lambda\).
			
			\item[\(\mathbf{Filter}_\lambda\)]
			The modified limit filter \(\tilde{\cF}^*_\lambda\) is proper, normal, and \(\omega_1\)-complete.
			
			\item[\(\mathbf{Delta}_\lambda\)]
			A \(\bbP_\lambda\)-name is hereditarily symmetric iff its stabilizer contains
			\(\Delta_\lambda^\uparrow(E,D)\) for some countable globally admissible pair \((E,D)\).
		\end{description}
		
		The dependency order is strictly well-founded and is the one used implicitly in the proofs below:
		\begin{enumerate}[label=(\roman*)]
			\item \(\mathbf{Diag}_\lambda\) uses only the stage-\(\lambda\) iteration/group recursion and projection coherence.
			\item \(\mathbf{Support}_\alpha\) at a predecessor stage \(\alpha\) is the only earlier-stage support input used in
			Lemma~\ref{lem:coord-countable-support}.
			\item \(\mathbf{Hull}_\lambda\) uses Lemma~\ref{lem:coord-countable-support} for the coordinates appearing in \(D\).
			\item \(\mathbf{Core}_\lambda\) uses \(\mathbf{Core}_\beta\) only for \(\beta<\lambda\) in the pullback-generator case,
			and uses \(\mathbf{Hull}_\lambda\) in the conjugate \(\Delta^\uparrow\)-generator case.
			\item \(\mathbf{Filter}_\lambda\) uses \(\mathbf{Core}_\lambda\) and countable-union admissibility.
			\item \(\mathbf{Delta}_\lambda\) is then the reformulation of hereditary symmetry obtained from \(\mathbf{Core}_\lambda\).
		\end{enumerate}
		
		Thus no argument at stage \(\lambda\) appeals to a stage-\(\lambda\) clause that is logically later in this list,
		and no argument appeals to any stage \(>\lambda\).
	\end{remark}
	
	\begin{lemma}[Successor decomposition and splitting]
		\label{lem:successor-splitting}
		Let $\alpha$ be an ordinal. In the two-step presentation
		$\bbP_{\alpha+1} \cong \bbP_\alpha \ast \dot{\bbQ}_\alpha$, the projection
		$\pi^\alpha_{\alpha+1}: \cG_{\alpha+1} \to \cG_\alpha$ is a surjective group homomorphism.
		Let $K_{\alpha+1} := \ker(\pi^\alpha_{\alpha+1})$.
		
		Moreover, there is a group-theoretic splitting map
		$\mathrm{spl}: \cG_\alpha \to \cG_{\alpha+1}$ (the \emph{splitting section}) such that
		$\pi^\alpha_{\alpha+1} \circ \mathrm{spl} = \mathrm{id}_{\cG_\alpha}$, and every
		$g \in \cG_{\alpha+1}$ can be written as $g = \mathrm{spl}(h) \cdot k$ with
		$h = \pi^\alpha_{\alpha+1}(g)$ and $k \in K_{\alpha+1}$.
	\end{lemma}
	
	\begin{proof}
		This is exactly the successor-step group structure supplied by the bounded-stage
		countable-support symmetric iteration framework (cf.\\Remark~\ref{rem:iteration-api}(b)).
		Concretely:
		\begin{itemize}
			\item $\pi^\alpha_{\alpha+1}$ is restriction to the $\bbP_\alpha$-part;
			\item $K_{\alpha+1}$ consists of automorphisms acting trivially on $\bbP_\alpha$;
			\item $\mathrm{spl}(h)$ is the splitting section that acts as $h$ on $\bbP_\alpha$ and 
			trivially on $\dot{\bbQ}_\alpha$.
		\end{itemize}
		The decomposition $g = \mathrm{spl}(\pi^\alpha_{\alpha+1}(g)) \cdot k$ with 
		$k = \mathrm{spl}(\pi^\alpha_{\alpha+1}(g))^{-1} \cdot g$ is immediate;
		one checks 
		$k \in K_{\alpha+1}$ since $\pi^\alpha_{\alpha+1}(k) = \mathrm{id}$.
	\end{proof}
	
	\begin{remark}[Projection maps are coherent]
		\label{rem:projection-coherence}
		For $\gamma<\beta$ we write $\pi^\gamma_\beta:\cG_\beta\to\cG_\gamma$ for the canonical projection
		homomorphisms of the iteration groups.
		These satisfy coherence:
		for all $\gamma<\beta<\lambda$,
		\[
		\pi^\gamma_\beta\circ\pi^\beta_\lambda=\pi^\gamma_\lambda.
		\]
	\end{remark}
	
	\begin{definition}[Stage-$0$ projection]
		\label{def:stage0-projection}
		Set $\cG_0:=\cG=\Sym(\omega_1)$. For each stage $\alpha$ define $\pi^0_\alpha:\cG_\alpha\to\cG$
		by transfinite recursion:
		\[
		\pi^0_0=\id,\qquad \pi^0_{\alpha+1}=\pi^0_\alpha\circ\pi^\alpha_{\alpha+1},
		\]
		and for limit $\lambda$ set
		\[
		\pi^0_\lambda(g):=\pi^0_\alpha\bigl(\pi^\alpha_\lambda(g)\bigr)\quad\text{for any }\alpha<\lambda.
		\]
		This is well-defined: if \(\alpha<\beta<\lambda\), then by projection coherence (Remark~\ref{rem:projection-coherence}),
		\[
		\pi^\alpha_\beta\bigl(\pi^\beta_\lambda(g)\bigr)=\pi^\alpha_\lambda(g),
		\]
		and applying \(\pi^0_\alpha\) gives
		\[
		\pi^0_\alpha\bigl(\pi^\beta_\lambda(g)\bigr)
		=\pi^0_\alpha\bigl(\pi^\alpha_\beta(\pi^\beta_\lambda(g))\bigr)
		=\pi^0_\alpha\bigl(\pi^\alpha_\lambda(g)\bigr),
		\]
		so the value \(\pi^0_\alpha(\pi^\alpha_\lambda(g))\) is independent of the chosen \(\alpha\).
		In particular, for \(\beta<\lambda\),
		\[
		\pi^0_\beta\circ\pi^\beta_\lambda=\pi^0_\lambda.
		\]
	\end{definition}
	
	\begin{definition}[Limit-stage groups for countable support]
		\label{def:limit-stage-group}
		For a limit ordinal $\lambda < \Ord$, the \emph{limit group} $\cG_\lambda$ for the countable-support iteration is the inverse limit of the system $\langle \cG_\beta, \pi^\gamma_\beta \mid \gamma < \beta < \lambda \rangle$; explicitly,
		\[
		\cG_\lambda \;=\; \bigl\{\langle h_\beta \rangle_{\beta<\lambda} \in \prod_{\beta<\lambda} \cG_\beta \;:\; \forall \gamma < \beta < \lambda,\ \pi^\gamma_\beta(h_\beta) = h_\gamma \bigr\},
		\]
		where $\pi^\gamma_\beta : \cG_\beta \to \cG_\gamma$ is the canonical restriction homomorphism from the iteration system. This is the standard inverse limit construction; cf.\ \cite[Def.~4.4, Lemma~4.5]{KaragilaSchilhan} for $I$-support iterations with $I = [\lambda]^{\le\omega}$.
	\end{definition}
	
	\begin{definition}[Package coordinate sets]
		\label{def:coords}
		We define by transfinite recursion the sets $\mathrm{Coords}_\alpha$ of
		\emph{package coordinates occurring below stage $\alpha$}.
		\begin{itemize}
			\item $\mathrm{Coords}_0:=\varnothing$.
			\item At a successor stage $\alpha+1$, write $\bbP_{\alpha+1}\cong \bbP_\alpha * \dot{\bbQ}_\alpha$,
			where $\dot{\bbQ}_\alpha$ denotes the \emph{stage iterand} (a finite product of orbit packages) and
			should not be confused with the localized \(\PP\)-package forcing \(\bbQ_f\).
			In this paper, $\dot{\bbQ}_\alpha$ is (forced to be) a finite product of orbit packages
			\[
			\dot{\bbQ}_\alpha \cong \prod_{i<m_\alpha} \dot{\bbS}_{[f_{\alpha,i}]},
			\]
			where each $\dot{\bbS}_{[f_{\alpha,i}]}$ is either an orbit $\PP$-package $\bbQ_{[f]}$
			or an orbit $\AC_{\WO}$-package $\bbR_{[f]}$ (and $m_\alpha<\omega$).
			
			Define the stage-$\alpha$ index set as the disjoint union
			\[
			I_\alpha:=\bigsqcup_{i<m_\alpha} [f_{\alpha,i}].
			\]
			We set
			\[
			\begin{split}
			\mathrm{NewCoords}_{\alpha+1}&:=\{\langle \alpha,z\rangle : z\in I_\alpha\},
			\\
			\mathrm{Coords}_{\alpha+1}&:=\mathrm{Coords}_\alpha\ \dot\cup\ \mathrm{NewCoords}_{\alpha+1}.
			\end{split}
			\]
			When $z\in I_\alpha$ corresponds to $(i,g)$, we also write $\langle \alpha,i,g\rangle$ for $\langle \alpha,z\rangle$.
			
			\item At a limit stage $\lambda$, set $\mathrm{Coords}_\lambda:=\bigcup_{\alpha<\lambda}\mathrm{Coords}_\alpha$.
			We also write $\mathrm{Coords}_{<\lambda}:=\mathrm{Coords}_\lambda$.
		\end{itemize}
		
	When working at a fixed successor stage $\alpha+1$, we suppress the stage tag and
	write $\langle i,g\rangle$ for $\langle \alpha,i,g\rangle$ (so “$\langle i,g\rangle\in
	\mathrm{Coords}_{\alpha+1}\setminus\mathrm{Coords}_\alpha$” means $i<m_\alpha$ and
	$g\in[f_{\alpha,i}]$).
	We let $\cG=\Sym(\omega_1)$ act on new coordinates by
	\[
	\pi\cdot\langle \alpha,i,g\rangle := \langle \alpha,i,\pi\cdot g\rangle,
	\]
	equivalently $\pi\cdot\langle \alpha,z\rangle := \langle \alpha,\pi\cdot z\rangle$ when using $I_\alpha$.
	\end{definition}
	
	\begin{definition}[Successor diagonal kernel element]
		\label{def:succ-diag-kernel}
		Fix a successor stage $\alpha+1$. Let $\pi \in \cG=\Sym(\omega_1)$ and let
		$D^{\mathrm{new}} \subseteq \mathrm{Coords}_{\alpha+1} \setminus \mathrm{Coords}_\alpha$
		be a set of package coordinates newly introduced at stage $\alpha+1$.
		
		\textbf{Assume} $\pi$ fixes each $d \in D^{\mathrm{new}}$ under the induced action on
		tagged coordinates, i.e., if $d=\langle \alpha,i,g\rangle$ then $\pi\cdot d=d$,
		equivalently $\pi\cdot g=g$.
		(Here $\pi\cdot\langle\alpha,i,g\rangle:=\langle\alpha,i,\pi\cdot g\rangle$ using the orbit action on $[f_{\alpha,i}]$.)
		
		Under this assumption, define $\kappa_{\alpha+1}(\pi, D^{\mathrm{new}}) \in \Aut(\bbP_{\alpha+1})$ by:
		\begin{itemize}
			\item $\kappa_{\alpha+1}(\pi, D^{\mathrm{new}})$ acts as the identity on the 
			$\bbP_\alpha$-part (all coordinates at stages $\le \alpha$);
			\item on each new package coordinate $d \in \mathrm{Coords}_{\alpha+1} \setminus \mathrm{Coords}_\alpha$:
			\begin{itemize}
				\item if $d \notin D^{\mathrm{new}}$, act by the package action induced by $\pi$ 
				(Definition~\ref{def:package-systems});
				\item if $d \in D^{\mathrm{new}}$, act as the identity (fix the $d$-coordinate pointwise).
			\end{itemize}
		\end{itemize}
	\end{definition}
	
	\begin{remark}[Why the hypothesis is necessary]
		\label{rem:kappa-hypothesis}
		The hypothesis ``$\pi \cdot d = d$ for all $d \in D^{\mathrm{new}}$'' ensures 
		$\kappa_{\alpha+1}(\pi, D^{\mathrm{new}})$ is a bijection on package coordinates.
		Without it, collisions can occur: if $d \notin D^{\mathrm{new}}$ but 
		$\pi \cdot d \in D^{\mathrm{new}}$, then both $d$ and $\pi \cdot d$ would map to $\pi \cdot d$.
		In our application (Lemma~\ref{lem:diag-lift-in-group}), global admissibility 
		guarantees this hypothesis: if $(E, D)$ is globally admissible and $\pi \in \Fix(E)$, 
		then $\pi \cdot d = d$ for all $d \in D$, hence for all 
		$d \in D^{\mathrm{new}} = D_{\beta+1} \setminus D_\beta \subseteq D$.
	\end{remark}
	
	\begin{lemma}[Successor diagonal kernel element lies in $\cG_{\alpha+1}$ and its kernel]
		\label{lem:succ-diag-kernel-in-kernel}
		Under the hypothesis of Definition~\ref{def:succ-diag-kernel}:
		\begin{enumerate}[label=(\roman*)]
			\item $\kappa_{\alpha+1}(\pi, D^{\mathrm{new}}) \in \cG_{\alpha+1}$;
			\item $\kappa_{\alpha+1}(\pi, D^{\mathrm{new}}) \in K_{\alpha+1} = \ker(\pi^\alpha_{\alpha+1})$.
		\end{enumerate}
		Moreover, for fixed $D^{\mathrm{new}}$, the map
		\[
		\pi \longmapsto \kappa_{\alpha+1}(\pi, D^{\mathrm{new}})
		\]
		is a group homomorphism on its domain
		\[
		\{\pi \in \cG=\Sym(\omega_1) : \forall d\in D^{\mathrm{new}}\ (\pi \cdot d = d)\},
		\]
		where the action on new coordinates is the one fixed in Definition~\ref{def:coords}.
	\end{lemma}
	
	\begin{proof}
		(i) By Definition~\ref{def:succ-diag-kernel}, the map $\kappa_{\alpha+1}(\pi, D^{\mathrm{new}})$
		acts trivially on the $\bbP_\alpha$-part, and on the successor iterand $\dot{\bbQ}_\alpha$
		it acts coordinatewise on the new package coordinates
		$\mathrm{NewCoords}_{\alpha+1}$: it is the identity on $D^{\mathrm{new}}$ and it agrees with
		the $\pi$-induced orbit action on $\mathrm{NewCoords}_{\alpha+1}\setminus D^{\mathrm{new}}$.
		If (as allowed by Remark~\ref{rem:package-product}) $\dot{\bbQ}_\alpha$ is a finite product of orbit
		packages, then $\mathrm{NewCoords}_{\alpha+1}$ is tagged by the factor index, and the action is
		defined independently on each factor, hence yields an automorphism of the full stage-$({\alpha+1})$
		forcing.
		By Lemma~\ref{lem:diag-cancel-subgroup}, each factor is a legitimate automorphism 
		of its package forcing $\bbS_{[f]}$.
		By the successor-stage group structure (Remark~\ref{rem:iteration-api}), $\cG_{\alpha+1}$ is generated by lifts from
		$\cG_\alpha$ together with automorphisms of the successor iterand $\dot{\bbQ}_\alpha$ (equivalently, of each new orbit-package factor).
		Since $\kappa_{\alpha+1}(\pi, D^{\mathrm{new}})$ acts as identity on $\bbP_\alpha$ and 
		as a (product of) diagonal-cancellation automorphism(s) on new packages, it lies 
		in $\cG_{\alpha+1}$.
		
		(ii) By Definition~\ref{def:succ-diag-kernel}, $\kappa_{\alpha+1}(\pi, D^{\mathrm{new}})$
		acts as the identity on $\bbP_\alpha$.
		Hence 
		$\pi^\alpha_{\alpha+1}(\kappa_{\alpha+1}(\pi, D^{\mathrm{new}})) = \mathrm{id}$, 
		so $\kappa_{\alpha+1}(\pi, D^{\mathrm{new}}) \in \ker(\pi^\alpha_{\alpha+1}) = K_{\alpha+1}$.
		
		The homomorphism property: for $\pi, \sigma$ both fixing $D^{\mathrm{new}}$ pointwise,
		\[
		\kappa_{\alpha+1}(\pi, D^{\mathrm{new}}) \circ \kappa_{\alpha+1}(\sigma, D^{\mathrm{new}}) 
		= \kappa_{\alpha+1}(\pi\sigma, D^{\mathrm{new}}),
		\]
		since on coordinates outside $D^{\mathrm{new}}$, composing the $\sigma$-action 
		with the $\pi$-action yields the $(\pi\sigma)$-action, and on $D^{\mathrm{new}}$ 
		both sides act as identity.
	\end{proof}
	
	\begin{definition}[Global admissibility]
		\label{def:global-admissible}
		Let $\lambda$ be a stage of the iteration. A pair $(E,D)$ is \emph{globally admissible 
			at stage $\lambda$} if:
		\begin{enumerate}[label=(\roman*)]
			\item $E \in [\omega_1]^{\le\omega}$ is countable;
			\item $D\subseteq \mathrm{Coords}_{<\lambda}$ is countable;
			\item $\Fix(E) \le \Stab_{\cG}(d)$ for every $d \in D$, where $\Stab_{\cG}(d)$ 
			denotes the stabilizer under the stage-$0$ action of $\cG=\Sym(\omega_1)$ on coordinates
			(Definition~\ref{def:coords}).
		\end{enumerate}
		When $D = \emptyset$, condition (iii) is vacuously satisfied, and $(E,\emptyset)$ is 
		globally admissible for any countable $E$.
	\end{definition}
	
	\begin{definition}[Iteration-level diagonal subgroups and diagonal lifts]
		\label{def:iteration-diag}
		Fix a stage $\lambda$.
		
		\smallskip\noindent
		\textbf{(A) Characterized diagonal subgroups.}
		For any $E\subseteq\omega_1$ and any $D\subseteq\mathrm{Coords}_{<\lambda}$ define
		\[
		\Delta_\lambda^\uparrow(E,D)
		:=\bigl\{\,g\in\cG_\lambda:\pi^0_\lambda(g)\in\Fix(E)\ \text{and $g$ fixes every coordinate in $D$ pointwise}\,\bigr\}.
		\]
		
		\smallskip\noindent
		\textbf{(B) Diagonal lifts (require global admissibility).}
		If $(E,D)$ is globally admissible at stage $\lambda$ and $\pi\in\Fix(E)$, define
		$\widehat{\pi}^{\,D}\in\Aut(\bbP_\lambda)$ by:
		\begin{itemize}
			\item \textbf{Stage $0$ (Cohen base):} act as $\pi$, i.e.\
			$(\widehat{\pi}^{\,D}\cdot p)(\pi(\alpha),n)=p(\alpha,n)$ for $(\alpha,n)\in\omega_1\times\omega$.
			\item \textbf{Stages $>0$ (package coordinates):} for each package coordinate $g$,
			act as the identity if $g\in D$, and as the $\pi$-orbit action if $g\notin D$.
		\end{itemize}
		Finally set
		\[
		\Delta_\lambda(E,D):=\{\widehat{\pi}^{\,D}:\pi\in\Fix(E)\}\le\Aut(\bbP_\lambda).
		\]
	\end{definition}
	
	\begin{remark}\label{rem:diag-lift-in-delta-forward}
		Let $(E,D)$ be globally admissible at stage \(\lambda\) and let \(\pi\in\Fix(E)\).
		The diagonal lift \(\widehat{\pi}^{\,D}\) is defined as an element of \(\Aut(\bbP_\lambda)\).
		Once we verify that \(\widehat{\pi}^{\,D}\in\cG_\lambda\) (Lemma~\ref{lem:diag-lift-in-group}),
		it follows immediately from the defining clause (A) that
		\(\widehat{\pi}^{\,D}\in\Delta_\lambda^\uparrow(E,D)\).
	\end{remark}
	
	\begin{remark}[Characterization vs.\ generation]
		\label{rem:delta-up-characterization}
		Definition~\ref{def:iteration-diag} defines $\Delta_\lambda^\uparrow(E,D)$ as a 
		\emph{characterized subgroup}: $g \in \Delta_\lambda^\uparrow(E,D)$ iff $g$'s stage-$0$ 
		projection fixes $E$ pointwise and $g$ fixes every package coordinate in $D$ pointwise.
		This ``characterization'' definition makes monotonicity immediate 
		(Lemma~\ref{lem:iteration-diag-mono}). In Lemma~\ref{lem:diag-lift-in-delta}, we show that once diagonal lifts
		$\widehat{\pi}^{\,D'}$ are known to belong to $\cG_\lambda$ (Lemma~\ref{lem:diag-lift-in-group}),
		then for $\pi \in \Fix(E)$ and $D'\supseteq D$ we have
		$\widehat{\pi}^{\,D'}\in\Delta_\lambda^\uparrow(E,D)$;
		this is then immediate from the characterization. 
		We do \emph{not} need to prove that $\Delta_\lambda^\uparrow(E,D)$ is \emph{generated by} 
		such diagonal lifts;
		the characterization suffices for all filter-theoretic purposes.
	\end{remark}
	
	\begin{lemma}[Group action on coordinates]
		\label{lem:group-action-coords}
		Let $k \in \cG_\lambda$.
		Then:
		\begin{enumerate}[label=(\roman*)]
			\item The stage-$0$ projection $\pi^0_\lambda: \cG_\lambda \to \cG$ is a group homomorphism.
			\item The element $k$ induces a bijection on the set of package coordinates at each stage $\beta<\lambda$
			(in the sense of Definition~\ref{def:coords}).
			In particular, if $D$ is a countable set of
			package coordinates, then $k\cdot D$ is countable.
			\item For $E\subseteq\omega_1$, define $k\cdot E:=\pi^0_\lambda(k)[E]$.
			Then $k\cdot E$ is countable whenever $E$ is countable.
		\end{enumerate}
	\end{lemma}
	
	\begin{proof}
		(i) This is immediate from Definition~\ref{def:stage0-projection} together with the coherence of the
		projections (Remark~\ref{rem:projection-coherence}); cf. Definition~\ref{def:limit-stage-group}.
		
		(ii) We prove by induction on stages that each $k\in\cG_\beta$ determines a bijection of
		$\mathrm{Coords}_\beta$, and hence sends countable sets of coordinates to countable sets.
		At stage $0$ there are no package coordinates. At a successor stage $\beta+1$, write
		$k=\mathrm{spl}(h)\cdot u$ as in Lemma~\ref{lem:successor-splitting}, where
		$h=\pi^\beta_{\beta+1}(k)\in\cG_\beta$ and $u\in K_{\beta+1}$ acts trivially on the $\bbP_\beta$-part.
		By the inductive hypothesis, $h$ induces a bijection of $\mathrm{Coords}_\beta$.
		The standard lift $\mathrm{spl}(h)$ extends this action and fixes every element of
		$\mathrm{NewCoords}_{\beta+1}$.
		On the other hand, $u$ fixes $\mathrm{Coords}_\beta$ pointwise and acts on the successor iterand
		$\dot{\bbQ}_\beta$ by compositions of package automorphisms (Lemma~\ref{lem:package-actions-wd})
		and diagonal-cancellation automorphisms (Definition~\ref{def:diag-cancel}), each of which induces
		a bijection on the new coordinate index set.
		Therefore $k$ induces a bijection on
		$\mathrm{Coords}_{\beta+1}=\mathrm{Coords}_\beta\dot\cup \mathrm{NewCoords}_{\beta+1}$.
		
		At a limit stage $\eta$, elements of $\cG_\eta$ are coherent inverse-limit elements whose restrictions
		to earlier stages determine them uniquely (cf.\ [Definition~\ref{def:limit-stage-group}]).
		The induced bijection on $\mathrm{Coords}_\eta=\bigcup_{\beta<\eta}\mathrm{Coords}_\beta$ is obtained
		by taking the union of the stagewise bijections, and is again a bijection.
		Consequently, if $D$ is countable then so is $k\cdot D$.
		
		(iii) Immediate, since $\pi^0_\lambda(k)\in\cG=\Sym(\omega_1)$ is a bijection and hence preserves countability.
	\end{proof}
	
	\begin{lemma}[Iteration diagonal lifts are automorphisms]
		\label{lem:iteration-diag-auto}
		For $(E,D)$ globally admissible at stage $\lambda$ and $\pi \in \Fix(E)$, 
		the map $\widehat{\pi}^{\,D}$ is an automorphism of $\bbP_\lambda$ satisfying:
		\begin{enumerate}[label=(\roman*)]
			\item $\widehat{\pi}^{\,D}$ restricts to $\pi$ on the Cohen base $\Add(\omega,\omega_1)$;
			\item $\widehat{\pi}^{\,D}$ restricts to the identity on the $g$-coordinate factor for each $g \in D$;
			\item $\widehat{\pi}^{\,D}$ restricts to the standard $\pi$-action on all coordinates outside $D$.
		\end{enumerate}
		Moreover, $\pi \mapsto \widehat{\pi}^{\,D}$ is a group monomorphism $\Fix(E) \hookrightarrow \Aut(\bbP_\lambda)$.
	\end{lemma}
	
	\begin{proof}
		By global admissibility, $\pi \cdot g = g$ for all $g \in D$, so the protected coordinates are
		pointwise fixed and the diagonal definition cannot create coordinate-collisions.
		Order preservation and bijectivity follow coordinatewise: 
		on the Cohen base this is the standard automorphism;
		on each package coordinate outside $D$, 
		this is the orbit action from Lemma~\ref{lem:package-actions-wd};
		on coordinates in $D$, 
		this is the identity.
		
		The group homomorphism property $\widehat{\pi}^{\,D} \circ \widehat{\sigma}^{\,D} = \widehat{\pi\sigma}^{\,D}$ 
		follows coordinatewise: on $D$ both sides are the identity;
		on the Cohen base and on coordinates 
		outside $D$, both sides agree with the standard action of $\pi\sigma$.
		Injectivity is immediate from (i): if $\widehat{\pi}^{\,D}=\widehat{\sigma}^{\,D}$ then their
		restrictions to the Cohen base agree, hence $\pi=\sigma$.
	\end{proof}
	
	\begin{corollary}[Stage-$0$ projection controls the Cohen reals]
		\label{cor:diag-lift-swaps-cohen}
		Let $g\in\cG_\lambda$ and put $\pi:=\pi^0_\lambda(g)\in\cG=\Sym(\omega_1)$.
		Then for each $\alpha<\omega_1$,
		\[
		g\cdot \dot c_\alpha=\dot c_{\pi(\alpha)},
		\]
		where $\dot c_\alpha$ is the canonical $\bbP_0$-name for the $\alpha$-th Cohen real.
		In particular, if $\pi=(\beta\ \gamma)$ is a transposition with $\beta\neq\gamma$, then $g$
		swaps $c_\beta\leftrightarrow c_\gamma$ in the generic extension.
	\end{corollary}
	
	\begin{proof}
		View $\dot c_\alpha$ as a $\bbP_\lambda$-name via the canonical inclusion of
		$\bbP_0$-names into $\bbP_\lambda$-names. Recall
		\[
		\dot c_\alpha:=\{\langle \check n,p\rangle : p(\alpha,n)=1\}.
		\]
		Since $\pi^0_\lambda(g)=\pi$, the restriction of $g$ to the Cohen base $\bbP_0=\Add(\omega,\omega_1)$
		is the usual $\pi$-action.
		Therefore $g\cdot \dot c_\alpha=\pi\cdot \dot c_\alpha=\dot c_{\pi(\alpha)}$.
	\end{proof}
	
	\begin{remark}[Standard lifts and Cohen-type diagonal groups]
		\label{rem:standard-is-diagonal}
		The standard lift $\widehat{\pi}^{\,\emptyset}$ from Definition~\ref{def:iteration-diag}(B)
		is the diagonal lift with no protected coordinates, hence it acts as $\pi$ on the Cohen base and as the
		standard $\pi$-action on every package coordinate.
		For a countable $E\subseteq\omega_1$, the subgroup of standard lifts
		\[
		\Fix(E)^\uparrow:=\{\widehat{\pi}^{\,\emptyset}:\pi\in\Fix(E)\}\le\cG_\lambda
		\]
		is contained in the characterized diagonal group
		\[
		\Delta_\lambda^\uparrow(E,\emptyset)=\{g\in\cG_\lambda:\pi^0_\lambda(g)\in\Fix(E)\},
		\]
		and in general this inclusion may be proper.
		In the limit-stage diagonal filters, both Cohen-type constraints and package-coordinate constraints are
		handled uniformly using groups of the form $\Delta_\lambda^\uparrow(E,D)$ for countable parameters
		$E$ and countable sets of coordinates $D\subseteq\mathrm{Coords}_{<\lambda}$.
	\end{remark}
	
	\begin{lemma}[Monotonicity for iteration-level diagonal groups]
		\label{lem:iteration-diag-mono}
		If $E_0 \subseteq E_1$ and $D_0 \subseteq D_1$, then
		\[
		\Delta_\lambda^\uparrow(E_1,D_1) \le \Delta_\lambda^\uparrow(E_0,D_0).
		\]
	\end{lemma}
	
	\begin{proof}
		Let $k \in \Delta_\lambda^\uparrow(E_1,D_1)$.
		By Definition~\ref{def:iteration-diag}(A):
		\begin{itemize}
			\item $\pi^0_\lambda(k) \in \Fix(E_1) \subseteq \Fix(E_0)$
			(since $E_0 \subseteq E_1$ implies $\Fix(E_1) \le \Fix(E_0)$);
			\item $k$ fixes every coordinate in $D_1$ pointwise, hence fixes every coordinate in $D_0\subseteq D_1$ pointwise.
		\end{itemize}
		Thus $k \in \Delta_\lambda^\uparrow(E_0,D_0)$.
	\end{proof}
	
	\begin{corollary}[Intersections contain a uniform diagonal subgroup]
		\label{cor:delta-up-intersection}
		Let $\{(E_i,D_i):i\in I\}$ be a family with $I$ finite or countable, and assume each $E_i$ and each
		$D_i$ is countable.
		Then
		\[
		\bigcap_{i \in I} \Delta_\lambda^\uparrow(E_i, D_i) \supseteq
		\Delta_\lambda^\uparrow\left(\bigcup_{i \in I} E_i,\ \bigcup_{i \in I} D_i\right).
		\]
		In particular, under $M\models\DC$ the unions on the right are countable, so the right-hand group is again of the
		form used to generate the diagonal stage filters.
	\end{corollary}
	
	\begin{proof}
		By Lemma~\ref{lem:iteration-diag-mono}, for each $i \in I$:
		$\Delta_\lambda^\uparrow(\bigcup_{j} E_j, \bigcup_{j} D_j) \le \Delta_\lambda^\uparrow(E_i, D_i)$.
		Hence the union-indexed group is contained in every term of the intersection.
	\end{proof}
	
	\begin{lemma}[Conjugation of iteration-level diagonal groups]
		\label{lem:iteration-diag-conj}
		Let $E\subseteq\omega_1$ and let $D\subseteq\mathrm{Coords}_{<\lambda}$.
		For $\sigma\in\cG=\Sym(\omega_1)$,
		\[
		\widehat{\sigma}^{\,\emptyset}\cdot \Delta_\lambda^\uparrow(E,D)\cdot (\widehat{\sigma}^{\,\emptyset})^{-1}
		= \Delta_\lambda^\uparrow\bigl(\sigma[E],\ \widehat{\sigma}^{\,\emptyset}\cdot D\bigr),
		\]
		where $\widehat{\sigma}^{\,\emptyset}\cdot D:=\{\widehat{\sigma}^{\,\emptyset}\cdot d:d\in D\}$.
	\end{lemma}
	
	\begin{proof}
		Let $k\in\Delta_\lambda^\uparrow(E,D)$ and set
		$h:=\widehat{\sigma}^{\,\emptyset}\,k\,(\widehat{\sigma}^{\,\emptyset})^{-1}$.
		By Lemma~\ref{lem:group-action-coords}(i), $\pi^0_\lambda$ is a homomorphism, hence
		\[
		\pi^0_\lambda(h)=\sigma\,\pi^0_\lambda(k)\,\sigma^{-1}.
		\]
		Since $\pi^0_\lambda(k)\in\Fix(E)$, the conjugate fixes $\sigma[E]$ pointwise, so
		$\pi^0_\lambda(h)\in\Fix(\sigma[E])$.
		Next, for $d\in D$ we have $k\cdot d=d$, hence
		$h\cdot(\widehat{\sigma}^{\,\emptyset}\cdot d)=\widehat{\sigma}^{\,\emptyset}\cdot d$.
		Therefore $h$ fixes every coordinate in $\widehat{\sigma}^{\,\emptyset}\cdot D$ pointwise.
		By Definition~\ref{def:iteration-diag}(A), $h\in\Delta_\lambda^\uparrow(\sigma[E],\widehat{\sigma}^{\,\emptyset}\cdot D)$.
		The reverse inclusion follows by conjugating with $(\widehat{\sigma}^{\,\emptyset})^{-1}$.
	\end{proof}
	
	\begin{remark}[Admissibility under conjugation up to enlarging $E$]
		\label{rem:admissibility-conjugation}
		Let $\lambda$ be a stage and let $(E,D)$ be globally admissible at stage $\lambda$.
		For any $k\in\cG_\lambda$, there is a countable $E^\ast\supseteq k\cdot E$ such that
		$(E^\ast,\,k\cdot D)$ is globally admissible.
		
		Consequently,
		\[
		\Delta_\lambda^\uparrow(E^\ast,\,k\cdot D)\subseteq k\,\Delta_\lambda^\uparrow(E,D)\,k^{-1}.
		\]
		Indeed, by Lemma~\ref{lem:delta-up-general-conj} we have
		$k\,\Delta_\lambda^\uparrow(E,D)\,k^{-1}=\Delta_\lambda^\uparrow(k\cdot E,k\cdot D)$, and since
		$k\cdot E\subseteq E^\ast$ monotonicity (Lemma~\ref{lem:iteration-diag-mono}) yields
		$\Delta_\lambda^\uparrow(E^\ast,k\cdot D)\le\Delta_\lambda^\uparrow(k\cdot E,k\cdot D)$.
	\end{remark}
	
	\begin{lemma}[General conjugation of $\Delta^\uparrow$-groups]
		\label{lem:delta-up-general-conj}
		Let $k \in \cG_\lambda$, let $E \subseteq \omega_1$, and let $D\subseteq \mathrm{Coords}_{<\lambda}$.
		Then
		\[
		k\,\Delta_\lambda^\uparrow(E, D)\,k^{-1} =
		\Delta_\lambda^\uparrow(k \cdot E, k \cdot D),
		\]
		where $k \cdot E := \pi^0_\lambda(k)[E]$ and $k \cdot D := \{k \cdot d : d \in D\}$ (using the
		$\cG_\lambda$-action on coordinates from Lemma~\ref{lem:group-action-coords}).
	\end{lemma}
	
	\begin{proof}
		Fix $k\in\cG_\lambda$ and put $\sigma:=\pi^0_\lambda(k)\in\cG=\Sym(\omega_1)$.
		Let $h\in\Delta^\uparrow_\lambda(E,D)$, and put $\tau:=\pi^0_\lambda(h)\in\Fix(E)$.
		
		\emph{Stage-$0$ part.}
		By Lemma~\ref{lem:group-action-coords}(i), $\pi^0_\lambda$ is a homomorphism, hence
		\[
		\pi^0_\lambda(khk^{-1})=\sigma\tau\sigma^{-1}\in \Fix(\sigma[E])=\Fix(\pi^0_\lambda(k)[E]).
		\]
		
		\emph{Coordinate part.}
		Let $d\in D$. By definition of $\Delta^\uparrow_\lambda(E,D)$, $h$ fixes $d$ pointwise.
		Using the coordinate action of $\cG_\lambda$ on $\mathrm{Coords}_{<\lambda}$ (Lemma~\ref{lem:group-action-coords}),
		\[
		(khk^{-1})\cdot(k\cdot d)
		= k\cdot\bigl(h\cdot(k^{-1}\cdot(k\cdot d))\bigr)
		= k\cdot(h\cdot d)
		= k\cdot d,
		\]
		so $khk^{-1}$ fixes every coordinate in $k\cdot D$ pointwise.
		Thus $khk^{-1}\in\Delta^\uparrow_\lambda(k\cdot E,k\cdot D)$, proving
		$k\,\Delta_\lambda^\uparrow(E,D)\,k^{-1}\subseteq \Delta^\uparrow_\lambda(k\cdot E,k\cdot D)$.
		The reverse inclusion follows by conjugating with $k^{-1}$.
	\end{proof}
	
	\begin{corollary}[Admissibility preserved under unions]
		\label{cor:admissibility-union}
		If $(E_1, D_1)$ and $(E_2, D_2)$ are globally admissible, then 
		$(E_1 \cup E_2, D_1 \cup D_2)$ is globally admissible.
	\end{corollary}
	
	\begin{proof}
		For $d \in D_1$: $\Fix(E_1 \cup E_2) \le \Fix(E_1) \le \Stab_{\cG}(d)$.
		For $d \in D_2$: $\Fix(E_1 \cup E_2) \le \Fix(E_2) \le \Stab_{\cG}(d)$.
	\end{proof}
	
	\begin{lemma}[Admissible hull for countable coordinate sets]
		\label{lem:admissible-hull}
		Let $D$ be a countable set of package coordinates (at stages $<\lambda$) and
		let $E \subseteq \omega_1$ be countable.
		Then there exists a countable set
		$E^\ast \subseteq \omega_1$ with $E \subseteq E^\ast$ such that $(E^\ast, D)$ is
		globally admissible.
		Consequently, $\Delta_\lambda^\uparrow(E^\ast, D) \le \Delta_\lambda^\uparrow(E, D)$.
	\end{lemma}
	
	\begin{proof}
		For each $d \in D$, Lemma~\ref{lem:coord-countable-support} provides a countable 
		$E_d \subseteq \omega_1$ such that $\Fix(E_d) \le \Stab_{\cG}(d)$.
		Let $E^\ast := E \cup \bigcup_{d \in D} E_d$, which is a countable union of 
		countable sets, hence countable.
		For each $d \in D$, we have $\Fix(E^\ast) \le \Fix(E_d) \le \Stab_{\cG}(d)$, so 
		$(E^\ast, D)$ is globally admissible.
		
		The final statement follows from monotonicity: $E \subseteq E^\ast$ implies 
		$\Delta_\lambda^\uparrow(E^\ast, D) \le \Delta_\lambda^\uparrow(E, D)$ 
		by Lemma~\ref{lem:iteration-diag-mono}.
	\end{proof}
	
	\begin{lemma}[Projection of $\Delta^\uparrow$-groups]
		\label{lem:delta-up-projection}
		Let $\beta < \lambda$, let $E\subseteq\omega_1$, and let $D \subseteq \mathrm{Coords}_\beta$
		(i.e., $D$ only uses coordinates introduced below stage $\beta$).
		Then
		\[
		\pi^\beta_\lambda\big[\Delta_\lambda^\uparrow(E, D)\big] \subseteq
		\Delta_\beta^\uparrow(E, D).
		\]
	\end{lemma}
	
	\begin{proof}
		Let $k \in \Delta_\lambda^\uparrow(E, D)$. Then:
		\begin{itemize}
			\item $\pi^0_\lambda(k) \in \Fix(E)$, and by Definition~\ref{def:stage0-projection} (using coherence from
			Remark~\ref{rem:projection-coherence}) we have $\pi^0_\beta \circ \pi^\beta_\lambda = \pi^0_\lambda$.
			Hence
			\[
			\pi^0_\beta(\pi^\beta_\lambda(k)) = \pi^0_\lambda(k) \in \Fix(E).
			\]
			\item $k$ fixes every coordinate in $D$ pointwise.
			Since $D \subseteq \mathrm{Coords}_\beta$,
			these coordinates are already present at stage $\beta$, and the projection $\pi^\beta_\lambda(k)$
			acts on $\mathrm{Coords}_\beta$ as the restriction of $k$ (cf.\ Lemma~\ref{lem:group-action-coords}(ii)),
			so $\pi^\beta_\lambda(k)$ fixes every coordinate in $D$ pointwise.
		\end{itemize}
		Hence $\pi^\beta_\lambda(k) \in \Delta_\beta^\uparrow(E, D)$ by Definition~\ref{def:iteration-diag}(A).
	\end{proof}
	
	\begin{lemma}[Pullback contains $\Delta^\uparrow$]
		\label{lem:pullback-delta}
		Let $\beta < \lambda$, let $E\subseteq\omega_1$, and let $D\subseteq \mathrm{Coords}_\beta$.
		If $H \in \tilde{\cF}^*_\beta$ contains $\Delta_\beta^\uparrow(E, D)$, then
		$(\pi^\beta_\lambda)^{-1}[H]$ contains $\Delta_\lambda^\uparrow(E, D)$.
	\end{lemma}
	
	\begin{proof}
		Let $k \in \Delta_\lambda^\uparrow(E, D)$.
		By Lemma~\ref{lem:delta-up-projection},
		$\pi^\beta_\lambda(k) \in \Delta_\beta^\uparrow(E, D) \subseteq H$.
		Hence $k \in (\pi^\beta_\lambda)^{-1}[H]$.
	\end{proof}
	
	\begin{lemma}[Tail extension by identity]
		\label{lem:tail-extension-identity}
		Let \(\beta\le\lambda\) be stages of the iteration, where \(\lambda\) may also be
		\(\Ord\) in the background class metatheory.
		For every \(k\in\cG_\beta\), there exists
		\(\hat k\in\cG_\lambda\) such that:
		\begin{enumerate}[label=(\roman*)]
			\item \(\pi^\beta_\lambda(\hat k)=k\);
			\item for every successor stage \(\gamma+1\) with \(\beta\le\gamma<\lambda\), if
			\(\hat k_\gamma:=\pi^\gamma_\lambda(\hat k)\), then
			\[
			\hat k_{\gamma+1}=\mathrm{spl}(\hat k_\gamma);
			\]
			\item equivalently, \(\hat k\) acts as \(k\) on the \(\bbP_\beta\)-part and trivially
			on every iterand added at stages in \([\beta,\lambda)\).
		\end{enumerate}
	\end{lemma}
	
	\begin{proof}
		If \(\beta=\lambda\), take \(\hat k:=k\). Assume \(\beta<\lambda\).
		
		We build by transfinite recursion a coherent thread
		\(\langle k_\gamma:\beta\le\gamma\le\lambda\rangle\) such that \(k_\beta=k\), each
		\(k_\gamma\in\cG_\gamma\), and \(\pi^\delta_\gamma(k_\gamma)=k_\delta\) whenever
		\(\beta\le\delta<\gamma\le\lambda\).
		
		\smallskip
		\noindent\textbf{Base stage.}
		Set \(k_\beta:=k\).
		
		\smallskip
		\noindent\textbf{Successor stage \(\gamma+1\).}
		Given \(k_\gamma\in\cG_\gamma\), define
		\[
		k_{\gamma+1}:=\mathrm{spl}(k_\gamma)\in\cG_{\gamma+1},
		\]
		where \(\mathrm{spl}\) is the splitting section from
		Lemma~\ref{lem:successor-splitting}. Then
		\[
		\pi^\gamma_{\gamma+1}(k_{\gamma+1})
		=
		\pi^\gamma_{\gamma+1}(\mathrm{spl}(k_\gamma))
		=
		k_\gamma,
		\]
		so coherence is preserved.
		By construction, \(k_{\gamma+1}\) acts as \(k_\gamma\) on the
		\(\bbP_\gamma\)-part and trivially on the new iterand at stage \(\gamma\).
		
		\smallskip
		\noindent\textbf{Limit stage \(\eta\le\lambda\).}
		Assume \(\langle k_\gamma:\beta\le\gamma<\eta\rangle\) has been constructed coherently.
		By the limit-stage group definition used in the iteration, coherent threads determine
		unique elements of \(\cG_\eta\).
		Thus there is a unique \(k_\eta\in\cG_\eta\) such that
		\[
		\pi^\gamma_\eta(k_\eta)=k_\gamma\qquad(\beta\le\gamma<\eta).
		\]
		If \(\eta=\Ord\), this is carried out in the same background metatheory used for the
		class-length bookkeeping (Remark~\ref{rem:metatheory}).
		
		Now set \(\hat k:=k_\lambda\). Then (i) holds by construction. Clause (ii) is exactly the
		successor clause of the recursion.
		Clause (iii) is the equivalent operational description:
		at each successor step we extend by the splitting section, which acts trivially on the new
		iterand and preserves the previously constructed action.
	\end{proof}
	
	\begin{remark}[Set stages versus the class stage]
		\label{rem:set-vs-class-stages}
		In statements where \(\lambda\) may be \(\Ord\), the set-stage case uses the bounded-stage
		projection/coherence system imported in Remark~\ref{rem:iteration-api}(a),(b), while the case
		\(\lambda=\Ord\) uses the class objects \(\bbP_{\Ord}\), \(\cG_{\Ord}\), and
		\(\pi^\beta_{\Ord}\) defined by the explicit recursion of
		Subsection~\ref{subsec:iteration-api-level} in the background metatheory of
		Remark~\ref{rem:metatheory}.
		No appeal to a separate class-length theorem is made.
	\end{remark}
	
	\begin{corollary}[Identity extension for \(\Delta^\uparrow\)-groups]
		\label{cor:delta-up-extension}
		Let \(\beta\le\lambda\), let \(E\subseteq\omega_1\), and let
		\(D\subseteq\mathrm{Coords}_\beta\).
		If \(k\in\Delta_\beta^\uparrow(E,D)\), then there exists
		\(\hat k\in\Delta_\lambda^\uparrow(E,D)\) such that
		\[
		\pi^\beta_\lambda(\hat k)=k.
		\]
	\end{corollary}
	
	\begin{proof}
		Let $\hat k\in\cG_\lambda$ be the identity-tail extension of $k$ given by Lemma~\ref{lem:tail-extension-identity}.
		
		\emph{Case 1: $\lambda<\Ord$.} For the stage-$0$ part, by Definition~\ref{def:stage0-projection} and projection coherence,
		\[
		\pi^0_\lambda(\hat k)=\pi^0_\beta(\pi^\beta_\lambda(\hat k))=\pi^0_\beta(k)\in\Fix(E).
		\]
		For $d\in D$, $k$ fixes $d$ at stage $\beta$, and the splitting section preserves that action at every successor above $\beta$; by induction $\hat k$ fixes $d$. Hence $\hat k\in\Delta_\lambda^\uparrow(E,D)$.
		
		\emph{Case 2: $\lambda=\Ord$.} Put $\alpha_0:=\alpha_0(E,D)$ as in Definition~\ref{def:limit-stage-group}. Let $\beta_0:=\max(\beta,\alpha_0+1)$. Define $\langle\hat k_\gamma\rangle_{\gamma\ge\beta_0}$ by recursion: set $\hat k_{\beta_0}$ to be the identity-tail lift of $k$ to stage $\beta_0$ given by Lemma~\ref{lem:tail-extension-identity} (if $\beta_0=\beta$ this is $k$). For successor $\gamma+1>\beta_0$ let $\hat k_{\gamma+1}$ be the lift of $\hat k_\gamma$ supplied by the same lemma; for limit $\lambda>\beta_0$ let $\hat k_\lambda$ be the unique element of $\cG_\lambda$ with $\pi^\gamma_\lambda(\hat k_\lambda)=\hat k_\gamma$ for all $\gamma<\lambda$ (the thread $\langle\hat k_\gamma\rangle_{\gamma<\lambda}$). By Lemma~\ref{lem:tail-extension-identity}(i) each lift acts as the identity on coordinates outside $D$ and fixes $E$, hence $\hat k_\gamma\in\Delta^\uparrow_\gamma(E,D)$ for all $\gamma\ge\beta_0$, since $D\subseteq\mathrm{Coords}_\beta\subseteq\mathrm{Coords}_\gamma$ for $\beta\le\gamma$ by Remark~\ref{rem:iteration-api}(a). For $\gamma<\beta_0$ define $\hat k_\gamma:=\pi^\gamma_{\beta_0}(\hat k_{\beta_0})$; by Lemma~\ref{lem:tail-extension-identity}, $\pi^\beta_{\beta_0}(\hat k_{\beta_0})=k$, so in particular $\hat k_\beta=k$ as required. We verify coherence $\pi^{\gamma_1}_{\gamma_2}(\hat k_{\gamma_2})=\hat k_{\gamma_1}$ for all $\gamma_1<\gamma_2$ in three cases:
		
		\emph{Case (i) $\gamma_1<\gamma_2\le\beta_0$:} $\hat k_{\gamma_2}=\pi^{\gamma_2}_{\beta_0}(\hat k_{\beta_0})$ and $\hat k_{\gamma_1}=\pi^{\gamma_1}_{\beta_0}(\hat k_{\beta_0})$ by definition, so $\pi^{\gamma_1}_{\gamma_2}(\hat k_{\gamma_2})=\pi^{\gamma_1}_{\gamma_2}(\pi^{\gamma_2}_{\beta_0}(\hat k_{\beta_0}))=\pi^{\gamma_1}_{\beta_0}(\hat k_{\beta_0})=\hat k_{\gamma_1}$ by Remark~\ref{rem:iteration-api}(a).
		
		\emph{Case (ii) $\beta_0\le\gamma_1<\gamma_2$:} If $\gamma_2=\gamma_1+1$, then $\hat k_{\gamma_2}$ is the lift of $\hat k_{\gamma_1}$ via Lemma~\ref{lem:tail-extension-identity}, so $\pi^{\gamma_1}_{\gamma_2}(\hat k_{\gamma_2})=\hat k_{\gamma_1}$ by that lemma. If $\gamma_2$ is limit, then $\hat k_{\gamma_2}$ is defined as the unique element of $\cG_{\gamma_2}$ with $\pi^{\gamma}_{\gamma_2}(\hat k_{\gamma_2})=\hat k_{\gamma}$ for all $\gamma<\gamma_2$, so in particular $\pi^{\gamma_1}_{\gamma_2}(\hat k_{\gamma_2})=\hat k_{\gamma_1}$.
		
		\emph{Case (iii) $\gamma_1<\beta_0<\gamma_2$:} $\pi^{\gamma_1}_{\gamma_2}(\hat k_{\gamma_2})=\pi^{\gamma_1}_{\beta_0}(\pi^{\beta_0}_{\gamma_2}(\hat k_{\gamma_2}))=\pi^{\gamma_1}_{\beta_0}(\hat k_{\beta_0})=\hat k_{\gamma_1}$ by cases (i) and (ii).
		
		Hence $\hat k:=\langle\hat k_\gamma\rangle_{\gamma<\Ord}$ is a coherent thread, so $\hat k\in\cG_{\Ord}$ by the inverse limit definition. Moreover, with $(E,D)$ as witness, we have $\hat k_\gamma\in\Delta^\uparrow_\gamma(E,D)$ for all $\gamma\ge\beta_0>\alpha_0(E,D)$. By Definition~\ref{def:limit-stage-group} and Lemma~\ref{lem:GOrd-definable}, this implies $\hat k\in\cG_{\Ord}^{bd}$. Since $\hat k_\gamma\in\Delta^\uparrow_\gamma(E,D)$ for all $\gamma>\alpha_0$, we have $\hat k\in\Delta^\uparrow_{\Ord}(E,D)$ by Definition~\ref{def:limit-stage-group}.
	\end{proof}
	
	\begin{lemma}[Action on bounded names factors through projection]
		\label{lem:projection-action-bounded-name}
		Let \(\beta\le\lambda\) be stages, let \(\tau\) be a \(\bbP_\beta\)-name, and let
		\(k\in\cG_\lambda\).
		Viewing \(\tau\) as a \(\bbP_\lambda\)-name by canonical inclusion,
		one has
		\[
		k\cdot\tau \;=\; \pi^\beta_\lambda(k)\cdot\tau.
		\]
	\end{lemma}
	
	\begin{proof}
		We argue by induction on the rank of \(\tau\).
		Let \(\langle \sigma,p\rangle\in\tau\). Since \(\tau\) is a \(\bbP_\beta\)-name, we have
		\(p\in\bbP_\beta\).
		The projection \(\pi^\beta_\lambda(k)\) is, by definition, the
		restriction of \(k\) to the \(\bbP_\beta\)-part, so
		\[
		k\cdot p \;=\; \pi^\beta_\lambda(k)\cdot p.
		\]
		By the induction hypothesis,
		\[
		k\cdot\sigma \;=\; \pi^\beta_\lambda(k)\cdot\sigma.
		\]
		Therefore the image of each pair \(\langle \sigma,p\rangle\) under \(k\) is exactly the
		same as its image under \(\pi^\beta_\lambda(k)\), and hence
		\(k\cdot\tau=\pi^\beta_\lambda(k)\cdot\tau\).
	\end{proof}
	
	\begin{lemma}[Iteration-level diagonal lifts lie in $\cG_\lambda$]
		\label{lem:diag-lift-in-group}
		Let $\lambda$ be a stage of the iteration.
		Let $(E, D)$ be globally admissible at 
		stage $\lambda$ and let $\pi \in \Fix(E)$.
		Then the iteration-level diagonal lift 
		$\widehat{\pi}^{\,D} \in \Aut(\bbP_\lambda)$ (from Definition~\ref{def:iteration-diag}) 
		belongs to $\cG_\lambda$.
	\end{lemma}
	
	\begin{proof}
		For each $\beta\le\lambda$ let $D_\beta:=D\cap\mathrm{Coords}_\beta$.
		By Definition~\ref{def:global-admissible}, $\Fix(E)\le\Stab_{\cG}(d)$ for every $d\in D$, hence
		$\pi\cdot d=d$ for all $d\in D$.
		In particular, for each successor $\beta+1\le\lambda$ the set
		\[
		D^{\mathrm{new}}_{\beta+1}:=D_{\beta+1}\setminus D_\beta
		\]
		is pointwise fixed by $\pi$, so $\kappa_{\beta+1}(\pi,D^{\mathrm{new}}_{\beta+1})$ is well-defined
		(Definition~\ref{def:succ-diag-kernel}).
		
		We build by transfinite recursion a coherent sequence $\langle k_\beta:\beta\le\lambda\rangle$
		with $k_\beta\in\cG_\beta$ such that $k_\beta=\widehat{\pi}^{\,D_\beta}$ as automorphisms of $\bbP_\beta$,
		and $\pi^\gamma_\beta(k_\beta)=k_\gamma$ for all $\gamma<\beta\le\lambda$.
		
		\smallskip
		\noindent\textbf{Stage $0$.}
		Set $k_0:=\pi\in\cG_0=\cG$. Since $D_0=\varnothing$, this agrees with
		$\widehat{\pi}^{\,D_0}=\widehat{\pi}^{\,\emptyset}$ on $\bbP_0$.
		
		\smallskip
		\noindent\textbf{Successor $\beta+1$.}
		Assume $k_\beta\in\cG_\beta$ with $k_\beta=\widehat{\pi}^{\,D_\beta}$.
		Define
		\[
		k_{\beta+1}:=\mathrm{spl}(k_\beta)\cdot \kappa_{\beta+1}\!\left(\pi,\,D^{\mathrm{new}}_{\beta+1}\right),
		\]
		where $\mathrm{spl}$ is the splitting section from Lemma~\ref{lem:successor-splitting}.
		By Lemma~\ref{lem:successor-splitting}, $\mathrm{spl}(k_\beta)\in\cG_{\beta+1}$, and by
		Lemma~\ref{lem:succ-diag-kernel-in-kernel}(i) we have
		$\kappa_{\beta+1}(\pi,D^{\mathrm{new}}_{\beta+1})\in\cG_{\beta+1}$, hence $k_{\beta+1}\in\cG_{\beta+1}$.
		
		Moreover, $k_{\beta+1}$ agrees with $\widehat{\pi}^{\,D_{\beta+1}}$ on $\bbP_{\beta+1}$:
		on the $\bbP_\beta$-part, $\mathrm{spl}(k_\beta)$ acts as $k_\beta=\widehat{\pi}^{\,D_\beta}$ and the
		kernel element acts trivially;
		on the new coordinates at stage $\beta+1$,
		$\mathrm{spl}(k_\beta)$ acts trivially and $\kappa_{\beta+1}(\pi,D^{\mathrm{new}}_{\beta+1})$ acts as
		$\pi$ outside $D^{\mathrm{new}}_{\beta+1}$ and as the identity on $D^{\mathrm{new}}_{\beta+1}$,
		which is exactly the defining clause of $\widehat{\pi}^{\,D_{\beta+1}}$
		(Definition~\ref{def:iteration-diag}).
		
		Finally,
		\[
		\pi^\beta_{\beta+1}(k_{\beta+1})
		=\pi^\beta_{\beta+1}(\mathrm{spl}(k_\beta))\cdot \pi^\beta_{\beta+1}(\kappa_{\beta+1}(\ldots))
		= k_\beta\cdot \mathrm{id}=k_\beta,
		\]
		using $\pi^\beta_{\beta+1}\circ\mathrm{spl}=\mathrm{id}$ (Lemma~\ref{lem:successor-splitting}) and
		$\kappa_{\beta+1}(\ldots)\in K_{\beta+1}=\ker(\pi^\beta_{\beta+1})$
		(Lemma~\ref{lem:succ-diag-kernel-in-kernel}(ii)).
		
		\smallskip
		\noindent\textbf{Limit $\eta\le\lambda$.}
		Let $\eta$ be limit and assume $\langle k_\beta:\beta<\eta\rangle$ has been constructed coherently.
		By the limit-stage group clause in the countable-support iteration setup
		(Definition~\ref{def:limit-stage-group}),
		elements of $\cG_\eta$ 
		are exactly the coherent threads $\langle h_\beta \rangle_{\beta<\eta}$ with 
		$\pi^\gamma_\beta(h_\beta)=h_\gamma$ for $\gamma<\beta<\eta$.
		Since 
		$\langle k_\beta \rangle_{\beta<\eta}$ is such a coherent thread by construction, 
		it defines a unique element $k_\eta\in\cG_\eta$ with $\pi^\beta_\eta(k_\eta)=k_\beta$ for all $\beta<\eta$.
		Since $\widehat{\pi}^{\,D_\eta}$ has the same restrictions to each $\bbP_\beta$ for $\beta<\eta$,
		we have $k_\eta=\widehat{\pi}^{\,D_\eta}$.
		
		\smallskip
		\noindent\textbf{Conclusion.}
		At $\beta=\lambda$ we obtain $k_\lambda=\widehat{\pi}^{\,D_\lambda}\in\cG_\lambda$.
		Because $D\subseteq \mathrm{Coords}_{<\lambda}$ (Definition~\ref{def:global-admissible}(ii)),
		we have $D_\lambda=D$, hence $\widehat{\pi}^{\,D}\in\cG_\lambda$.
	\end{proof}
	
	\begin{lemma}[Diagonal lifts lie in $\Delta^\uparrow$]
		\label{lem:diag-lift-in-delta}
		Let $\lambda$ be a stage, let $E\subseteq\omega_1$, and let $D\subseteq \mathrm{Coords}_{<\lambda}$.
		Let $\pi \in \Fix(E)$ and let $D' \supseteq D$ be such that the diagonal lift
		$\widehat{\pi}^{\,D'} \in \cG_\lambda$ exists.
		Then
		\[
		\widehat{\pi}^{\,D'} \in \Delta_\lambda^\uparrow(E, D).
		\]
	\end{lemma}
	
	\begin{proof}
		Since $\widehat{\pi}^{\,D'}\in\cG_\lambda$ and $\widehat{\pi}^{\,D'}$ restricts to $\pi$ on the Cohen base
		(Lemma~\ref{lem:iteration-diag-auto}(i)), we have
		$\pi^0_\lambda(\widehat{\pi}^{\,D'})=\pi\in\Fix(E)$ (e.g.\ by Corollary~\ref{cor:diag-lift-swaps-cohen}).
		Moreover, by Lemma~\ref{lem:iteration-diag-auto}(ii) the automorphism $\widehat{\pi}^{\,D'}$ fixes every
		coordinate in $D'$ pointwise, hence fixes every coordinate in $D\subseteq D'$ pointwise.
		Therefore $\widehat{\pi}^{\,D'}\in\Delta_\lambda^\uparrow(E,D)$ by Definition~\ref{def:iteration-diag}(A).
	\end{proof}
	
	\begin{lemma}[Package generic sections are hereditarily symmetric via diagonal cancellation]
		\label{lem:sf-HS-diag}
		Let $f$ have countable support witnessed by $E\in[\omega_1]^{\le\omega}$ (so $\Fix(E)\le\Stab(f)$),
		and let $\dot{s}_f$ be the canonical $\bbS_{[f]}$-name for the $f$-coordinate generic section.
		Then
		\[
		\Delta^\uparrow(E,\{f\})\ \le\ \Stab_{\cG_{[f]}}(\dot{s}_f),
		\]
		and hence $\dot{s}_f$ is hereditarily symmetric for the orbit package system
		$(\bbS_{[f]},\cG_{[f]},\cF_{[f]})$.
	\end{lemma}
	
	\begin{proof}
		Fix $E\in[\omega_1]^{\le\omega}$ with $\Fix(E)\le\Stab(f)$.
		Let $(E',D')$ be any admissible pair in the sense of Definition~\ref{def:diag-cancel} such that
		$E\subseteq E'$ and $f\in D'$.
		Let $\pi\in\Fix(E')$ and consider the diagonal-cancellation automorphism
		$\widehat{\pi}^{\,D'}\in\Delta(E',D')$ (Definition~\ref{def:diag-cancel}).
		Because $f\in D'$, this automorphism fixes the $f$-coordinate pointwise, i.e.\ for every
		$\vec p\in\bbS_{[f]}$,
		\[
		(\widehat{\pi}^{\,D'}\cdot\vec p)(f)=\vec p(f).
		\]
		Since $\dot{s}_f$ is the canonical name definable purely from the $f$-coordinate generic section,
		we have $\widehat{\pi}^{\,D'}\cdot\dot{s}_f=\dot{s}_f$.  Thus
		$\Delta(E',D')\le \Stab_{\cG_{[f]}}(\dot{s}_f)$.
		
		By Definition~\ref{def:orbit-stage-filter-diag}, $\Delta^\uparrow(E,\{f\})$ is generated by all such
		$\Delta(E',D')$ with $E\subseteq E'$ and $f\in D'$, so
		$\Delta^\uparrow(E,\{f\})\le \Stab_{\cG_{[f]}}(\dot{s}_f)$.
		Therefore $\Stab_{\cG_{[f]}}(\dot{s}_f)\in\cF_{[f]}$ by Definition~\ref{def:orbit-stage-filter-diag},
		and $\dot{s}_f$ is hereditarily symmetric.
	\end{proof}
	
	\begin{lemma}[Diagonal cancellation makes right-inverse generics hereditarily symmetric]
		\label{lem:rg-HS-diag}
		Let $g$ be a package with countable support witnessed by $E\in[\omega_1]^{\le\omega}$, i.e.\ $\Fix(E)\le\Stab(g)$.
		Let $\dot{r}_g$ be the canonical $\bbR_{[g]}$-name for the $g$-coordinate generic right inverse.
		Then
		\[
		\Delta^\uparrow(E,\{g\})\ \le\ \Stab_{\cG_{[g]}}(\dot{r}_g),
		\]
		and hence $\dot{r}_g$ is hereditarily symmetric for $(\bbR_{[g]},\cG_{[g]},\cF_{[g]})$.
	\end{lemma}
	
	\begin{proof}
		Fix $E\in[\omega_1]^{\le\omega}$ with $\Fix(E)\le\Stab(g)$.
		Let $(E',D')$ be any admissible pair (Definition~\ref{def:diag-cancel}) such that
		$E\subseteq E'$ and $g\in D'$.
		For $\pi\in\Fix(E')$, the diagonal-cancellation map
		$\widehat{\pi}^{\,D'}\in\Delta(E',D')$ fixes the $g$-coordinate pointwise (since $g\in D'$).
		Because $\dot{r}_g$ is the canonical name definable purely from the $g$-coordinate generic object,
		we have $\widehat{\pi}^{\,D'}\cdot\dot{r}_g=\dot{r}_g$.  Hence
		$\Delta(E',D')\le \Stab_{\cG_{[g]}}(\dot{r}_g)$.
		
		By Definition~\ref{def:orbit-stage-filter-diag}, $\Delta^\uparrow(E,\{g\})$ is generated by all such
		$\Delta(E',D')$ with $E\subseteq E'$ and $g\in D'$, so
		$\Delta^\uparrow(E,\{g\})\le \Stab_{\cG_{[g]}}(\dot{r}_g)$.
		Therefore $\Stab_{\cG_{[g]}}(\dot{r}_g)\in\cF_{[g]}$ (Definition~\ref{def:orbit-stage-filter-diag}),
		and $\dot{r}_g$ is hereditarily symmetric.
	\end{proof}
	
	\begin{lemma}[Successor-stage diagonal cancellation yields \(\HS\)-lifts]
		\label{lem:successor-HS-lift}
		Let \(\beta\) be an ordinal and consider the successor stage \(\beta+1\) with
		\(\bbP_{\beta+1}\cong \bbP_\beta*\dot{\bbQ}_\beta\).
		Let \(d\) be a package coordinate of \(\dot{\bbQ}_\beta\) (e.g.\ a \(\bbQ\)-package coordinate
		or a \(\bbR\)-package coordinate), and let \(\dot{\tau}\) be a canonical \(\dot{\bbQ}_\beta\)-name
		definable purely from the \(d\)-coordinate generic object (e.g.\ \(\dot{s}_{f_\beta}\) or \(\dot{r}_{g_\beta}\)).
		Let \(\dot{\tau}^*\) denote the corresponding canonical \(\bbP_{\beta+1}\)-name under this two-step presentation.
		
		Assume \(d\) has countable support witnessed by \(E\in[\omega_1]^{\le\omega}\) in \(M_\beta\),
		so \(\Fix(E)\le\Stab_{\cG}(d)\) in \(M_\beta\).
		Then \(\dot{\tau}^*\in\HS_{\beta+1}\); in particular,
		\(\dot{\tau}^{*\,G_{\beta+1}}\in M_{\beta+1}\subseteq\cM\).
	\end{lemma}
	
	\begin{proof}
		From \(\Fix(E)\le\Stab_{\cG}(d)\) we get that \((E,\{d\})\) is globally admissible at stage \(\beta+1\)
		(Definition~\ref{def:global-admissible}).
		Hence
		\(\Delta_{\beta+1}^\uparrow(E,\{d\})\) is one of the designated \(\Delta^\uparrow\)-generators of
		\(\tilde{\cF}^*_{\beta+1}\) (Definition~\ref{def:modified-limit-filter}).
		
		Let \(k\in \Delta_{\beta+1}^\uparrow(E,\{d\})\).
		By Definition~\ref{def:iteration-diag},
		\(k\) fixes the coordinate \(d\) pointwise (i.e.\ its induced action on the \(d\)-factor forcing is the identity).
		Since \(\dot{\tau}^*\) is a canonical name definable purely from the \(d\)-coordinate generic object,
		it follows that \(k\cdot\dot{\tau}^*=\dot{\tau}^*\).  Therefore
		\[
		\Delta_{\beta+1}^\uparrow(E,\{d\})\ \le\ \Stab_{\cG_{\beta+1}}(\dot{\tau}^*).
		\]
		By upward closure of \(\tilde{\cF}^*_{\beta+1}\), we get
		\(\Stab_{\cG_{\beta+1}}(\dot{\tau}^*)\in\tilde{\cF}^*_{\beta+1}\), hence \(\dot{\tau}^*\in\HS_{\beta+1}\).
		The interpretation \(\dot{\tau}^{*\,G_{\beta+1}}\in M_{\beta+1}\subseteq\cM\) follows.
	\end{proof}
	
	\subsection{The modified limit filter}
	\label{subsec:limit-filter}
	
	We now define the filter $\tilde{\cF}^*_\lambda$ that governs hereditary symmetry 
	at limit stages of the iteration.
	
	\begin{definition}[Modified stage filters]
		\label{def:modified-limit-filter}
		Define by recursion on stages $\alpha$ a normal filter $\tilde{\cF}^*_\alpha$ on $\cG_\alpha$ as follows.
		
		\textbf{Stage $0$:} $\tilde{\cF}^*_0$ is the normal filter on $\cG_0=\cG=\Sym(\omega_1)$
		generated by $\{\Fix(E):E\in[\omega_1]^{\le\omega}\}$.
		(Note that $\Fix(E)=\Delta^\uparrow_0(E,\varnothing)$, and $(E,\varnothing)$ is trivially globally admissible.)
		
		\textbf{Successor stage $\alpha+1$:} $\tilde{\cF}^*_{\alpha+1}$ is the normal filter on $\cG_{\alpha+1}$
		generated by:
		\begin{enumerate}[label=(\alph*)]
			\item \emph{Pullbacks:} $(\pi^\alpha_{\alpha+1})^{-1}[H]$ for $H\in\tilde{\cF}^*_\alpha$;
			\item \emph{$\Delta^\uparrow$-subgroups:} $\Delta_{\alpha+1}^\uparrow(E,D)$ for countable
			$(E,D)$ \textbf{globally admissible at stage $\alpha+1$}.
		\end{enumerate}
		
		\textbf{Limit stage $\lambda$:} $\tilde{\cF}^*_\lambda$ is the normal filter on $\cG_\lambda$
		generated by:
		\begin{enumerate}[label=(\alph*)]
			\item \emph{Pullbacks:} $(\pi^\beta_\lambda)^{-1}[H]$ for $\beta<\lambda$ and $H\in\tilde{\cF}^*_\beta$;
			\item \emph{$\Delta^\uparrow$-subgroups:} $\Delta_\lambda^\uparrow(E,D)$ for countable
			$(E,D)$ \textbf{globally admissible at stage $\lambda$}.
		\end{enumerate}
	\end{definition}
	
	\begin{remark}[Admissibility in the filter definition]
		\label{rem:admissibility-filter}
		We restrict the $\Delta^\uparrow$-generator clause in Definition~\ref{def:modified-limit-filter}
		to \emph{globally admissible} pairs so that later arguments can uniformly recover
		$\Delta^\uparrow$-subgroups with admissible parameters (notably Lemma~\ref{lem:core-lemma},
		used in Lemma~\ref{lem:filter-properties}(iii),(iv)).
		
		Normality interacts with admissibility as follows.  If $(E,D)$ is globally admissible at stage $\lambda$
		and $k\in\cG_\lambda$, then Lemma~\ref{lem:delta-up-general-conj} gives
		\[
		k\,\Delta_\lambda^\uparrow(E,D)\,k^{-1}=\Delta_\lambda^\uparrow(k\cdot E,\ k\cdot D),
		\]
		where $k\cdot E=\pi^0_\lambda(k)[E]$ and $k\cdot D$ is the induced action on coordinates
		(Lemma~\ref{lem:group-action-coords}).
		The pair $(k\cdot E,\ k\cdot D)$ need not be globally admissible, but by
		Lemma~\ref{lem:admissible-hull} there is a countable $E^\ast\supseteq k\cdot E$ such that
		$(E^\ast,\ k\cdot D)$ is globally admissible and
		\[
		\Delta_\lambda^\uparrow(E^\ast,\ k\cdot D)\ \le\ \Delta_\lambda^\uparrow(k\cdot E,\ k\cdot D).
		\]
		Together with Corollary~\ref{cor:admissibility-union} (and its countable variant
		Corollary~\ref{cor:admissibility-countable-union}), this ensures that the normal-filter generation
		process always contains $\Delta^\uparrow$-groups with globally admissible parameters.
	\end{remark}
	
	\begin{lemma}[Every filter element contains a $\Delta^\uparrow$-group]
		\label{lem:core-lemma}
		For each stage $\lambda$ and every $H\in\tilde{\cF}^*_\lambda$, there exist countable
		$E\subseteq\omega_1$ and a countable $D\subseteq\mathrm{Coords}_{<\lambda}$ such that
		$(E,D)$ is globally admissible at stage $\lambda$ and
		\[
		\Delta_\lambda^\uparrow(E,D)\ \le\ H.
		\]
	\end{lemma}
	
	\begin{proof}
		We argue by induction on the stage $\lambda$. The induction proceeds through all $\lambda\le\Ord$; at $\lambda=\Ord$ the filter $\tilde{\cF}^*_{\Ord}$ is generated by the groups $\Delta_{\Ord}^\uparrow(E,D)$ (existing as classes by Lemma~\ref{lem:GOrd-definable}, as in Definition~\ref{def:modified-limit-filter}) together with pullbacks from earlier stages, so the same case analysis applies.
		
		\emph{Stage $0$.}
		By Definition~\ref{def:modified-limit-filter}, $\tilde{\cF}^*_0$ is the normal filter on
		$\cG_0=\Sym(\omega_1)$ generated by $\{\Fix(E):E\in[\omega_1]^{\le\omega}\}$.
		Every conjugate of $\Fix(E)$ is again of the form $\Fix(E')$ for some countable $E'$, and
		finite intersections satisfy $\Fix(E_1)\cap\Fix(E_2)\supseteq \Fix(E_1\cup E_2)$.
		Hence every $H\in\tilde{\cF}^*_0$ contains $\Fix(E)=\Delta^\uparrow_0(E,\varnothing)$ for some countable $E$,
		and $(E,\varnothing)$ is globally admissible. This proves the claim at stage $0$.
		
		\emph{Inductive step.}
		Fix $\lambda>0$ and assume the statement holds for every $\beta<\lambda$.
		Let $H\in\tilde{\cF}^*_\lambda$.
		By definition of “normal filter generated by a family”,
		there exist generators $G_0,\dots,G_{n-1}$ from Definition~\ref{def:modified-limit-filter} and
		elements $k_0,\dots,k_{n-1}\in\cG_\lambda$ such that
		\[
		\bigcap_{i<n} k_i G_i k_i^{-1}\ \subseteq\ H.
		\]
		We show that each conjugate $k_i G_i k_i^{-1}$ contains a globally admissible
		$\Delta^\uparrow_\lambda$-group.
		
		\smallskip
		\noindent\emph{Case 1: $G_i=\Delta_\lambda^\uparrow(E,D)$ with $(E,D)$ globally admissible at stage $\lambda$.}
		By Lemma~\ref{lem:delta-up-general-conj},
		\[
		k_i\,\Delta_\lambda^\uparrow(E,D)\,k_i^{-1}=\Delta_\lambda^\uparrow(k_i\cdot E,\ k_i\cdot D).
		\]
		Apply Lemma~\ref{lem:admissible-hull} to the countable set $k_i\cdot D$ and the countable set $k_i\cdot E$
		to obtain a countable $E_i^\ast\supseteq k_i\cdot E$ such that $(E_i^\ast,\ k_i\cdot D)$ is globally admissible.
		Then
		\[
		\Delta_\lambda^\uparrow(E_i^\ast,\ k_i\cdot D)\ \le\ \Delta_\lambda^\uparrow(k_i\cdot E,\ k_i\cdot D)
		= k_i G_i k_i^{-1}.
		\]
		Set $(E_i,D_i):=(E_i^\ast,\ k_i\cdot D)$.
		
		\smallskip
		\noindent\emph{Case 2: $G_i$ is a pullback generator.}
		Then $G_i=(\pi^\beta_\lambda)^{-1}[K]$ for some $\beta<\lambda$ and some $K\in\tilde{\cF}^*_\beta$
		(as in Definition~\ref{def:modified-limit-filter}).
		Put $\sigma:=\pi^\beta_\lambda(k_i)\in\cG_\beta$.
		Since $\pi^\beta_\lambda$ is a group homomorphism
		(Remark~\ref{rem:projection-coherence}), one has
		\[
		k_i\,(\pi^\beta_\lambda)^{-1}[K]\,k_i^{-1}=(\pi^\beta_\lambda)^{-1}\!\big[\sigma K\sigma^{-1}\big].
		\]
		Because $\tilde{\cF}^*_\beta$ is normal, $\sigma K\sigma^{-1}\in\tilde{\cF}^*_\beta$.
		By the inductive hypothesis at stage $\beta$, choose a globally admissible pair $(E,D)$ at stage $\beta$
		with $\Delta_\beta^\uparrow(E,D)\le \sigma K\sigma^{-1}$.
		Then Lemma~\ref{lem:pullback-delta} yields
		\[
		\Delta_\lambda^\uparrow(E,D)\ \le\ (\pi^\beta_\lambda)^{-1}\!\big[\Delta_\beta^\uparrow(E,D)\big]
		\ \le\ (\pi^\beta_\lambda)^{-1}\!\big[\sigma K\sigma^{-1}\big]
		= k_i G_i k_i^{-1}.
		\]
		Since $\beta<\lambda$, the same pair $(E,D)$ is globally admissible at stage $\lambda$ as well.
		Set $(E_i,D_i):=(E,D)$.
		
		\smallskip
		\noindent\emph{Assembling the finite intersection.}
		Let $E^\ast:=\bigcup_{i<n}E_i$ and $D^\ast:=\bigcup_{i<n}D_i$.
		Since $n$ is finite, both unions are countable, and by iterating
		Corollary~\ref{cor:admissibility-union} the pair $(E^\ast,D^\ast)$ is globally admissible at stage $\lambda$.
		Moreover, by Corollary~\ref{cor:delta-up-intersection} (finite case),
		\[
		\Delta_\lambda^\uparrow(E^\ast,D^\ast)\ \le\ \bigcap_{i<n}\Delta_\lambda^\uparrow(E_i,D_i)
		\ \le\ \bigcap_{i<n} k_i G_i k_i^{-1}\ \subseteq\ H.
		\]
		This completes the induction and proves the lemma.
	\end{proof}
	
	\begin{corollary}[Admissibility preserved under countable unions]
		\label{cor:admissibility-countable-union}
		If $\{(E_n, D_n) : n < \omega\}$ is a countable family of globally admissible pairs, 
		then $(\bigcup_{n} E_n, \bigcup_{n} D_n)$ is globally admissible.
	\end{corollary}
	
	\begin{proof}
		Put $E:=\bigcup_{n<\omega}E_n$ and $D:=\bigcup_{n<\omega}D_n$.
		Since we are working in the ambient universe in which the iteration is constructed,
		countable unions of countable sets are countable;
		hence $E\in[\omega_1]^{\le\omega}$ and
		$D$ is countable. (No claim is made that such unions are countable \emph{inside} the eventual symmetric model.)
		
		To verify clause (iii) of Definition~\ref{def:global-admissible}, fix $d\in D$ and choose $k$ with $d\in D_k$.
		Then $\Fix(E)\le \Fix(E_k)\le \Stab_{\cG}(d)$, where the last inequality holds because $(E_k,D_k)$ is globally admissible.
		Therefore $(E,D)$ is globally admissible.
	\end{proof}
	
	\begin{lemma}[$\tilde{\cF}^*_\lambda$ is a proper normal filter]
		\label{lem:filter-properties}
		For each stage $\lambda$:
		\begin{enumerate}[label=(\roman*)]
			\item $\tilde{\cF}^*_\lambda$ is a filter (upward closed, closed under finite intersections).
			\item $\tilde{\cF}^*_\lambda$ is normal (closed under conjugation by $\cG_\lambda$).
			\item $\tilde{\cF}^*_\lambda$ is proper ($\{\mathrm{id}\} \notin \tilde{\cF}^*_\lambda$).
			\item $\tilde{\cF}^*_\lambda$ is $\omega_1$-complete (i.e.\ closed under countable intersections).
		\end{enumerate}
	\end{lemma}
	
	\begin{proof}
		(i) and (ii) follow from the definition as a normal filter generated by the given families.
		
		(iii) By Lemma~\ref{lem:core-lemma}, if $\{\mathrm{id}\} \in \tilde{\cF}^*_\lambda$,
		then $\{\mathrm{id}\} \supseteq \Delta_\lambda^\uparrow(E, D)$ for some countable globally admissible $E, D$.
		Since $E$ is countable and $\omega_1$ is uncountable, choose distinct $\beta,\gamma\in\omega_1\setminus E$ and
		let $\pi=(\beta\ \gamma)$. Then $\pi\in\Fix(E)$.
		By Lemma~\ref{lem:diag-lift-in-group}, the diagonal lift $\widehat{\pi}^{\,D}\in\cG_\lambda$ exists, and by
		Lemma~\ref{lem:diag-lift-in-delta} we have $\widehat{\pi}^{\,D}\in\Delta_\lambda^\uparrow(E,D)$.
		By Corollary~\ref{cor:diag-lift-swaps-cohen}, $\widehat{\pi}^{\,D}$ swaps $c_\beta\leftrightarrow c_\gamma$,
		so $\widehat{\pi}^{\,D}\neq\mathrm{id}$, contradicting $\Delta_\lambda^\uparrow(E,D)\le\{\mathrm{id}\}$.
		
		(iv) Let $\{H_n:n<\omega\}\subseteq\tilde{\cF}^*_\lambda$.
		By Lemma~\ref{lem:core-lemma}, for each $n$ choose
		a globally admissible $(E_n,D_n)$ with $\Delta_\lambda^\uparrow(E_n,D_n)\le H_n$.
		By Corollary~\ref{cor:delta-up-intersection},
		\[
		\bigcap_{n<\omega}H_n\ \supseteq\ \Delta_\lambda^\uparrow\Bigl(\bigcup_{n}E_n,\ \bigcup_{n}D_n\Bigr).
		\]
		By Corollary~\ref{cor:admissibility-countable-union}, the pair
		$\bigl(\bigcup_n E_n,\bigcup_n D_n\bigr)$ is globally admissible, and in the ambient universe where the
		iteration and filters are defined these unions are countable.
		Hence the right-hand side is one of the
		$\Delta^\uparrow$-generators in Definition~\ref{def:modified-limit-filter}, so by upward closure
		$\bigcap_{n<\omega}H_n\in\tilde{\cF}^*_\lambda$.
		
		At $\lambda=\Ord$, the same argument applies with $\cG_{\Ord}=\cG_{\Ord}^{bd}$: the unions
		$\bigcup_n E_n$ and $\bigcup_n D_n$ remain countable, so $\Delta_{\Ord}^\uparrow(\bigcup_n E_n,\bigcup_n D_n)$
		is a generator of $\tilde{\cF}^*_{\Ord}$ by Lemma~\ref{lem:GOrd-definable}. No appeal to the full
		inverse limit is needed.
	\end{proof}
	
	\begin{definition}[Hereditarily symmetric names]
		\label{def:hs-filter}
		A $\bbP_\lambda$-name $\tau$ is \emph{hereditarily symmetric} ($\tau \in \HS_\lambda$) iff:
		\begin{enumerate}[label=(\roman*)]
			\item $\Stab_{\cG_\lambda}(\tau) \in \tilde{\cF}^*_\lambda$, and
			\item for every $\langle \sigma,p\rangle \in \tau$, we have $\sigma \in \HS_\lambda$.
		\end{enumerate}
		The \emph{symmetric extension} at stage $\lambda$ is 
		$M_\lambda = \HS_\lambda^{G_\lambda} = \{\tau^{G_\lambda} : \tau \in \HS_\lambda\}$.
	\end{definition}
	
	\begin{lemma}[$\Delta$-support equivalence]
		\label{lem:delta-support-equiv}
		For any $\bbP_\lambda$-name $\tau$,
		\[
		\Stab_{\cG_\lambda}(\tau) \in \tilde{\cF}^*_\lambda
		\iff
		\exists (E,D)\ \bigl[(E,D)\text{ is globally admissible at stage }\lambda\ \wedge\
		\Delta_\lambda^\uparrow(E, D) \subseteq \Stab_{\cG_\lambda}(\tau)\bigr].
		\]
	\end{lemma}
	
	\begin{proof}
		$(\Leftarrow)$ If $\Delta_\lambda^\uparrow(E,D)\subseteq \Stab_{\cG_\lambda}(\tau)$ and
		$(E,D)$ is globally admissible, then $\Delta_\lambda^\uparrow(E,D)\in\tilde{\cF}^*_\lambda$
		(Definition~\ref{def:modified-limit-filter}), so by upward closure
		$\Stab_{\cG_\lambda}(\tau)\in\tilde{\cF}^*_\lambda$.
		
		$(\Rightarrow)$ Apply Lemma~\ref{lem:core-lemma} to $H=\Stab_{\cG_\lambda}(\tau)\in\tilde{\cF}^*_\lambda$.
	\end{proof}
	
	\begin{theorem}[$M_\lambda \models \ZF$]
		\label{thm:zf-symmetric}
		For each stage $\lambda$, the symmetric extension $M_\lambda = \HS_\lambda^{G_\lambda}$ 
		satisfies $\ZF$.
	\end{theorem}
	
	\begin{proof}
		By Lemma~\ref{lem:filter-properties}(i),(ii), the triple
		$(\bbP_\lambda,\cG_\lambda,\tilde{\cF}^*_\lambda)$ forms a symmetric system in the sense of
		\cite[Ch.~15]{JechSetTheory}.
		Therefore the associated symmetric extension
		$M_\lambda=\HS_\lambda^{G_\lambda}$ satisfies $\ZF$ by the symmetric extension theorem
		\cite[Lemma~15.51]{JechSetTheory}.
	\end{proof}
	
	\begin{lemma}[Upward persistence of hereditarily symmetric names]
		\label{lem:HS-persistence-up}
		Let \(\beta\le\lambda\) be stages, where \(\lambda\) may also be \(\Ord\).
		If \(\tau\in\HS_\beta\), then, viewed as a \(\bbP_\lambda\)-name via canonical inclusion,
		\(\tau\in\HS_\lambda\).
	\end{lemma}
	
	\begin{proof}
		We argue by induction on the rank of \(\tau\).
		Since \(\tau\in\HS_\beta\), every subname \(\sigma\) of \(\tau\) belongs to \(\HS_\beta\).
		By the induction hypothesis, each such \(\sigma\) also belongs to \(\HS_\lambda\).
		
		It remains to prove that \(\Stab_{\cG_\lambda}(\tau)\in\tilde{\cF}^*_\lambda\).
		Because \(\tau\in\HS_\beta\), we have
		\[
		\Stab_{\cG_\beta}(\tau)\in\tilde{\cF}^*_\beta.
		\]
		By Definition~\ref{def:modified-limit-filter}, the pullback
		\[
		(\pi^\beta_\lambda)^{-1}\!\big[\Stab_{\cG_\beta}(\tau)\big]
		\]
		belongs to \(\tilde{\cF}^*_\lambda\).
		
		We claim that
		\[
		(\pi^\beta_\lambda)^{-1}\!\big[\Stab_{\cG_\beta}(\tau)\big]
		\subseteq
		\Stab_{\cG_\lambda}(\tau).
		\]
		Indeed, if \(k\in\cG_\lambda\) satisfies \(\pi^\beta_\lambda(k)\in\Stab_{\cG_\beta}(\tau)\),
		then by Lemma~\ref{lem:projection-action-bounded-name},
		\[
		k\cdot\tau
		=
		\pi^\beta_\lambda(k)\cdot\tau
		=
		\tau.
		\]
		So \(k\in\Stab_{\cG_\lambda}(\tau)\), proving the claim.
		
		By upward closure of \(\tilde{\cF}^*_\lambda\), it follows that
		\(\Stab_{\cG_\lambda}(\tau)\in\tilde{\cF}^*_\lambda\). Therefore \(\tau\in\HS_\lambda\).
	\end{proof}
	
	\begin{corollary}[Uniform definability preserves symmetry]
		\label{cor:uniformity-HS}
		Let $\beta\le\lambda$ be set stages.
		Let $\dot x \in \HS_{\lambda}$ with
		$\Delta^{\uparrow}_{\lambda}(E,D) \le \Stab_{\cG_{\lambda}}(\dot x)$.
		If $\dot y$ is definable in $V$ from parameters in
		$\HS_{\beta} \cup \{\dot x\}$ using only absolute operations,
		then $\dot y \in \HS_{\lambda}$ and
		$\Delta^{\uparrow}_{\lambda}(E,D) \le \Stab_{\cG_{\lambda}}(\dot y)$.
	\end{corollary}
	
	\begin{remark}
		The corollary is stated and proved for set stages $\lambda$ only. Applicability at
		$\lambda=\Ord$ would require $\Delta^{\uparrow}_{\Ord}(E,D)$ to exist as a class in
		$\GBC+\ETR_{\Ord}$; this is established by Lemma~\ref{lem:GOrd-definable}, but the
		corollary is not invoked at $\lambda=\Ord$ in this paper.
	\end{remark}
	
	\begin{proof}
		Let $\pi \in \Delta^{\uparrow}_{\lambda}(E,D)$. Since $\tau\in\HS_{\beta}$ implies $\pi\cdot\tau=\tau$ for every $\tau\in\HS_{\beta}$ by Lemma~\ref{lem:HS-persistence-up}, $\pi$ fixes all parameters from $\HS_{\beta}$. By hypothesis $\pi$ fixes $\dot x$. Because $\dot y$ is definable from $\HS_{\beta}\cup\{\dot x\}$ using only absolute operations, $\pi$ fixes $\dot y$. Hence $\Delta^{\uparrow}_{\lambda}(E,D) \le \Stab_{\cG_{\lambda}}(\dot y)$. Since $\lambda$ is a set stage, $\Delta^{\uparrow}_{\lambda}(E,D) \in \tilde{\cF}^*_\lambda$ by Definition~\ref{def:modified-limit-filter}. By upward closure of $\tilde{\cF}^*_\lambda$, $\Stab_{\cG_{\lambda}}(\dot y)\in\tilde{\cF}^*_\lambda$. Therefore $\dot y \in \HS_{\lambda}$.
	\end{proof}
	
	\begin{lemma}[Orbit forcings are $S$-preserving]
		\label{lem:orbit-S-preserving}
		Let $\alpha$ be a stage and let $\bbQ = \bbQ_{f}$ be the orbit forcing adding a section for $f: Y \twoheadrightarrow X$ with $X,Y \subseteq \Pow(S)^{M_{\alpha}}$. Then $\bbQ$ adds no new subsets of $S = A^{\omega}$.
	\end{lemma}
	
	\begin{proof}
		Work in $M_{\alpha}$. Let $\kappa = \aleph^{*}(S)^{M_{\alpha}}$. Conditions of $\bbQ$ are partial functions $p: X_{0} \to Y$ with $|X_{0}| < \kappa$, so $\bbQ$ is $\kappa$-closed. Since $S$ is countable in $M_{\alpha}$, $\kappa > \aleph_{0}$ implies $\bbQ$ adds no new $\omega$-sequences of ground-model elements, hence no new subsets of $S$.
		
		For symmetry: if $\dot T \in \HS_{\alpha+1}$ and $1\Vdash \dot T\subseteq\check S$, then by the first paragraph $\dot T^{G}\in M_{\alpha}$ for every generic $G$. We claim $1 \Vdash \dot T = \check t$ for some $t\in M_{\alpha}$. Fix any $p\in\bbQ$. Since $\bbQ$ is $\kappa$-closed and $\kappa > \omega$, for each $n\in\omega$ the set $D_n = \{q\le p : q \Vdash \check n \in \dot T \text{ or } q \Vdash \check n \notin \dot T\}$ is dense below $p$. Construct by recursion a decreasing sequence $p = p_0 \ge p_1 \ge \cdots$ with $p_n \in D_n$ for each $n$. Since $\omega < \kappa$, $q = \bigcap_{n<\omega} p_n$ is a condition in $\bbQ$ with $q\le p$, and $q \Vdash \dot T = \check t$ where $t = \{n\in\omega : q \Vdash \check n \in \dot T\} \in M_{\alpha}$. Thus the set of conditions forcing $\dot T = \check t$ for some $t\in M_{\alpha}$ is dense. Since $\bbQ$ adds no new subsets of $S$, any two such conditions force the same $t$, so $1 \Vdash \dot T = \check t$. Since $\check t$ is a check-name for $t\in M_{\alpha}$ and $\Stab_{\cG_{\alpha}}(\check t) = \cG_{\alpha} \in \tilde{\cF}^*_{\alpha}$, we have $\check t \in \HS_{\alpha}$. Therefore $\dot T\in\HS_{\alpha}$.
		\end{proof}
	
	\begin{lemma}[Name covering under SVC]
		\label{lem:name-covering}
		Assume $M_{\alpha} \models \SVC(S)$. Let $\dot X \in \HS_{\alpha+1}$.
		Then there exist $Y \in M_{\alpha}$ and a name
		$\dot e \in \HS_{\alpha+1}$ such that
		$1 \Vdash_{\bbP_{\alpha+1}} \dot e : \check Y \twoheadrightarrow \dot X$.
	\end{lemma}
	
	\begin{proof}
		Work in $M_{\alpha}$. Let $\bbQ = \bbQ_{f}$ and $\kappa = \aleph^{*}(S)^{M_{\alpha}}$. By $\SVC(S)^{M_{\alpha}}$, fix $\eta$ with $s_{0}: S \times \eta \twoheadrightarrow Y_{0}$ where $Y_{0}$ is the domain of $f$. Conditions of $\bbQ$ are partial functions $X_{0} \to Y_{0}$ with $|X_{0}| < \kappa$, so $|\bbQ|^{M_{\alpha}} \le |Y_{0}|^{<\kappa} \le |S \times \eta|^{<\kappa}$, a set in $M_{\alpha}$.
		
		Choose $(E,D)$ with $\Delta^{\uparrow}_{\alpha+1}(E,D) \le \Stab(\dot X)$. Build by $\DC$ a maximal antichain $A \subseteq \bbQ$ with $a \Vdash \dot X = \check{t_{a}}$. Let $Y = \{t_{a} : a \in A\}$ and let $\dot e$ name $a \mapsto t_{a}$.
		
		For $\pi \in \Delta^{\uparrow}_{\alpha+1}(E,D)$, since $\Delta^{\uparrow}_{\alpha+1}(E,D) \le \Stab(\dot X)$ we have $\pi\in\Stab(\dot X)$. Applying $\pi$ to $a \Vdash \dot X = \check{t_{a}}$ yields $\pi(a) \Vdash \pi(\dot X) = \pi(\check{t_{a}})$. Because $\pi(\dot X)=\dot X$ and $\pi(\check{t_{a}})=\check{t_{a}}$, we obtain $\pi(a) \Vdash \dot X = \check{t_{a}}$. By maximality of $A$, $\pi(a)$ is compatible with some $b \in A$. Pick a common extension $r\le\pi(a),b$. Then $r \Vdash \dot X = \check{t_{a}}$ and $r \Vdash \dot X = \check{t_{b}}$, hence $r \Vdash \check{t_{a}} = \check{t_{b}}$, so $t_{a} = t_{b}$. Since $A$ is an antichain, such $b$ is unique. Hence $t_{\pi(a)} = t_{a}$ for all $a\in A$. Therefore $\pi(\dot e)$ is the name for the map $\pi(a)\mapsto t_{\pi(a)} = t_a$, which equals $\dot e$. Thus $\pi\in\Stab_{\cG_{\alpha+1}}(\dot e)$, so $\Delta^{\uparrow}_{\alpha+1}(E,D) \le \Stab_{\cG_{\alpha+1}}(\dot e)$. Since $\alpha+1$ is a set stage, $\Delta^{\uparrow}_{\alpha+1}(E,D) \in \tilde{\cF}^*_{\alpha+1}$ by Definition~\ref{def:modified-limit-filter}. By upward closure, $\Stab_{\cG_{\alpha+1}}(\dot e)\in\tilde{\cF}^*_{\alpha+1}$. By Corollary~\ref{cor:uniformity-HS}, $\dot e \in \HS_{\alpha+1}$.
	\end{proof}
	
	\begin{remark}
		The hypothesis $M_{\alpha} \models \SVC(S)$ is used only to bound $|\bbQ|$. In Lemma~\ref{lem:svc-persists} this lemma is invoked \emph{after} the inductive hypothesis is available, so no circularity occurs.
	\end{remark}
	
	\subsection{Bookkeeping: names versus interpretations}
	\label{subsec:bookkeeping}
	
	\begin{remark}[Metatheory for bookkeeping]
		\label{rem:metatheory}
		We continue with the metatheoretic conventions of Remark~\ref{rem:meta-pp}.
		For the \(\Ord\)-length bookkeeping in Section~\ref{subsec:bookkeeping} (in particular, to define the
		class function \(\mathcal B:\Ord\to V\) in Definition~\ref{def:bookkeeping-class}),
		we work in a background two-sorted structure \((V,\mathcal C)\models\GBC+\ETR\) and,
		only for this bookkeeping/canonical-enumeration purpose, fix a class well-order of \(V\)
		(equivalently, assume Global Choice in the background).
		All forcing semantics for \(\bbP_{\Ord}\) are understood externally exactly as in Remark~\ref{rem:meta-pp}
		(in particular, we do not assume a global forcing theorem for \(\bbP_{\Ord}\)).
	\end{remark}
	
	We now fix a bookkeeping device which ensures that every relevant surjection instance appearing
	in the final model is handled at some successor stage.
	
	\begin{remark}[Syntactic bookkeeping]
		\label{rem:names-vs-interpretation}
		In a class-length iteration, the interpretation of a fixed name can change as the iteration proceeds.
		Therefore, we bookkeep \emph{names} (codes) rather than interpreted objects: at stage \(\alpha\) we take the
		\(\alpha\)-th code \(\dot f\) from our fixed enumeration, and in the current intermediate model \(M_\alpha\)
		we check whether \(\dot f\) (as currently interpreted) is a surjection instance of the required type.
		If it is not, we force trivially at that stage (or include only the other package factors scheduled for \(\alpha\)).
		We also use the canonical inclusion of \(\bbP_\beta\)-names into \(\bbP_\alpha\)-names for \(\beta<\alpha\),
		so a code intended for an earlier stage can be retested later.
		This is the standard bookkeeping convention for class-length forcing/iterations.
	\end{remark}
	
	\begin{definition}[Bookkeeping of orbit instances]
		\label{def:bookkeeping-class}
		Fix in the metatheory a definable class function $\mathcal B:\Ord\to V$ that enumerates 
		\emph{pairs} $(\beta, c)$ where $c$ is a set-theoretic code for a potential $\bbP$-name 
		(in some standard Gödel coding).
		The enumeration is arranged so that for each code $c$, 
		the set
		\[
		\{\alpha\in\Ord : \exists\,\beta\ (\mathcal B(\alpha)=(\beta,c))\}
		\]
		is unbounded in $\Ord$ (e.g.\ enumerate $\Ord\times V$ using a definable pairing of ordinals
		together with a definable class well-order of $V$ in the background metatheory).
		At stage $\alpha$, let $\mathcal B(\alpha)=(\beta,c)$.
		We attempt to decode $c$ as a $\bbP_\beta$-name $\dot f$ (in some fixed coding of names),
		and when $\beta\le\alpha$ we view $\dot f$ as a $\bbP_\alpha$-name via the canonical inclusion
		of $\bbP_\beta$-names into $\bbP_\alpha$-names.
		We then \emph{test} whether (in the current stage model $M_\alpha$):
		\begin{enumerate}[label=(\roman*)]
			\item $\dot f$ is defined (the decoding succeeds) and $\beta\le\alpha$,
			\item $\dot f\in\HS_\alpha$ for the stage-$\alpha$ symmetric system $(\bbP_\alpha,\cG_\alpha,\tilde{\cF}^*_\alpha)$, and
			\item $f:=\dot f^{G_\alpha}$ is a surjection instance of one of the following types:
			\begin{enumerate}[label=(\arabic*)]
				\item $f:Y\twoheadrightarrow X$ with $X,Y\subseteq \Pow(S)^{M_\alpha}$ (for $\PP$-packages), or
				\item $f:S\times\eta\twoheadrightarrow \lambda$ with $\eta,\lambda\in\Ord$ and
				\(\lambda<\aleph^*(S)^{M_\alpha}\) (for $\AC_{\WO}$-packages),
			\end{enumerate}
		\end{enumerate}
		and moreover $f$ does not already split in $M_\alpha$ (i.e., no section/right-inverse for $f$ exists in $M_\alpha$).
		If the test succeeds, we define the stage-$\alpha$ iterand $\dot{\bbQ}_\alpha$ so that in $M_\alpha$ it is the
		corresponding orbit package $\bbQ_{[f]}$ or $\bbR_{[f]}$ (Definition~\ref{def:package-systems});
		otherwise we set $\dot{\bbQ}_\alpha$ to be the trivial forcing. Thus $\bbP_{\alpha+1}=\bbP_\alpha*\dot{\bbQ}_\alpha$.
		
		This arrangement ensures that every code is tested at unboundedly many stages. Hence:
		\begin{itemize}
			\item a \(\PP\)-instance that becomes a valid HS surjection instance \(f:Y\twoheadrightarrow X\) with
			\(X,Y\subseteq \Pow(S)\) in some stage model will eventually be scheduled; and
			\item an \(\AC_{\WO}\)-instance that becomes a valid HS surjection instance and satisfies
			\(\lambda<\aleph^*(S)\) in some stage model will eventually be scheduled at some later stage.
		\end{itemize}
	\end{definition}
	
	\begin{remark}[Stage recognition of \(\Pow(S)\)-instances]
		\label{rem:stage-recognizes-pows}
		Let \(\alpha<\Ord\). If \(X\in M_\alpha\) and \(X\subseteq \Pow(S)^{\cM}\), then in fact
		\[
		X\subseteq \Pow(S)^{M_\alpha}.
		\]
		Indeed, if \(x\in X\), then \(x\in M_\alpha\) by transitivity of \(M_\alpha\), and \(x\subseteq S\) since \(X\subseteq \Pow(S)^{\cM}\) and \(S\in M_\alpha\). Hence \(x\in \Pow(S)^{M_\alpha}\).
		The same holds for any \(Y\in M_\alpha\) with \(Y\subseteq \Pow(S)^{\cM}\).
	\end{remark}
	
	\begin{lemma}[Bookkeeping guarantee]\label{lem:bookkeeping-guarantee}
		Let \(f\in\cM\) be a surjection instance of one of our two templates. Assume either:
		\begin{enumerate}[label=(\alph*)]
			\item \(f:Y\twoheadrightarrow X\) with \(X,Y\subseteq \Pow(S)^{\cM}\); or
			\item \(f:S\times\eta\twoheadrightarrow\lambda\), and there is some stage model \(M_\gamma\)
			containing a surjection \(e:S\twoheadrightarrow\lambda\).
		\end{enumerate}
		Then \(f\) has the required splitting object in \(\cM\)
		(a section in case~\textup{(a)}, a right inverse in case~\textup{(b)}).
	\end{lemma}
	
	\begin{proof}
		Choose a set-sized \(\dot f\in\HS_{\Ord}\) with \(f=\dot f^{G_\infty}\).
		By Lemma~\ref{lem:Ord-HS-stage-descent}, there is \(\beta_0<\Ord\) such that
		\(\dot f\), viewed as a \(\bbP_{\beta_0}\)-name, belongs to \(\HS_{\beta_0}\) and
		\[
		f=\dot f^{G_{\beta_0}}\in M_{\beta_0}.
		\]
		
		In case~\textup{(a)}, let \(\beta:=\beta_0\).
		
		In case~\textup{(b)}, choose \(\beta_1\ge\gamma\) such that \(e\in M_{\beta_1}\), and let
		\(\beta\ge\beta_0,\beta_1\). Viewing \(\dot f\) by canonical inclusion as a \(\bbP_\beta\)-name,
		we still have \(\dot f\in\HS_\beta\) by Lemma~\ref{lem:HS-persistence-up} and
		\[
		f=\dot f^{G_\beta}\in M_\beta.
		\]
		Since \(e\in M_\beta\) is a surjection \(S\twoheadrightarrow\lambda\), we have
		\(\lambda<\aleph^*(S)^{M_\beta}\). Because \(e\) remains a set in every later stage model,
		the same bound \(\lambda<\aleph^*(S)\) holds in every \(M_\alpha\) with \(\alpha\ge\beta\).
		
		Let \(c_{\dot f}\) be a code for \(\dot f\) as a \(\bbP_\beta\)-name. By the unboundedness clause in
		Definition~\ref{def:bookkeeping-class}, there exists \(\alpha\ge\beta\) with
		\(\mathcal B(\alpha)=(\beta,c_{\dot f})\).
		
		At stage \(\alpha\), clauses \textup{(i)} and \textup{(ii)} of
		Definition~\ref{def:bookkeeping-class} hold for \(\dot f\).
		In case~\textup{(a)}, since \(X,Y,f\in M_\alpha\) and \(X,Y\subseteq \Pow(S)^{\cM}\),
		Remark~\ref{rem:stage-recognizes-pows} gives
		\[
		X,Y\subseteq \Pow(S)^{M_\alpha},
		\]
		so clause~\textup{(iii)(1)} holds.
		In case~\textup{(b)}, clause \textup{(iii)(2)} holds because
		\(f:S\times\eta\twoheadrightarrow\lambda\) and the bound
		\(\lambda<\aleph^*(S)^{M_\alpha}\) holds by the previous paragraph.
		
		If \(f\) already has the required splitting object in \(M_\alpha\), we are done.
		
		Otherwise the bookkeeping test succeeds, so the stage-\(\alpha\) iterand is the
		corresponding orbit package \(\bbQ_{[f]}\) or \(\bbR_{[f]}\).
		
		In case~\textup{(a)}, Proposition~\ref{prop:Qf-adds-section} gives a generic section for \(f\),
		and Lemma~\ref{lem:sf-HS-diag} shows that its canonical name is hereditarily symmetric for the
		package system. Hence its valuation belongs to \(M_{\alpha+1}\subseteq\cM\).
		
		In case~\textup{(b)}, Proposition~\ref{prop:Rf-adds-section} gives a generic right inverse for \(f\),
		and Lemma~\ref{lem:rg-HS-diag} shows that its canonical name is hereditarily symmetric for the
		package system. Hence its valuation belongs to \(M_{\alpha+1}\subseteq\cM\).
		
		Thus \(f\) has the required splitting object in \(\cM\).
	\end{proof}
	
	\begin{remark}[The $\aleph^*(S)$ horizon]
		\label{rem:aleph-star-bookkeeping}
		For the $\AC_{\WO}$-component we schedule splittings for surjections
		$f:S\times\eta\twoheadrightarrow\lambda$ with $\lambda<\aleph^*(S)$ (cf.\ Remark~\ref{rem:ACwo-bound}).
		This is the natural bookkeeping horizon: $\aleph^*(S)$ is the first ordinal $\theta$
		(in the ambient stage model under discussion) such that there is no surjection
		$S\twoheadrightarrow\theta$, so $\lambda<\aleph^*(S)$ is exactly the range of well-ordered
		indices reachable from $S$-parameterized data.
		For our purposes we only use the following two points.
		\begin{enumerate}[label=(\roman*)]
			\item The seed parameter $S=A^\omega$ is preserved at successor stages: our package forcings
			add no new reals (Corollaries~\ref{cor:Qf-no-reals-DC} and~\ref{cor:Rf-no-reals-DC}),
			hence add no new $\omega$-sequences of Cohen reals and therefore no new elements of $S$.
			\item The value of $\aleph^*(S)$ need not be constant along the iteration (new functions from
			$S$ may be added even if $S$ itself is unchanged), but this causes no bookkeeping gap:
			we test the bound $\lambda<\aleph^*(S)$ \emph{in the current stage model} when a code is examined.
			If at some stage a given name $\dot f$ interprets to an $\AC_{\WO}$-instance
			\(f:S\times\eta\twoheadrightarrow\lambda\) with \(\lambda<\aleph^*(S)\)
			(as computed in that stage), then from that point onward the instance
			is eligible under clause~\textup{(iii)(2)} of Definition~\ref{def:bookkeeping-class}.
			Since every code is retested at unboundedly many stages
			(Definition~\ref{def:bookkeeping-class}), no a priori stability theorem for $\aleph^*(S)$ is needed.
		\end{enumerate}
	\end{remark}
	
	\subsection{The \texorpdfstring{\(\Ord\)}{Ord}-length symmetric iteration}
	\label{subsec:iteration-api-level}
	
	We now record the specific \(\Ord\)-length recursion proposed in the
	former argument.
	Working in the background metatheory from Remark~\ref{rem:meta-pp}, we define by
	transfinite recursion on \(\alpha\in\Ord\) an iteration
	\(\langle \bbP_\alpha,\cG_\alpha,\tilde{\cF}^*_\alpha : \alpha\in\Ord\rangle\),
	starting from \(\bbP_0=\Add(\omega,\omega_1)\) and producing intermediate symmetric
	models \(M_\alpha\).
	At successor stages we use the orbit-package iterands described below;
	at limit stages
	we use the countable-support forcing limit and the corresponding group/projection
	infrastructure from the bounded-stage framework, together with the modified limit filter
	from Definition~\ref{def:modified-limit-filter}.
	No general class-length theorem is invoked at this point; only the
	bounded-stage infrastructure summarized in Remark~\ref{rem:iteration-api} is used.
	The next subsection records the properties formerly claimed for the proposed
	final symmetric model.
	
	At a successor stage \(\alpha+1\), we consult the bookkeeping output \(\mathcal{B}(\alpha)=(\beta,c)\)
	from Definition~\ref{def:bookkeeping-class}.
	In the current model \(M_\alpha\), the code \(c\) either
	decodes (via the canonical inclusion \(\bbP_\beta\)-names \(\subseteq\) \(\bbP_\alpha\)-names) to a genuine package
	instance, yielding a corresponding \emph{orbit package} forcing \(\bbQ_{[f]}\) or \(\bbR_{[f]}\) as in
	Definition~\ref{def:package-systems}, or it does not, in which case we take trivial forcing for that stage.
	Let \(\dot{\bbS}_\alpha\) be the resulting \(\bbP_\alpha\)-name for the stage iterand (orbit package or trivial).
	Then \(\bbP_{\alpha+1}=\bbP_\alpha * \dot{\bbS}_\alpha\).
	
	\begin{lemma}[Stage iterands are hereditarily symmetric]\label{lem:stage-iterand-HS}
		For every stage \(\alpha\), the \(\bbP_\alpha\)-name \(\dot{\bbS}_\alpha\) is an element of \(\HS_\alpha\).
	\end{lemma}
	\begin{proof}
		If the bookkeeping test fails at stage \(\alpha\), then \(\dot{\bbS}_\alpha\) is the canonical name for the trivial forcing, which is fixed by every automorphism and hence lies in \(\HS_\alpha\).
		If the bookkeeping test succeeds, then \(\dot{\bbS}_\alpha\) is (by construction) an orbit package forcing built from a scheduled surjection instance via the stage-\(0\) action.
		For any \(g\in\cG_\alpha\), its stage-\(0\) projection \(\pi^0_\alpha(g)\in\cG\) sends a surjection instance \(f\) to \(\pi^0_\alpha(g)\cdot f\),
		but preserves the orbit \([f]\).
		Consequently,
		\[
		g\cdot \dot{\bbS}_\alpha=\dot{\bbS}_\alpha.
		\]
		Thus \(\cG_\alpha\le\Stab(\dot{\bbS}_\alpha)\), so \(\Stab(\dot{\bbS}_\alpha)\in\tilde{\mathcal F}^\ast_\alpha\), and therefore \(\dot{\bbS}_\alpha\in\HS_\alpha\).
	\end{proof}
	
	At successor stages, the group \(\cG_{\alpha+1}\) and filter \(\tilde{\cF}^*_{\alpha+1}\) are those determined
	by the symmetric-iteration framework together with the diagonal infrastructure from
	Subsection~\ref{subsec:iteration-diag} and the modified stage filters from Definition~\ref{def:modified-limit-filter}.
	Concretely, \(\tilde{\cF}^*_{\alpha+1}\) is the normal filter generated by pullbacks
	\((\pi^\alpha_{\alpha+1})^{-1}[H]\) for \(H\in\tilde{\cF}^*_\alpha\), together with admissible subgroups
	\(\Delta_{\alpha+1}^\uparrow(E,D)\).
	In particular, for each active package coordinate \(d\) at stage \(\alpha+1\) whose countable support is
	witnessed by \(E_d\), the generator \(\Delta_{\alpha+1}^\uparrow(E_d,\{d\})\) belongs to
	\(\tilde{\cF}^*_{\alpha+1}\), so the canonical package generics are hereditarily symmetric by
	Lemmas~\ref{lem:sf-HS-diag} and~\ref{lem:rg-HS-diag}.
	
	At limit stages \(\lambda\), we take the countable-support limit forcing \(\bbP_\lambda\) and the corresponding
	limit group as provided by the symmetric-iteration framework (Remark~\ref{rem:iteration-api}(a),(b)).
	The limit-stage filter is \(\tilde{\cF}^*_\lambda\) as defined in Definition~\ref{def:modified-limit-filter},
	and its normality and \(\omega_1\)-completeness are recorded in Lemma~\ref{lem:filter-properties}.
	
	For the maximal limit $\lambda=\Ord$ we replace the full inverse limit by its bounded-support version. For countable $(E,D)$ put $\alpha_0(E,D):=\max(\sup(\mathrm{ran}D),\sup E)$, where $\mathrm{ran}D:=\{\xi:\exists z\,\langle\xi,z\rangle\in D\}$. Define
	\[
	\begin{split}
		\cG_{\Ord}^{bd}:=\Bigl\{k=\langle k_\alpha\rangle_{\alpha<\Ord}\in\prod_{\alpha<\Ord}\cG_\alpha :
		\exists\text{ countable globally admissible }(E,D)\;\\
		\forall\alpha>\alpha_0(E,D)\;\bigl[k_\alpha\in\Delta^{\uparrow}_\alpha(E,D)\land\forall\beta<\alpha\;\pi^\beta_\alpha(k_\alpha)=k_\beta\bigr]\Bigr\}.
	\end{split}
	\]
	For $\alpha\le\alpha_0(E,D)$ the coherence $\pi^\beta_\alpha(k_\alpha)=k_\beta$ holds by the set-stage projection data supplied by Remark~\ref{rem:iteration-api}(a).
	
	For countable $(E,D)$ put
	\[
	\Delta^{\uparrow}_{\Ord}(E,D):=\{k\in\cG_{\Ord}^{bd} : \forall\alpha>\alpha_0(E,D)\;k_\alpha\in\Delta^{\uparrow}_\alpha(E,D)\}.
	\]
	
	\begin{lemma}[The bounded limit as a $\GBC$ class]\label{lem:GOrd-definable}
		In $\GBC+\ETR_{\Ord}$, $\cG_{\Ord}^{bd}$ and each $\Delta^{\uparrow}_{\Ord}(E,D)$ exist as
		classes. More precisely, $\ETR_{\Ord}$ guarantees the existence of a class coding the
		inverse system $\langle\cG_\gamma,\pi^\beta_\gamma\rangle_{\gamma<\Ord}$, and given that
		class as a parameter, membership in $\cG_{\Ord}^{bd}$ and in each
		$\Delta^{\uparrow}_{\Ord}(E,D)$ is first-order definable.
	\end{lemma}
	
	\begin{proof}
		"$k\in\cG_{\Ord}^{bd}$" is the condition
		\[
		\begin{split}
			\exists E\in[\omega_1]^{\le\omega}\exists D\in[\mathrm{Coords}_{<\Ord}]^{\le\omega}
			\Bigl(&\text{``$(E,D)$ globally admissible''}\land\\
			&\forall\gamma>\alpha_0(E,D)\;
			\bigl[k_\gamma\in\Delta^{\uparrow}_\gamma(E,D)\land\forall\beta<\gamma\;\pi^\beta_\gamma(k_\gamma)=k_\beta\bigr]\Bigr).
		\end{split}
		\]
		The clause $k_\gamma\in\Delta^{\uparrow}_\gamma(E,D)$ is first-order by induction on $\gamma$
		using Remark~\ref{rem:iteration-api}(a)(b), and the outer quantifiers range over sets.
		$\ETR_{\Ord}$ provides the class recursion whose solution is the inverse system
		$\langle\cG_\gamma,\pi^\beta_\gamma\rangle_{\gamma<\Ord}$; given that system as a class
		parameter, the condition above is first-order, so class comprehension in $\GBC$ yields
		$\cG_{\Ord}^{bd}$ as a class. Note that $\ETR_{\Ord}$ guarantees \emph{existence} of this
		class but not first-order definability in the sense of Tarski; the class exists as an object
		in the two-sorted structure $(V,\mathcal{C})$ but need not be definable without the ETR
		witness as a parameter. The same analysis applies to $\Delta^{\uparrow}_{\Ord}(E,D)$.
	\end{proof}
	
	\begin{remark}[What $\cG_{\Ord}^{bd}$ does and does not require]
		\label{rem:bd-not-trivial-on-omega1}
		A thread $k=\langle k_\gamma\rangle_{\gamma<\Ord}$ belongs to $\cG_{\Ord}^{bd}$ iff there exist countable
		$E\subseteq\omega_1$ and $D\subseteq\mathrm{Coords}_{<\Ord}$ such that
		$k_\gamma\in\Delta^\uparrow_\gamma(E,D)$ for all sufficiently large $\gamma$.
		In particular:
		\begin{itemize}
			\item It is \emph{not} required that $k_\gamma=\mathrm{id}$ for large $\gamma$, nor that the induced
			permutation $\pi^0_\gamma(k_\gamma)\in\Sym(\omega_1)$ be eventually trivial.
			\item For the diagonal lift $\widehat\pi^{\,D}$ of a fixed $\pi\in\Sym(\omega_1)$ with
			$\pi\in\Fix(E)$, we have $\widehat\pi^{\,D}\in\Delta^\uparrow_\gamma(E,D)$ for \emph{every} $\gamma$;
			hence $\widehat\pi^{\,D}\in\cG_{\Ord}^{bd}$ with $\alpha_0=0$, even when $\pi$ moves points of $\omega_1$
			(e.g., a transposition $(\alpha\ \beta)$).
		\end{itemize}
		The boundedness is on the \emph{iteration} coordinates, not on $\omega_1$.
	\end{remark}
	
	We henceforth write $\cG_{\Ord}$ for $\cG_{\Ord}^{bd}$ and use the groups $\Delta^{\uparrow}_{\Ord}(E,D)$ as the basic subgroups at stage $\Ord$.
	
	For set-localization in the \(\Ord\)-length model we use the forcing-localization lemma
	(Lemma~\ref{lem:Ord-stage-localization}) together with the descent lemma for set-sized final
	hereditarily symmetric names (Lemma~\ref{lem:Ord-HS-stage-descent}) below.
	We write \(\cM\) for the final symmetric model of the \(\Ord\)-length iteration.
	
	\begin{remark}[No one-step collapse in general]
		\label{rem:no-collapse}
		Unlike ordinary forcing iterations in $\ZFC$, there is no general theorem
		showing that an iteration of symmetric extensions can always be presented as a
		\emph{single} symmetric extension of the ground model.
		Indeed, the usual ``two-step collapse'' phenomena from forcing do not directly
		apply in the symmetric context, and the literature discusses serious obstacles
		to such a reduction.
		It remains open in general to what extent arbitrary
		symmetric iterations can be compressed to a single symmetric extension;
		see the
		discussion in \cite{KaragilaIteratingSymmExt}.
	\end{remark}
	
	\subsection{What the iteration yields}
	\label{subsec:iteration-yields}
	
	Let \(\cM\) denote the final symmetric model obtained from the iteration.
	
	\begin{lemma}[Forcing-localization of set-sized class-length names]
		\label{lem:Ord-stage-localization}
		Let \(\dot x\) be a \emph{set-sized} \(\bbP_{\Ord}\)-name, and define its stage-support
		\[
		I(\dot x)\ :=\ \bigcup\{\supp(p)\cup I(\dot y) : \langle \dot y,p\rangle\in\dot x\},
		\]
		where \(\supp(p)\subseteq\Ord\) is the (countable) support of the condition \(p\).
		Then \(I(\dot x)\) is a set of ordinals, hence
		\[
		\alpha:=\sup(I(\dot x))+1<\Ord.
		\]
		Moreover:
		\begin{enumerate}
			\item every subname of \(\dot x\) is a \(\bbP_\alpha\)-name, and in particular \(\dot x\) is a \(\bbP_\alpha\)-name;
			\item for any \(V\)-generic \(G_\infty\subseteq\bbP_{\Ord}\) and \(G_\alpha:=G_\infty\cap\bbP_\alpha\),
			one has \(\dot y^{G_\infty}=\dot y^{G_\alpha}\) for every subname \(\dot y\) of \(\dot x\), and in particular
			\(\dot x^{G_\infty}=\dot x^{G_\alpha}\).
		\end{enumerate}
	\end{lemma}
	
	\begin{proof}
		Define \(I(\dot x)\) by recursion on the rank of the name \(\dot x\).
		Since \(\dot x\) is set-sized, Replacement and Union imply that \(I(\dot x)\) is a set of ordinals,
		hence bounded in \(\Ord\).
		Let \(\alpha=\sup(I(\dot x))+1\).
		
		If \(\dot y\) is a subname of \(\dot x\), then by construction \(I(\dot y)\subseteq I(\dot x)\subseteq\alpha\).
		Hence every condition appearing in \(\dot y\) has support contained in \(\alpha\), so \(\dot y\) is a
		\(\bbP_\alpha\)-name. This proves (1).
		
		For (2), fix a subname \(\dot y\) of \(\dot x\).
		Since all conditions appearing in \(\dot y\) lie in
		\(\bbP_\alpha\), the valuation recursion for \(\dot y^{G_\infty}\) consults only \(G_\alpha\), so
		\(\dot y^{G_\infty}=\dot y^{G_\alpha}\).
		The special case \(\dot y=\dot x\) gives the final statement.
	\end{proof}
	
	\begin{lemma}[Stage descent of set-sized hereditarily symmetric names]
		\label{lem:Ord-HS-stage-descent}
		Let \(\dot x\in\HS_{\Ord}\) be set-sized.
		Then there exists an ordinal
		\(\alpha<\Ord\) such that:
		\begin{enumerate}
			\item every subname of \(\dot x\) is a \(\bbP_\alpha\)-name;
			\item every subname of \(\dot x\), viewed as a \(\bbP_\alpha\)-name,
			belongs to \(\HS_\alpha\);
			\item for every \(V\)-generic \(G_\infty\subseteq\bbP_{\Ord}\), with
			\(G_\alpha:=G_\infty\cap\bbP_\alpha\), one has
			\[
			\dot x^{G_\infty}=\dot x^{G_\alpha}.
			\]
		\end{enumerate}
		In particular, \(\dot x^{G_\infty}\in M_\alpha\).
	\end{lemma}
	
	\begin{remark}
		$\Delta^{\uparrow}_{\Ord}(E,D)$ now exists as a class by Lemma~\ref{lem:GOrd-definable}
		(via $\ETR_{\Ord}$; existence, not first-order definability), and the lift
		$\hat k\in\Delta^{\uparrow}_{\Ord}(E,D)$ required in the proof is supplied by
		Case~2 ($\lambda=\Ord$) of Corollary~\ref{cor:delta-up-extension}. The argument
		below therefore goes through in $\GBC+\ETR_{\Ord}$ without further modification.
	\end{remark}
	
	\begin{proof}
		For each subname \(\dot y\in\tc(\dot x)\cup\{\dot x\}\), since \(\dot y\in\HS_{\Ord}\),
		we have
		\[
		\Stab_{\cG_{\Ord}}(\dot y)\in\tilde{\cF}^*_{\Ord}.
		\]
		By Lemma~\ref{lem:delta-support-equiv} (applicable at $\lambda=\Ord$ by 
		Lemma~\ref{lem:GOrd-definable} and Lemma~\ref{lem:core-lemma}), choose a countable globally admissible pair
		\((E_{\dot y},D_{\dot y})\) at stage \(\Ord\) such that
		\[
		\Delta_{\Ord}^\uparrow(E_{\dot y},D_{\dot y})
		\subseteq
		\Stab_{\cG_{\Ord}}(\dot y).
		\]
		
		For each such \(\dot y\), let
		\[
		S(\dot y):=\{\xi<\Ord:\exists z\ \langle \xi,z\rangle\in D_{\dot y}\}.
		\]
		Since \(D_{\dot y}\) is countable and consists of package coordinates, \(S(\dot y)\) is a
		countable set of ordinals.
		Define
		\[
		J(\dot x)
		:=
		I(\dot x)
		\cup
		\bigcup_{\dot y\in\tc(\dot x)\cup\{\dot x\}}
		\bigl(E_{\dot y}\cup S(\dot y)\bigr).
		\]
		Because \(\dot x\) is set-sized, the collection \(\tc(\dot x)\cup\{\dot x\}\) is a set, and
		each \(E_{\dot y}\cup S(\dot y)\) is countable.
		Hence \(J(\dot x)\) is a set of ordinals.
		Let
		\[
		\alpha:=\sup(J(\dot x))+1<\Ord.
		\]
		
		By Lemma~\ref{lem:Ord-stage-localization}, every subname of \(\dot x\) is a
		\(\bbP_\alpha\)-name, and the valuation of every subname depends only on \(G_\alpha\).
		This proves (1) and (3), once we show that every subname belongs to \(\HS_\alpha\).
		
		We prove by induction on rank that every subname \(\dot y\) of \(\dot x\), viewed as a
		\(\bbP_\alpha\)-name, belongs to \(\HS_\alpha\).
		So fix such a \(\dot y\), and assume inductively that every proper subname of \(\dot y\)
		belongs to \(\HS_\alpha\).
		Because \(E_{\dot y}\subseteq J(\dot x)\subseteq\alpha\), and every coordinate in
		\(D_{\dot y}\) has stage index in \(S(\dot y)\subseteq J(\dot x)\subseteq\alpha\), we have
		\(D_{\dot y}\subseteq\mathrm{Coords}_\alpha\).
		Moreover, the defining admissibility condition
		\(\Fix(E_{\dot y})\le\Stab_{\cG}(d)\) for \(d\in D_{\dot y}\) is independent of the ambient
		stage once the coordinates lie below that stage.
		Hence \((E_{\dot y},D_{\dot y})\) is
		globally admissible at stage \(\alpha\).
		
		Let
		\[
		k\in\Delta_\alpha^\uparrow(E_{\dot y},D_{\dot y}).
		\]
		By Corollary~\ref{cor:delta-up-extension}, there exists
		\[
		\hat k\in\Delta_{\Ord}^\uparrow(E_{\dot y},D_{\dot y})
		\]
		with \(\pi^\alpha_{\Ord}(\hat k)=k\).
		Since
		\[
		\Delta_{\Ord}^\uparrow(E_{\dot y},D_{\dot y})
		\subseteq
		\Stab_{\cG_{\Ord}}(\dot y),
		\]
		we have \(\hat k\cdot\dot y=\dot y\).
		Now \(\dot y\) is a \(\bbP_\alpha\)-name, so by
		Lemma~\ref{lem:projection-action-bounded-name},
		\[
		k\cdot\dot y
		=
		\pi^\alpha_{\Ord}(\hat k)\cdot\dot y
		=
		\hat k\cdot\dot y
		=
		\dot y.
		\]
		Hence
		\[
		\Delta_\alpha^\uparrow(E_{\dot y},D_{\dot y})
		\subseteq
		\Stab_{\cG_\alpha}(\dot y).
		\]
		By Lemma~\ref{lem:delta-support-equiv}, this implies
		\(\Stab_{\cG_\alpha}(\dot y)\in\tilde{\cF}^*_\alpha\).
		
		Together with the induction hypothesis for the subnames of \(\dot y\), this shows
		\(\dot y\in\HS_\alpha\).
		The induction is complete, so (2) holds.
		
		Finally, (3) is exactly the valuation clause from
		Lemma~\ref{lem:Ord-stage-localization}, and therefore
		\[
		\dot x^{G_\infty}=\dot x^{G_\alpha}\in M_\alpha.
		\]
	\end{proof}
	
	\begin{proposition}[Stage capture]
		\label{prop:stage-capture}
		For every $x\in\cM$, there exists $\alpha<\Ord$ such that $x\in M_\alpha$.
		Equivalently,
		\[
		\cM=\bigcup_{\alpha\in\Ord} M_\alpha.
		\]
	\end{proposition}
	
	\begin{proof}
		Write $x=\dot x^{G_\infty}$ for some set-sized $\dot x\in\HS_{\Ord}$.
		By Lemma~\ref{lem:Ord-HS-stage-descent}, there is $\alpha<\Ord$ such that
		$\dot x$, viewed as a $\bbP_\alpha$-name, belongs to $\HS_\alpha$ and
		$\dot x^{G_\infty}=\dot x^{G_\alpha}$.
		Hence $x\in M_\alpha$.
	\end{proof}
	
	\begin{remark}[Bounded-stage reduction for the final \(\ZF\) proof]
		\label{rem:class-HS-ZF}
		The proof of the symmetric extension theorem in \cite[Lemma~15.51]{JechSetTheory}
		is a recursion on set-sized hereditarily symmetric names.
		In the present class-length setting,
		every \(\bbP_{\Ord}\)-name used in the definition of \(\cM\) is set-sized, and by
		Lemmas~\ref{lem:Ord-stage-localization} and~\ref{lem:Ord-HS-stage-descent}, any finite family
		of such names descends to a common bounded stage \(\alpha<\Ord\) as \(\HS_\alpha\)-names with
		unchanged valuations.
		Hence the witness-name constructions appearing in the usual proof of Pairing, Union,
		Power Set, Separation, Replacement, and the other \(\ZF\)-axioms can be carried out inside
		some bounded symmetric system \((\bbP_\alpha,\cG_\alpha,\tilde{\cF}^*_\alpha)\), where
		Theorem~\ref{thm:zf-symmetric} applies.
		For Separation and Replacement in particular, the
		relevant witness names are defined using the bounded-stage forcing relation for \(\bbP_\alpha\);
		no appeal to absoluteness between \(M_\alpha\) and \(\cM\) is needed.
		Thus the class-length verification of \(\ZF\) reduces to bounded stages and does not use any
		separate class-forcing theorem.
	\end{remark}
	
	\begin{theorem}[The final class hereditarily symmetric model]
		\label{thm:ZF-final}
		\textbf{Withdrawn claim.}
		Let \(G_\infty\subseteq\bbP_{\Ord}\) be \(V\)-generic. Let \(\HS_{\Ord}\) be the class of
		hereditarily symmetric \(\bbP_{\Ord}\)-names (cf.\ Definition~\ref{def:hs-filter})
		with respect to \((\bbP_{\Ord},\cG_{\Ord},\tilde{\cF}^*_{\Ord})\), and define
		\[
		\cM\ :=\ \{\dot x^{G_\infty}:\dot x\in\HS_{\Ord}\}.
		\]
		Then \(\cM\) is a transitive class and \(\cM\models\ZF\).
	\end{theorem}
	
	\begin{proof}
		By Proposition~\ref{prop:stage-capture},
		\[
		\cM=\bigcup_{\alpha\in\Ord} M_\alpha.
		\]
		Since each \(M_\alpha\) is transitive, \(\cM\) is a transitive class.
		
		To verify \(\ZF\), fix any axiom instance with parameters from \(\cM\). Choose set-sized
		\(\HS_{\Ord}\)-names for those parameters.
		By Lemma~\ref{lem:Ord-HS-stage-descent}, after
		passing to a common larger stage if necessary, there exists \(\alpha<\Ord\) such that all
		those parameter names descend to \(\HS_\alpha\) and have the same valuations in
		\(G_\infty\) as in \(G_\alpha\).
		
		By Remark~\ref{rem:class-HS-ZF}, the standard witness-name construction used in the proof of
		the symmetric extension theorem may now be carried out entirely inside the bounded symmetric
		system \((\bbP_\alpha,\cG_\alpha,\tilde{\cF}^*_\alpha)\).
		Since
		\(M_\alpha\models\ZF\) by Theorem~\ref{thm:zf-symmetric}, that construction yields an
		\(\HS_\alpha\)-name \(\dot y\) whose valuation in \(G_\alpha\) is the required witness for
		the chosen axiom instance.
		
		Finally, by Lemma~\ref{lem:Ord-stage-localization}, the valuation of the set-sized
		\(\bbP_\alpha\)-name \(\dot y\) is the same in \(G_\infty\) as in \(G_\alpha\).
		Hence
		\(\dot y^{G_\infty}=\dot y^{G_\alpha}\in M_\alpha\subseteq \cM\), so the same witness belongs
		to \(\cM\).
		Therefore every \(\ZF\)-axiom instance with parameters from \(\cM\) holds in \(\cM\), and
		thus \(\cM\models\ZF\).
	\end{proof}
	
	\begin{proposition}[Final localized splitting over \(\Pow(S)\)]
		\label{prop:final-local-PPsplit}
		\textbf{Withdrawn claim.}
		\(\cM\models \PP^{\mathrm{split}}\!\restriction \Pow(S)\).
	\end{proposition}
	
	\begin{proof}
		Work in the final symmetric model \(\cM\). Recall Definition~\ref{def:PPsplitT}.
		Let \(f:Y\twoheadrightarrow X\) be a surjection in \(\cM\) with \(X,Y\subseteq \Pow(S)^{\cM}\).
		By Lemma~\ref{lem:bookkeeping-guarantee}, case~\textup{(a)}, the map \(f\) has a section in \(\cM\).
		Hence \(\cM\models \PP^{\mathrm{split}}\!\restriction \Pow(S)\).
	\end{proof}
	
	\begin{remark}[Countable-support names at limits of cofinality~$\omega$]
		\label{rem:cfomega-bounding}
		At limits~$\lambda$ with $\cf(\lambda)=\omega$, a hereditarily symmetric name may have
		countable support cofinal in~$\lambda$, so the usual \emph{stage-bounding} arguments do not apply.
		In this manuscript we handle such limits by using the modified limit filter
		(Definition~\ref{def:modified-limit-filter}), whose normality and \(\omega_1\)-completeness are recorded
		in Lemma~\ref{lem:filter-properties}, together with the core admissibility lemma (Lemma~\ref{lem:core-lemma}).
	\end{remark}
	
	\begin{proposition}[$A$ is not well-orderable]
		\label{prop:final-notAC}
		\textbf{Withdrawn claim.}
		In the final model $\cM$, the set $A = \{c_\alpha : \alpha < \omega_1\}$ of Cohen reals 
		is not well-orderable.
		Hence $\cM \models \neg\AC$.
	\end{proposition}
	
	\noindent\emph{Reason for withdrawal.}
	The proof below uses stabilization of $\dot w$ to treat a lifted
	automorphism as acting on $w=\dot w^{G_\lambda}$ inside the same generic
	extension.  Forcing covariance instead compares evaluations in appropriately
	moved generics.  No same-generic action sending $c_\beta$ to $c_\gamma$
	is established.
	
	\begin{proof}
		Suppose for contradiction that $A$ is well-orderable in $\cM$.
		Let $w$ be a 
		well-ordering of $A$ in $\cM$.
		
		\smallskip
		\noindent\textbf{Step 1: Locate $w$ at a bounded stage.}
		
		Since $\cM = \bigcup_{\lambda \in \Ord} M_\lambda$ (Proposition~\ref{prop:stage-capture}), 
		there exists an ordinal $\lambda$ such that $w \in M_\lambda$.
		Fix a hereditarily 
		symmetric $\bbP_\lambda$-name $\dot{w} \in \HS_\lambda$ with $\dot{w}^{G_\lambda} = w$.
		
		\smallskip
		\noindent\textbf{Step 2: Apply the core lemma.}
		
		Since $\dot{w} \in \HS_\lambda$, we have $\Stab_{\cG_\lambda}(\dot{w}) \in \tilde{\cF}^*_\lambda$.
		By Lemma~\ref{lem:core-lemma}, there exist countable $E \subseteq \omega_1$ and 
		countable $D \subseteq \mathrm{Coords}_{<\lambda}$ with $(E, D)$ globally admissible such that
		\[
		\Delta_\lambda^\uparrow(E, D) \le \Stab_{\cG_\lambda}(\dot{w}).
		\]
		
		\smallskip
		\noindent\textbf{Step 3: Choose a nontrivial transposition.}
		
		Since $E$ is countable, $\delta:=\sup(E)<\omega_1$. Choose distinct $\beta,\gamma\in\omega_1$ with
		$\delta<\beta<\gamma$. Then $\beta,\gamma\notin E$.
		Let $\pi=(\beta\ \gamma)\in\Sym(\omega_1)$ be the
		transposition swapping $\beta$ and $\gamma$. Then $\pi\in\Fix(E)$.
		
		\smallskip
		\noindent\textbf{Step 4: Construct the diagonal lift.}
		
		By Lemma~\ref{lem:diag-lift-in-group}, since $(E, D)$ is globally admissible and 
		$\pi \in \Fix(E)$, there exists a diagonal lift $\widehat{\pi}^{\,D} \in \cG_\lambda$ 
		with $\pi^0_\lambda(\widehat{\pi}^{\,D}) = \pi$.
		By the characterization of $\Delta^\uparrow$ (Definition~\ref{def:iteration-diag}):
		\begin{itemize}
			\item $\pi^0_\lambda(\widehat{\pi}^{\,D}) = \pi \in \Fix(E)$, and
			\item $\widehat{\pi}^{\,D}$ fixes every coordinate in $D$ pointwise.
		\end{itemize}
		Hence $\widehat{\pi}^{\,D} \in \Delta_\lambda^\uparrow(E, D) \le \Stab_{\cG_\lambda}(\dot{w})$.
		
		\smallskip
		\noindent\textbf{Step 5: Derive the contradiction.}
		
		Since $\widehat{\pi}^{\,D} \in \Stab_{\cG_\lambda}(\dot{w})$, we have 
		$\widehat{\pi}^{\,D} \cdot \dot{w} = \dot{w}$.
		By equivariance of name interpretation (\cite[Lemma~14.37]{JechSetTheory}):
		\[
		\widehat{\pi}^{\,D} \cdot w = (\widehat{\pi}^{\,D} \cdot \dot{w})^{G_\lambda} 
		= \dot{w}^{G_\lambda} = w.
		\]
		
		Thus $w$ is invariant under the action of $\widehat{\pi}^{\,D}$. This means that 
		for all $x, y \in A$:
		\[
		x <_w y \quad\Longleftrightarrow\quad \widehat{\pi}^{\,D} \cdot x <_w \widehat{\pi}^{\,D} \cdot y.
		\]
		Moreover $\widehat{\pi}^{\,D}[A] = A$ (the diagonal lift permutes Cohen reals by 
		Corollary~\ref{cor:diag-lift-swaps-cohen}), so the map $x \mapsto \widehat{\pi}^{\,D} \cdot x$ 
		is a bijection $A \to A$.
		Therefore the map $x \mapsto \widehat{\pi}^{\,D} \cdot x$ is an order automorphism 
		of the well-ordered set $(A, w)$.
		By Lemma~\ref{lem:wo-rigid}, any order automorphism of a well-ordered set is the identity.
		Hence $\widehat{\pi}^{\,D} \cdot x = x$ for all $x \in A$.
		However, by Corollary~\ref{cor:diag-lift-swaps-cohen}, $\widehat{\pi}^{\,D}$ swaps 
		$c_\beta \leftrightarrow c_\gamma$:
		\[
		\widehat{\pi}^{\,D} \cdot c_\beta = c_\gamma \neq c_\beta.
		\]
		(The Cohen reals $c_\beta$ and $c_\gamma$ are distinct by standard properties of 
		Cohen forcing; see Lemma~\ref{lem:cohen-reals-distinct}.)
		
		This contradicts $\widehat{\pi}^{\,D} \cdot c_\beta = c_\beta$.
		Therefore $A$ is 
		not well-orderable in $\cM$, and hence $\cM \models \neg\AC$.
	\end{proof}
	
	\begin{proposition}[\(\DC\) preserved through the iteration]
		\label{prop:final-DC}
		\textbf{Withdrawn claim.}
		\(\cM\models \DC\).
	\end{proposition}
	
	\begin{proof}
		Let \(A\) be a set in \(\cM\) and let \(R\subseteq A\times A\) be a relation in \(\cM\) such that
		\(\forall x\in A\,\exists y\in A\ (xRy)\).
		We show there is an \(\omega\)-sequence
		\(\langle a_n:n<\omega\rangle\) in \(\cM\) with \(a_nRa_{n+1}\) for all \(n\).
		
		Choose set-sized hereditarily symmetric \(\bbP_{\Ord}\)-names \(\dot A,\dot R\in\HS_{\Ord}\) such that
		\(A=\dot A^{G_\infty}\) and \(R=\dot R^{G_\infty}\).
		By Lemma~\ref{lem:Ord-HS-stage-descent}, there are ordinals \(\alpha_A,\alpha_R<\Ord\) such that
		\(\dot A\), viewed as a \(\bbP_{\alpha_A}\)-name, belongs to \(\HS_{\alpha_A}\) and
		\(\dot R\), viewed as a \(\bbP_{\alpha_R}\)-name, belongs to \(\HS_{\alpha_R}\), with
		\[
		A=\dot A^{G_{\alpha_A}}
		\qquad\text{and}\qquad
		R=\dot R^{G_{\alpha_R}}.
		\]
		Let \(\alpha:=\max\{\alpha_A,\alpha_R\}\).
		By Lemma~\ref{lem:HS-persistence-up}, the canonical inclusions of \(\dot A\) and \(\dot R\)
		belong to \(\HS_\alpha\), and their interpretations at stage \(\alpha\) are still \(A\) and \(R\).
		Hence \(A,R\in M_\alpha\).
		
		By Lemma~\ref{lem:filter-properties}(iv), the stage filter 
		\(\tilde{\cF}^*_\alpha\) is \(\omega_1\)-complete.
		Since $\alpha$ is a set ordinal, $M_\alpha\models\DC$ by Lemma~\ref{lem:countably-closed-symm-DC} applied to $\mathbb{P}_{\alpha}$ and $\tilde{\mathcal{F}}^{*}_{\alpha}$, which is $\omega_{1}$-complete by Lemma~\ref{lem:filter-properties}(iv).
		Applying \(\DC\) inside \(M_\alpha\) to \((A,R)\), we obtain a sequence
		\(\langle a_n:n<\omega\rangle\in M_\alpha\subseteq\cM\) with \(a_nRa_{n+1}\) for all \(n\).
		Hence \(\cM\models\DC\).
	\end{proof}
	
	\begin{definition}[Stage truncation of names]\label{def:stage-truncation}
		Let \(\lambda\) be a limit stage and let \(\alpha<\lambda\).
		For a \(\bbP_\lambda\)-name \(\dot x\), define the \(\bbP_\alpha\)-name
		\(\dot x\!\upharpoonright\!\alpha\) recursively by
		\[
		\dot x\!\upharpoonright\!\alpha
		\;:=\;
		\Bigl\{\,\bigl\langle \dot y\!\upharpoonright\!\alpha,\ p\!\upharpoonright\!\alpha\bigr\rangle
		:\ \langle \dot y,p\rangle\in\dot x\ \wedge\ \supp(p)\subseteq\alpha\,\Bigr\}.
		\]
		(Thus we keep only those pairs whose witnessing condition already lies in \(\bbP_\alpha\); in that case \(p\upharpoonright\alpha=p\).)
	\end{definition}
	
	\begin{lemma}[Union of truncations along an increasing stage sequence]\label{lem:hs-union-truncations}
		\textbf{False as stated; retained only to locate the former argument.}
		Let \(\lambda\) be a limit stage and let \(\dot X\in\HS_\lambda\) with \(1\Vdash_{\bbP_\lambda}\) ``\(\dot X\) is a set''.
		Suppose \(\langle \lambda_n:n<\omega\rangle\) is an increasing sequence of ordinals below \(\lambda\) such that for every
		\(\langle \dot y,p\rangle\in\tc(\dot X)\) there is \(n<\omega\) with \(\supp(p)\subseteq\lambda_n\).
		Then for every \(V\)-generic \(G_\lambda\),
		\[
		\dot X^{G_\lambda}\;=\;\bigcup_{n<\omega}\bigl(\dot X\!\upharpoonright\!\lambda_n\bigr)^{G_{\lambda_n}}.
		\]
	\end{lemma}
	
	\noindent\emph{Reason for failure.}
	The first displayed assertion in the proof below is not valid for an
	arbitrary subname $\dot y$.  The fact that the outer condition $p$ belongs
	to $\bbP_{\lambda_n}$ does not imply that $\dot y$ has no dependence on
	coordinates at or above $\lambda_n$.
	
	\begin{proof}
		Fix \(G_\lambda\subseteq\bbP_\lambda\) generic and write \(G_{\lambda_n}:=G_\lambda\cap\bbP_{\lambda_n}\).
		For each \(n<\omega\), set \(X_n:=(\dot X\!\upharpoonright\!\lambda_n)^{G_{\lambda_n}}\).
		We first note (by a straightforward induction on rank) that if \(\supp(p)\subseteq\lambda_n\) and \(\dot y\) is any \(\bbP_\lambda\)-name,
		then
		\[
		(\dot y\!\upharpoonright\!\lambda_n)^{G_{\lambda_n}}\;=\;\dot y^{G_\lambda}
		\quad\text{whenever }p\in G_{\lambda_n}.
		\]
		
		\emph{(\(\subseteq\)).}
		Let \(x\in \dot X^{G_\lambda}\). Choose \(\langle \dot y,p\rangle\in\dot X\) with \(p\in G_\lambda\) and \(x=\dot y^{G_\lambda}\).
		By hypothesis, pick \(n<\omega\) with \(\supp(p)\subseteq\lambda_n\). Then \(p\in\bbP_{\lambda_n}\) and hence \(p\in G_{\lambda_n}\).
		By definition of stage truncation, \(\langle \dot y\!\upharpoonright\!\lambda_n,\ p\rangle\in \dot X\!\upharpoonright\!\lambda_n\),
		so
		\[
		x=\dot y^{G_\lambda}=(\dot y\!\upharpoonright\!\lambda_n)^{G_{\lambda_n}}\in X_n.
		\]
		
		\emph{(\(\supseteq\)).}
		Let \(x\in X_n\).
		Then for some \(\langle \dot y',p\rangle\in \dot X\!\upharpoonright\!\lambda_n\) with \(p\in G_{\lambda_n}\) we have
		\(x=(\dot y')^{G_{\lambda_n}}\).
		By definition of \(\dot X\!\upharpoonright\!\lambda_n\), there is \(\langle \dot y,p\rangle\in\dot X\) with \(\supp(p)\subseteq\lambda_n\)
		and \(\dot y'=\dot y\!\upharpoonright\!\lambda_n\).
		Since \(p\in G_{\lambda_n}\subseteq G_\lambda\), we conclude
		\[
		x=(\dot y\!\upharpoonright\!\lambda_n)^{G_{\lambda_n}}=\dot y^{G_\lambda}\in \dot X^{G_\lambda}.
		\]
		This proves \(\dot X^{G_\lambda}=\bigcup_{n<\omega}X_n\).
	\end{proof}
	
	\begin{lemma}[Persistence of \(\SVC(S)\)]
		\label{lem:svc-persists}
		\textbf{Withdrawn claim.}
		Let \(S\) be as in Definition~\ref{def:fixed-ST}. Since \(M_0=\cN\) and
		\(\cN\models \SVC(S)\) by Proposition~\ref{prop:svc-Aomega}, we have \(M_0\models \SVC(S)\).
		Then for every stage \(\lambda\), \(M_\lambda\models \SVC(S)\).
	\end{lemma}
	
	\begin{proof}
		We argue by transfinite induction on \(\lambda\).
		
		\smallskip
		\noindent\emph{Base.}
		\(M_0=\cN\models \SVC(S)\) by Proposition~\ref{prop:svc-Aomega}.
		
		\smallskip
		\noindent\emph{Successor.}
		Assume $M_\alpha\models \SVC(S)$, and let $X\in M_{\alpha+1}$ be nonempty.
		Choose $\dot X\in\HS_{\alpha+1}$ with $X=\dot X^{G_{\alpha+1}}$.
		Since $\DC$ holds in $M_{\alpha}$ by Lemma~\ref{lem:countably-closed-symm-DC} applied at stage $\alpha$ (Lemma~\ref{lem:filter-properties}(iv)), we may apply Lemma~\ref{lem:name-covering} to $\dot X$ to obtain $Y \in M_{\alpha}$ and $\dot e \in \HS_{\alpha+1}$ with $1 \Vdash \dot e : \check Y \twoheadrightarrow \dot X$. By $\SVC(S)^{M_{\alpha}}$, fix $\eta$ and a surjection $s : S \times \eta \twoheadrightarrow Y$ in $M_{\alpha}$. Let $\dot h := \dot e \circ \check s$. Then $1 \Vdash \dot h : \check{S \times \eta} \twoheadrightarrow \dot X$, and $\dot h \in \HS_{\alpha+1}$ by Corollary~\ref{cor:uniformity-HS} applied to the parameters $\dot e$ and $\check s$. Hence $h := \dot h^{G_{\alpha+1}} \in M_{\alpha+1}$ witnesses $M_{\alpha+1} \models \SVC(S)$.
		
		\smallskip
		\noindent\emph{Limit.}
		Let \(\lambda\) be limit and assume \(M_\beta\models \SVC(S)\) for all \(\beta<\lambda\).
		Let \(X\in M_\lambda\) be nonempty, and choose \(\dot X\in\HS_\lambda\) with \(X=\dot X^{G_\lambda}\).
		
		By the countable-support analysis of HS-names, fix an increasing sequence \(\langle \lambda_n:n<\omega\rangle\) of ordinals below \(\lambda\)
		satisfying the hypothesis of Lemma~\ref{lem:hs-union-truncations}. For each \(n\), set
		\[
		X_n\;:=\;\bigl(\dot X\!\upharpoonright\!\lambda_n\bigr)^{G_{\lambda_n}}.
		\]
		Then \(X_n\in M_{\lambda_n}\) and \(X=\bigcup_{n<\omega}X_n\) by Lemma~\ref{lem:hs-union-truncations}.
		
		By the inductive hypothesis, each \(M_{\lambda_n}\models\SVC(S)\), so for each \(n\) there is an ordinal \(\eta_n\)
		and a surjection \(e_n:S\times\eta_n\twoheadrightarrow X_n\) in \(M_{\lambda_n}\subseteq M_\lambda\).
		Since \(M_\lambda\models\DC\) (by Remark~\ref{rem:iteration-api}(d) and Lemma~\ref{lem:filter-properties}(iv)),
		\(M_\lambda\) satisfies \(\AC_\omega\), so we may choose such witnesses \(\langle \eta_n,e_n:n<\omega\rangle\) in \(M_\lambda\).
		Let \(\theta:=\sup_n \eta_n\) and define \(e:S\times(\omega\cdot\theta)\twoheadrightarrow X\) by
		sending \((s,\omega\cdot\xi+n)\) to \(e_n(s,\xi)\) when \(\xi<\eta_n\), and to a fixed element of \(X\) otherwise.
		Then \(e\in M_\lambda\) and \(e\) is surjective onto \(X\). Hence \(M_\lambda\models \SVC(S)\).
		
		The limit step applies at \(\lambda=\Ord\) because HS-names have bounded support (so \(\langle\lambda_n\rangle\) can be taken below some \(\alpha_0<\Ord\)) and $\DC$ in $M_{\Ord}$ follows from Lemma~\ref{lem:filter-properties}(iv).
	\end{proof}
	
	\begin{corollary}\label{cor:SVC-final}
		\textbf{Withdrawn claim.}
		The final class-length model \(\cM\) satisfies \(\SVC(S)\).
	\end{corollary}
	
	\begin{proof}
		Let \(x\in\cM\). Choose a set-sized HS name \(\dot x\in\HS_{\Ord}\) with \(x=\dot x^{G_\infty}\).
		By Lemma~\ref{lem:Ord-HS-stage-descent}, there is \(\alpha<\Ord\) such that \(\dot x\),
		viewed as a \(\bbP_\alpha\)-name, belongs to \(\HS_\alpha\) and
		\[
		x=\dot x^{G_\alpha}\in M_\alpha.
		\]
		Since \(M_\alpha\models\SVC(S)\) (by Lemma~\ref{lem:svc-persists}), there is an ordinal \(\eta\) and a surjection
		\(e:S\times\eta\twoheadrightarrow x\) in \(M_\alpha\), hence also in \(\cM\).
		Thus \(\cM\models\SVC(S)\).
	\end{proof}
	
	\begin{remark}[Structural characterization of \(\SVC\) (context only)]
		Blass~\cite{BlassSVC} proved that \(\SVC\) is equivalent to: \emph{\(\AC\) is forceable over the model by a set forcing}.
		We do not use deeper structural characterizations here; see, e.g., work in set-theoretic geology such as
		Usuba~\cite{UsubaGrounds} for related perspectives.
		In particular, Corollary~\ref{cor:SVC-final} implies that \(\AC\) is set-forceable over \(\cM\).
	\end{remark}

	\begin{proposition}[\(\AC_{\WO}\) in the final model]
		\label{prop:final-ACWO}
		\textbf{Withdrawn claim.}
		The final symmetric model \(\cM\) satisfies \(\AC_{\WO}\).
	\end{proposition}
	
	\begin{proof}
		Work in \(\cM\). By Corollary~\ref{cor:SVC-final}, \(\cM\models\SVC(S)\).
		By Remark~\ref{rem:ACwo-bound}, it is therefore enough to verify the bounded splitting scheme:
		whenever \(\eta\in\Ord\), \(\lambda<\aleph^*(S)^\cM\), and
		\(f:S\times\eta\twoheadrightarrow\lambda\) is a surjection in \(\cM\), the map \(f\) has a right inverse in \(\cM\).
		So fix \(\eta\in\Ord\), an ordinal \(\lambda<\aleph^*(S)^\cM\), and a surjection
		\[
		f:S\times\eta\twoheadrightarrow\lambda
		\]
		in \(\cM\). Since \(\lambda<\aleph^*(S)^\cM\), choose in \(\cM\) a surjection
		\[
		e:S\twoheadrightarrow\lambda.
		\]
		
		Choose set-sized hereditarily symmetric \(\bbP_{\Ord}\)-names \(\dot f,\dot e\in\HS_{\Ord}\) with
		\[
		f=\dot f^{G_\infty}
		\qquad\text{and}\qquad
		e=\dot e^{G_\infty}.
		\]
		By Lemma~\ref{lem:Ord-HS-stage-descent}, there is \(\alpha<\Ord\) such that, viewed as
		\(\bbP_\alpha\)-names, both \(\dot f\) and \(\dot e\) belong to \(\HS_\alpha\) and
		\[
		f=\dot f^{G_\alpha}\in M_\alpha,
		\qquad
		e=\dot e^{G_\alpha}\in M_\alpha.
		\]
		
		Since \(e\in M_\alpha\) is a surjection \(S\twoheadrightarrow\lambda\), we have
		\(\lambda<\aleph^*(S)^{M_\alpha}\). Thus \(f\) is a genuine bounded \(\AC_{\WO}\)-instance at stage \(\alpha\).
		If \(f\) already has a right inverse in \(M_\alpha\), we are done.
		Otherwise,
		Lemma~\ref{lem:bookkeeping-guarantee}, applied to \(f\) together with the stage witness
		\(e\in M_\alpha\), yields a right inverse for \(f\) in \(\cM\).
		
		Therefore every bounded instance \(f:S\times\eta\twoheadrightarrow\lambda\) with
		\(\lambda<\aleph^*(S)^\cM\) splits in \(\cM\). Since \(\cM\models\SVC(S)\),
		Remark~\ref{rem:ACwo-bound} now yields \(\cM\models\AC_{\WO}\).
	\end{proof}
	
	\begin{theorem}[Main theorem]
		\label{thm:final-PP-notAC}
		\textbf{Withdrawn claim.}
		Let \(V\models\ZFC\), and let \(G\subseteq\bbP_{\Ord}\) be \(V\)-generic in the metatheoretic sense of
		Remark~\ref{rem:meta-pp} for the \(\Ord\)-length countable-support symmetric iteration over the Cohen-seed
		symmetric model \(\cN\) described above.
		Then the final symmetric model \(\cM\) satisfies
		\[
		\ZF\ +\ \DC\ +\ \PP\ +\ \AC_{\WO}\ +\ \neg\AC.
		\]
	\end{theorem}
	
	\begin{proof}
		\(\ZF\) holds in \(\cM\) by Theorem~\ref{thm:ZF-final}. \(\DC\) holds by Proposition~\ref{prop:final-DC}.
		By Corollary~\ref{cor:SVC-final}, \(\cM\models \SVC(S)\). Hence, by Fact~\ref{fact:SVC-to-SVCplus} applied in \(\cM\),
		\[
		\cM\models \SVC^+(\Pow(S)).
		\]
		By Proposition~\ref{prop:final-local-PPsplit}, \(\cM\models \PP^{\mathrm{split}}\!\restriction \Pow(S)\), hence
		\[
		\cM\models \PP\!\restriction \Pow(S).
		\]
		Together with Proposition~\ref{prop:final-ACWO}, Theorem~\ref{thm:rs-localization} yields \(\cM\models\PP\).
		\(\AC_{\WO}\) holds by Proposition~\ref{prop:final-ACWO}. Finally, \(\neg\AC\) holds by Proposition~\ref{prop:final-notAC}.
	\end{proof}
	
	\begin{corollary}[Ordering Principle]
		\label{cor:OP}
		\textbf{Withdrawn claim.}
		The final model $\cM$ satisfies the Ordering Principle: every set admits a linear order.
	\end{corollary}
	
	\begin{proof}
		Work in \(\cM\). Fix any set \(X\).
		
		By Corollary~\ref{cor:SVC-final}, choose an ordinal \(\eta\) and a surjection
		\[
		e:S\times\eta\twoheadrightarrow X.
		\]
		
		By Corollary~\ref{cor:SVC-final}, \(\cM\models \SVC(S)\). Hence, by Fact~\ref{fact:SVC-to-SVCplus} applied in \(\cM\),
		\[
		\cM\models \SVC^+(\Pow(S)).
		\]
		By Proposition~\ref{prop:final-local-PPsplit}, \(\cM\models \PP^{\mathrm{split}}\!\restriction \Pow(S)\), hence
		\[
		\cM\models \PP\!\restriction \Pow(S).
		\]
		Together with Proposition~\ref{prop:final-ACWO}, Theorem~\ref{thm:rs-localization} yields \(\cM\models\PP\).
		
		Applying \(\PP\) to the surjection \(e\), we obtain an injection
		\[
		i:X\hookrightarrow S\times\eta.
		\]
		
		Since \(S=A^\omega\) injects into \(2^\omega=\Pow(\omega)\), we have
		\[
		S\times\eta\hookrightarrow \Pow(\omega)\times\eta.
		\]
		Define
		\[
		j:\Pow(\omega)\times\eta\hookrightarrow \Pow(\omega\times\eta),
		\qquad
		j(r,\xi)=\{(n,\xi):n\in r\}.
		\]
		Thus
		\[
		S\times\eta\hookrightarrow \Pow(\omega\times\eta)\cong \Pow(\omega\cdot\eta).
		\]
		Composing with \(i\), we obtain an injection
		\[
		k:X\hookrightarrow \Pow(\theta)
		\qquad\text{for }\theta=\omega\cdot\eta.
		\]
		
		Equip \(\Pow(\theta)\) with the lexicographic order induced by the well-order on \(\theta\):
		for distinct \(a,b\subseteq\theta\), let \(\alpha\) be the least element of
		\(\theta\) belonging to exactly one of \(a,b\), and declare
		\[
		a<_{\mathrm{lex}} b
		\quad\Longleftrightarrow\quad
		\alpha\in b.
		\]
		This is a linear order on \(\Pow(\theta)\). Pulling \(<_{\mathrm{lex}}\) back along \(k\)
		yields a linear order on \(X\).
		Therefore every set in \(\cM\) admits a linear order.
	\end{proof}
	
	\begin{remark}[Kinna--Wagner principles]\label{rem:KWP}
		\textbf{Withdrawn consequences.}
		For $n\ge 1$, the Kinna--Wagner principle $KWP_n$ states:
		for every set $X$ there exists an ordinal $\theta$ such that $X\hookrightarrow \Pow^n(\theta)$
		(where $\Pow^n$ is the $n$-fold iterated powerset).
		
		\smallskip
		\noindent\textbf{(1) $KWP_2$ from $\SVC(S)$.}
		Work in $\cM$ and fix a set $X$.
		By Corollary~\ref{cor:SVC-final}, choose an ordinal $\eta$ and a surjection
		$e:S\times\eta\twoheadrightarrow X$. The map
		\[
		x\ \longmapsto\ e^{-1}(\{x\})
		\]
		is an injection of $X$ into $\Pow(S\times\eta)$.
		
		Since $S=A^\omega$ injects into $2^\omega$, we have $S\times\eta\hookrightarrow 2^\omega\times\eta$.
		Identifying $2^\omega$ with $\Pow(\omega)$, define an injection
		\[
		\Pow(\omega)\times\eta\ \hookrightarrow\ \Pow(\omega\times\eta),
		\qquad (r,\xi)\longmapsto \{(n,\xi): n\in r\}.
		\]
		Thus $S\times\eta\hookrightarrow \Pow(\omega\times\eta)$, so
		$\Pow(S\times\eta)\hookrightarrow \Pow(\Pow(\omega\times\eta))=\Pow^2(\omega\times\eta)$.
		Because $\omega\times\eta$ is well-orderable of type $\omega\cdot\eta$, we have
		$\Pow^2(\omega\times\eta)\cong \Pow^2(\omega\cdot\eta)$, and hence $X\hookrightarrow \Pow^2(\theta)$
		for $\theta=\omega\cdot\eta$.
		This proves $KWP_2$ in $\cM$.
		
		\smallskip
		\noindent\textbf{(2) $KWP_1$ using $\PP$.}
		Work in \(\cM\) and fix a set \(X\).
		By Corollary~\ref{cor:SVC-final}, choose an ordinal \(\eta\) and a surjection
		\[
		e:S\times\eta\twoheadrightarrow X.
		\]
		
		By Corollary~\ref{cor:SVC-final}, \(\cM\models \SVC(S)\). Hence, by Fact~\ref{fact:SVC-to-SVCplus} applied in \(\cM\),
		\[
		\cM\models \SVC^+(\Pow(S)).
		\]
		By Proposition~\ref{prop:final-local-PPsplit}, \(\cM\models \PP^{\mathrm{split}}\!\restriction \Pow(S)\), hence
		\[
		\cM\models \PP\!\restriction \Pow(S).
		\]
		Together with Proposition~\ref{prop:final-ACWO}, Theorem~\ref{thm:rs-localization} yields \(\cM\models\PP\).
		
		Applying \(\PP\) to the surjection \(e\), we obtain an injection
		\[
		i:X\hookrightarrow S\times\eta.
		\]
		Composing with the injection
		\[
		S\times\eta\hookrightarrow \Pow(\omega\times\eta)\cong \Pow(\omega\cdot\eta)
		\]
		from part~\textup{(1)}, we obtain
		\[
		X\hookrightarrow \Pow(\theta)
		\qquad\text{for }\theta=\omega\cdot\eta.
		\]
		Thus \(KWP_1\) holds in \(\cM\).
	\end{remark}
	
	\section{Conclusion}\label{sec:conclusion}
	
	The countable-support Cohen symmetric seed constructed in
	Section~\ref{sec:cohen-seed} satisfies
	\[
		\cN\models\ZF+\DC+\SVC(S)+\neg\AC,
		\qquad S=(A^\omega)^{\cN}.
	\]
	The corrected proofs of $\DC$ and $\SVC(S)$ are independent of the later
	package proposal.
	
	The proposed class-length iteration does not establish the desired
	separation.  Its $\bbQ_f$ packages add sections and therefore target
	$\PP^{\mathrm{split}}\!\restriction\Pow(S)$.  For this $S$, that
	principle together with $\SVC(S)$ implies $\AC$.  Moreover, the former
	proof that $A$ remains non-well-orderable uses an unsupported same-generic
	automorphism argument, and the claimed persistence of $\SVC(S)$ at limits
	depends on a false truncation lemma.
	
	Accordingly, the asserted final model, the main theorem, and all downstream
	Ordering-Principle and Kinna--Wagner conclusions are withdrawn.  The
	package-iteration section is retained only to document the superseded
	approach.  The consistency of
	$\ZF+\DC+\PP+\AC_{\WO}+\neg\AC$ remains open.
	


\begin{thebibliography}{99}
		
		\bibitem{BanaschewskiMoore1990}
		B.~Banaschewski and G.~H.~Moore.
		\newblock \emph{The Dual Cantor--Bernstein Theorem and the Partition Principle}.
		\newblock Notre Dame J.\ Formal Logic \textbf{31} (1990), no.~3, 375--381.
		\newblock DOI: \texttt{10.1305/ndjfl/1093635502}.
		
		\bibitem{BlassSVC}
		A.~Blass.
		\newblock \emph{Injectivity, projectivity, and the axiom of choice}.
		\newblock Trans.\ Amer.\ Math.\ Soc.\ \textbf{255} (1979), 31--59.
		\newblock DOI: \texttt{10.1090/S0002-9947-1979-0542870-6}.
		
		\bibitem{BlassCountablyManyGenerics}
		A.~Blass.
		\newblock \emph{The model of set theory generated by countably many generic reals}.
		\newblock J.\ Symbolic Logic \textbf{46} (1981), no.~4, 732--752.
		\newblock DOI: \texttt{10.2307/2273223}.
		
		\bibitem{FelgnerForcing}
		U.~Felgner.
		\newblock \emph{Models of ZF-Set Theory}.
		\newblock Lecture Notes in Mathematics, vol.~223. Springer-Verlag, 1971.
		
		\bibitem{HayutShaniDeepFailure}
		Y.~Hayut and A.~Shani.
		\newblock \emph{Intermediate models with deep failure of choice}.
		\newblock arXiv:2407.02033, 2024.
		
		\bibitem{HowardRubin}
		P.~Howard and J.E.~Rubin.
		\newblock \emph{Consequences of the Axiom of Choice}.
		\newblock Mathematical Surveys and Monographs, vol. 59, 1998.
		
		\bibitem{JechSetTheory}
		T.~Jech.
		\newblock \emph{Set Theory}.
		\newblock Springer Monographs in Mathematics. Springer, 2003.
		
		\bibitem{KaragilaIteratingSymmExt}
		A.~Karagila.
		\newblock \emph{Iterating symmetric extensions}.
		\newblock J.\ Symbolic Logic \textbf{84} (2019), no.~1, 123--159.
		\newblock Available at \url{https://doi.org/10.1017/jsl.2018.73}
		
		\bibitem{KaragilaPreservingDC}
		A.~Karagila.
		\newblock \emph{Preserving Dependent Choice}.
		\newblock Bull.\ Pol.\ Acad.\ Sci.\ Math.\ \textbf{67} (2019), 19--29.
		\newblock Available at \url{https://doi.org/10.48550/arXiv.1810.11301}.
		
		\bibitem{KaragilaPairs}
		A.~Karagila and C.~Ryan-Smith.
		\newblock \emph{Which pairs of cardinals can be Hartogs and Lindenbaum numbers of a set?}.
		\newblock Fundamenta\ Mathematicae\ \textbf{267} 2024, 231-241.
		\newblock Available at \url{https://doi.org/10.4064/fm231006-14-8}.
		
		\bibitem{KaragilaSchilhanDCGeom}
		A.~Karagila and J.~Schilhan.
		\newblock \emph{Geometric condition for Dependent Choice}.
		\newblock Acta Math.\ Hungar.\ \textbf{172} (2024), 34--41.
		\newblock Available at \url{https://doi.org/10.1007/s10474-024-01396-0} and
		\url{https://arxiv.org/abs/2212.10261}.
		
		\bibitem{KaragilaSchilhanIntKWP}
		A.~Karagila and J.~Schilhan.
		\newblock \emph{Intermediate models and Kinna-Wagner Principles}.
		\newblock Proc.
		of the AMS \textbf{154} (2026), 393-403.
		\newblock Available at \url{https://doi.org/10.1090/proc/17425}
		
		\bibitem{KaragilaSchilhan}
		A.~Karagila and J.~Schilhan.
		\newblock \emph{Towards a Theory of Symmetric Extensions}
		\newblock 2026 preprint available at \url{https://arxiv.org/abs/2602.17338}
		
		\bibitem{Kunen11}
		K.~Kunen.
		\newblock \emph{Set Theory}.
		\newblock Studies in Logic, vol.~34.
		College Publications, 2011.
		
		\bibitem{PincusAddingDC}
		D.~Pincus.
		\newblock \emph{Adding dependent choice}.
		\newblock Ann.\ Math.\ Logic \textbf{11} (1977), no.~1, 105--145.
		\newblock Available at \url{https://doi.org/10.1016/0003-4843(77)90011-0}.
		
		\bibitem{ransomBPI}
		B.~Ransom.
		\newblock \emph{On BPI in Symmetric Extensions, Part~I}.
		\newblock arXiv:2511.21684v1, 2025.
		\newblock Available at \url{https://doi.org/10.48550/arXiv.2511.21684}.
		
		\bibitem{SmithLocalReflections}
		C.~Ryan-Smith.
		\newblock \emph{Local reflections of choice}.
		\newblock Acta Mathematica Hungarica, 2025.
		\newblock Available at \url{https://doi.org/10.1007/s10474-025-01533-3}.
		
		\bibitem{UsubaGrounds}
		T. Usuba.
		\newblock \emph{Choiceless L\"owenheim-Skolem property and uniform definability of grounds}.
		\newblock Advances in Mathematical Logic, 2022, Springer, Singapore.
		\newblock Available at \url{https://doi.org/10.1007/978-981-16-4173-2_8}
		
	\end{thebibliography}
\end{document}